\documentclass[11pt]{amsart}
\usepackage[letterpaper,margin=1in]{geometry}
\usepackage{amsmath,amssymb,amsthm,mathtools}
\usepackage{booktabs}
\usepackage{tikz}
\usetikzlibrary{positioning,arrows.meta,fit,calc}
\usepackage{cite}
\usepackage[colorlinks,linkcolor=blue,citecolor=blue,urlcolor=blue]{hyperref}

\newtheorem{theorem}{Theorem}[section]
\newtheorem{lemma}[theorem]{Lemma}
\newtheorem{proposition}[theorem]{Proposition}
\newtheorem{corollary}[theorem]{Corollary}
\theoremstyle{definition}
\newtheorem{definition}[theorem]{Definition}
\newtheorem{examplex}[theorem]{Example}
\newenvironment{example}{\pushQED{\qed}\examplex}{\popQED\endexamplex}
\theoremstyle{remark}
\newtheorem{remark}[theorem]{Remark}
\numberwithin{equation}{section}
\makeatletter
\def\paragraph{\@startsection{paragraph}{4}\z@{.3\baselineskip\@plus.1\baselineskip}{-\fontdimen2\font}{\normalfont\itshape}}
\makeatother

\newcommand{\R}{\mathbb{R}}
\newcommand{\GF}{\mathrm{GF}}
\newcommand{\PG}{\mathrm{PG}}
\newcommand{\AG}{\mathrm{AG}}
\newcommand{\aff}{\operatorname{aff}}
\newcommand{\conv}{\operatorname{conv}}
\newcommand{\rk}{\operatorname{rk}}
\makeatletter
\renewcommand*{\top}{%
  {\mathpalette\@transpose{}}%
}
\newcommand*{\@transpose}[2]{%
  \raisebox{\depth}{$\m@th#1\scriptscriptstyle\mathsf{T}$}%
}
\makeatother

\title[Dimension of assignment and transversal polytopes]{The dimension of multi-index assignment polytopes and\\ \mbox{\spaceskip=0.8\fontdimen2\font Latin-square transversal polytopes, via association schemes}}
\author{Jon Lee}
\address{University of Michigan, Ann Arbor, MI, USA}
\email{jonxlee@umich.edu}
\date{\today}
\subjclass[2020]{Primary 90C27, 52B12; Secondary 05B15, 05E30}
\thanks{This work was supported in part by U.S. Office of Naval Research grant N00014-21-1-2135.}
\keywords{integer programming, polytope dimension, stable-set polytope, multi-index assignment, Latin-square transversal, association scheme}

\begin{document}
\begin{abstract}
We determine the dimensions, and the implicit equations, of several polytopes from integer programming whose vertices are the largest stable sets of highly symmetric graphs. Our main tools are two new dimension formulae within Lee's spectral method, which reads the dimension off of an association scheme: for amorphic schemes it becomes a threshold formula, and for point sets in $\GF(q)^d$ invariant under the maps $z\mapsto az+b$ it is $q-1$ times the number of hyperplane directions along which some point is not equidistributed. A single integer point, a witness, then certifies a whole eigenspace. For the multi-index assignment polytopes of Appa, Magos, and Mourtos at prime-power orders, the dimension equals that of the linear relaxation except in two infinite families, which are, to our knowledge, the first counterexamples to their expectation that the two always coincide. The missing equations are explicit clique equations that the standard integer-programming formulation omits. For every abelian group with a transversal, the Latin-square transversal polytope of its Cayley table, equivalently the polytope of rainbow perfect matchings of the corresponding edge-colored $K_{n,n}$\,, has the maximal dimension ${(n-1)(n-2)}$; the witnesses come from new constructions of complete mappings. We also treat transversals of Desarguesian nets and the largest sets of symmetric $3\times3$ and $4\times4$ matrices over $\GF(q)$ whose differences have rank at least~$3$, of sizes $q^3$ and $q^5$.
\end{abstract}
\maketitle
\tolerance=300\relax\emergencystretch=3em\relax
\setlength{\parfillskip}{0pt plus 0.45\textwidth}

\section{Introduction}\label{sec:intro}

\paragraph{Dimension in integer linear programming.}
Let $E$ be a finite set and let $X\subseteq\R^E$ be the incidence vectors of some family of combinatorial objects on the ground set $E$: perfect matchings, Hamiltonian tours, orthogonal arrays, transversals of a Latin square. The dimension of the polytope $\conv(X)$ is the first thing one needs in order to optimize over $X$ by polyhedral methods: it counts the independent linear equations that hold on all of $X$ (equations an integer-programming formulation can and should contain), and the facets, which carry the cutting planes, are the faces of dimension one less. When $\dim(\conv(X))$ is less than the dimension of the natural linear relaxation, the formulation is missing equations, and adding them cuts off fractional points.

A theme of this paper is that many such families are the \emph{largest stable sets} of a graph $\Gamma$ on $E$: a transversal of a Latin square is a stable set of size $n$ in the graph in which two cells are adjacent when they share a row, a column, or a symbol; an orthogonal array of strength $s$ and index one is a code of maximum size for its minimum distance; a set of symmetric matrices with pairwise differences of rank at least $d$ is a clique in a rank-distance graph. The convex hull of the largest stable sets is the face of the stable-set polytope of $\Gamma$ on which the rank inequality $x(E)\le\alpha(\Gamma)$ is tight, and optimizing a weight over it is a maximum-weight problem over the largest stable sets. In the examples of this paper, $\alpha(\Gamma)$ is a bound of Delsarte or Hoffman type coming from an association scheme, the face is far from being a facet, and we ask which equations hold on it.

\paragraph{The spectral method.}
Lee~\cite{L1,L2,L3} introduced a method that answers this question when $X$ has enough symmetry. With $\bar x$ the average of $X$, the dimension is the rank of the symmetric positive semidefinite matrix $B\coloneqq\sum_{x\in X}(x-\bar x)(x-\bar x)^\top$, which is indexed by the (small) ground set $E$ rather than by the (large) family $X$. If $X$ is \emph{concordant} with an association scheme on $E$, which happens whenever a group of symmetries of $X$ has the relations of the scheme as its orbits on pairs, then $B$ lies in the algebra spanned by the relation matrices of the scheme (its Bose--Mesner algebra), so it is a combination $\sum_j\beta_jE_j$ of the primitive idempotents, and $\dim(\conv(X))$ is the sum of the ranks of the $E_j$ with $\beta_j\neq0$ (Proposition~\ref{prop:principle}). The $\beta_j$ are computed from the eigenmatrix of the scheme and a few numbers $\tau_i$ that record how often two elements of $E$ in the $i$th relation occur together in a member of $X$. The kernel of $B$, which is the sum of the images of the $E_j$ with $\beta_j=0$, is exactly the space of implicit equations of $\conv(X)$: $By=0$ if and only if $y^\top x=y^\top\bar x$ for every $x\in X$; so any basis of the kernel is a maximal independent set of equations for $X$. The method was developed for the Johnson schemes of complete graphs and hypergraphs, extended to coherent configurations by Lee and Lee~\cite{LeeLee}, and used to find facets of the acyclic subgraph and linear ordering polytopes by Leung and Lee~\cite{LL90,LL26}. Those schemes have few classes and are tied to graphs. Here we use schemes whose ground set carries algebraic structure.

\paragraph{Two new dimension formulae.}
In general the spectral method needs the full eigenmatrix of the scheme, which is rarely available. Our first contribution is two formulae that avoid it; every result below applies one of them. The first is for \emph{amorphic} schemes, those in which every fusion of relations is again a scheme; with three or more classes, all their relations are strongly regular graphs of Latin-square or negative Latin-square type~\cite{vDM}. Such schemes are plentiful and natural: the parallel classes of any net, including the rows, columns, and symbols of any Latin square, form one, as do many cyclotomic schemes of finite fields~\cite{vDM}. Their eigenmatrix is determined by the valencies alone, and the dimension becomes a \emph{threshold formula}: a relation contributes its full rank or nothing, according as its count $\tau_j$ differs from one critical value or not (Theorem~\ref{thm:amorphic}). The second is for point sets in $\GF(q)^d$ invariant under the affine maps $z\mapsto az+b$, where no amorphic scheme need be available: the dimension is $q-1$ times the number of hyperplane directions along which some member of $X$ is not equidistributed, plus one if the coordinate sum varies (Theorem~\ref{thm:affine}). In both cases a lower bound needs only one member of $X$ for each eigenspace, with a nonzero component there, which we call a \emph{witness}; for the second formula, a member that is not equidistributed along the given direction. Linear codes and complete mappings of groups (bijections $\sigma$ for which $a\mapsto a+\sigma(a)$ is also a bijection) will supply witnesses. A witness is a certificate of the cheapest kind, one integer point violating a candidate equation, while the scheme supplies the matching upper bound.

\paragraph{Results.}
\begin{enumerate}
\item[(a)] \emph{Transversals of nets.} The polytope of transversals of a Desarguesian $k$-net of order $q$ has dimension $(q-1)(q+1-k)$ for $2\le k\le q$ (Corollary~\ref{cor:nets}). For $k=2$ this is the Birkhoff polytope; for $k=3$ it is the transversal polytope of the Cayley table of $(\GF(q),+)$.
\item[(b)] \emph{Multi-index assignment polytopes.} The $(k,s)$ assignment polytope $P^{(k,s)}_n$ of Appa, Magos, and Mourtos~\cite{AMM} is the convex hull of the orthogonal arrays $\mathrm{OA}(n^s,k,n,s)$ of index one; it is the polytope of the multi-index assignment problems, including the axial and planar three-index assignment problems and the problem of pairs of orthogonal Latin squares. For prime powers $q$ and $1\le s<k\le q+1$ we determine its dimension: it equals that of the linear relaxation except when $(k,s)=(q+1,q-1)$, where the arrays are perfect codes and the dimension is less by exactly $(q-1)^{q+1}$ (Theorem~\ref{thm:assign}). A second exceptional family occurs at $k=q+2$ for $q\ge4$ a power of $2$ (Theorem~\ref{thm:qplus2}), and products give full-dimensional cases at orders that are not prime powers (Proposition~\ref{prop:products}). As far as we can tell, the exceptional cases are the first counterexamples to the expectation of~\cite{AMM} that the two dimensions always coincide. The missing equations are clique equations that the standard formulation omits (\S\ref{sec:ip}).
\item[(c)] \emph{Transversals of group tables.} For every abelian group with a transversal, the transversal polytope of its Cayley table, equivalently the polytope of rainbow perfect matchings of the corresponding edge-colored $K_{n,n}\,$, has the maximal dimension $(n-1)(n-2)$: the line sums are its only implicit equations (Theorem~\ref{thm:main}). For each candidate equation not forced by the line sums, the proof exhibits a transversal that violates it (\S\ref{sec:latin}).
\item[(d)] \emph{Symmetric matrices.} On the space of symmetric $3\times3$ or $4\times4$ matrices over $\GF(q)$, the differences of rank at least $3$ form a relation of Latin-square type, and the largest sets of matrices whose pairwise differences have rank at least $3$, of size $n=q^3$ and $n=q^5$ respectively by Schmidt's bound~\cite{Sch}, span a polytope of dimension $(n-q^2)(n-1)$, for every $q$ (Propositions~\ref{prop:symodd} and~\ref{prop:sym4}). The proofs count points and lines on quadrics.
\end{enumerate}

\paragraph{Relation to previous work.}
Appa, Magos, and Mourtos~\cite{AMM} determine the dimension of the axial polytopes ($s=1$) for $n\ge3$ and of the planar polytopes ($s=2$) for $n\ge\max\{5,k\}$, and regard their results as strong evidence that $\dim(P^{(k,s)}_n)$ always equals the dimension $L^{(k,s)}_n$ of the relaxation; the orthogonal Latin-square polytope was treated in~\cite{AMMJ}, the Latin-square polytope in~\cite{EBG}, and the axial three-index polytope in~\cite{BS}. Computationally, exact algorithms for these problems separate clique and other set-packing inequalities~\cite{BS91,AMM04,MM09}, and integer programming has been used to classify orthogonal arrays~\cite{BM08}. Bulutoglu~\cite{Bul} proved by representation theory that the deficit can only consist of whole levels of the analysis-of-variance decomposition, and that it vanishes unless $n$ divides a certain binomial coefficient; his criterion settles $2994$ of the $3217$ cases $(k,s,q)$ with $q\le32$ a prime power and $k\le q+1$, and he names $(n,k,s)=(10,6,2)$ as the smallest open candidate. Theorem~\ref{thm:assign} settles all cases with $k\le q+1$ at prime-power orders and shows that the deficit his theorem allows is attained in the exceptional family. The difference in method is that his symmetry group $S_n\wr S_k\,$, the permutations of $[n]^k$ that permute the $k$ coordinates and relabel the $n$ symbols separately in each coordinate, sees only whole levels, while the affine group sees each direction inside them (\S\ref{sec:lit}).

For group tables, the transversal polytope of the Cayley table of $(\GF(q),+)$ is a special case of (a); for general abelian groups the dimension does not seem to have been determined before. By Hall and Paige~\cite{HP}, an abelian group has a transversal exactly when its Sylow $2$-subgroup is trivial or not cyclic, and these are exactly the groups of Theorem~\ref{thm:main}. Nonabelian groups are not covered, because their symmetries give no commutative scheme; still, computation confirms the same dimension for all $35$ groups of order at most $21$ with a transversal.

For the direct counts of \S\ref{sec:symm} we use Schmidt's bound and constructions for symmetric rank-metric codes~\cite{Sch}; for even $q$, the common fissions of van Dam, Koolen, and Xiong~\cite{vDKX} give a further dimension formula. The complexity results that motivate (c) are those for the axial three-index assignment problem, which contains three-dimensional matching~\cite{Karp} (see Spieksma~\cite{Spi} for applications), and for rainbow matchings, which are NP-hard to maximize in edge-colored bipartite graphs~\cite{IRT}, even on properly colored paths~\cite{LP}; for the complete properly colored graphs that come from Latin squares we know of no such result. For a group table the transversals are the graphs of the complete mappings $\sigma$, and the orthomorphisms $a\mapsto a+\sigma(a)$ have been used to build block substitutions for encryption~\cite{Mit}.

\paragraph{Organization.}
\S\ref{sec:schemes} develops the spectral method and the threshold formula for amorphic schemes, \S\ref{sec:affine} the formula for affinely invariant point sets and the transversals of Desarguesian nets, \S\ref{sec:assign} the assignment polytopes, \S\ref{sec:latin} the transversals of group tables, \S\ref{sec:almost} almost amorphic schemes and the symmetric matrices, and \S\ref{sec:ip} the consequences for integer-programming formulations. All constructions were verified computationally.

\section{The spectral method and amorphic schemes}\label{sec:schemes}

This section reviews the spectral method (see~\cite{L1,L2,L3}) from the linear algebra up, then the association schemes through
which it is applied, and then amorphic schemes, whose explicit spectra are what make the method effective for the polytopes of
the later sections. Each stage comes with a small worked example.

\subsection{Dimension as a rank}\label{sec:rank}

Let $X\subseteq\R^E$ be finite, with average $\bar x$. The vectors
$x-\bar x$, $x\in X$, span the space parallel to $\aff(X)$, so $\dim(\conv(X))=\dim(\aff(X))$ is the rank of the matrix having them as
rows, and hence of
\[
B\coloneqq\sum_{x\in X}(x-\bar x)(x-\bar x)^\top,
\]
because $\operatorname{rank}(M^\top M)=\operatorname{rank}(M)$. The matrix $B$ is symmetric and positive semidefinite, so the dimension is the
number of its nonzero eigenvalues, counted with multiplicity.

Two features make this useful. First, $B$ is indexed by the ground set $E$, not by $X$: for the Hamiltonian tours of the complete
graph $K_n$\,, $|E|=\binom n2$ while $|X|=(n-1)!/2$. Second, $B$ inherits the symmetries of $X$: if a group $\Gamma$ of permutations
of $E$ maps $X$ onto itself, then $\Gamma$ fixes $\bar x$, so $B$ commutes with the permutation matrices of $\Gamma$. The more
commutative the algebra generated by those matrices, the more explicitly $B$ can be diagonalized.

\begin{example}[cyclic symmetry]\label{ex:cyclic}
Let $E\coloneqq\mathbb Z_6$\,, and let $X$ consist of the incidence vectors of the six pairs $\{a,a+1\}$, that is, the six cyclic shifts of
$v\coloneqq(1,1,0,0,0,0)$. Cyclic shifts map $X$ onto itself, so $B$ commutes with the shift, which means that $B$ is a circulant matrix.
Circulant matrices are diagonalized by the characters $\chi_j(a)\coloneqq\omega^{ja}$, where $\omega\coloneqq e^{2\pi i/6}$, and a standard computation demonstrates
that the eigenvalue of $B$ at $\chi_j$ is $|\hat u(j)|^2$, where $u\coloneqq v-\bar x$ and $\hat u(j)\coloneqq\sum_au_a\omega^{-ja}$. For $j=0$ this is
$0$, because $\sum_au_a=0$; for $j\neq0$ it is $|1+\omega^{-j}|^2$, which vanishes only for $j=3$. So $\operatorname{rank}(B)=4$, and $\dim(\conv(X))=4$.
Each vanishing eigenvalue is an equation holding on $X$: here the alternating sum $\sum_a(-1)^ax_a$ takes the same value ($0$) on
every member.
\end{example}

\subsection{Association schemes and concordance}\label{sec:principle}

Cyclic symmetry is the simplest instance of a structure that has been studied for its own sake. An \emph{association scheme} on $E$
with $d$ classes consists of symmetric $0/1$ matrices $A_0=I,A_1\,,\dots,A_d$ indexed by $E$, with $\sum_iA_i=J$ (the all-ones matrix),
whose linear span is closed under multiplication. We think of $A_i$ as the adjacency matrix of a graph on $E$, the $i$th
\emph{relation}; the relations partition the pairs of elements of $E$. The span is the \emph{Bose--Mesner algebra}. Because the
$A_i$ are symmetric and span an algebra, they commute, and the algebra has a basis of orthogonal projections
$E_0=\frac1{|E|}J,E_1\,,\dots,E_d$\,, the \emph{primitive idempotents}, with $E_jE_l=\delta_{jl}E_j$ and $\sum_jE_j=I$. Write
$m_j\coloneqq\operatorname{rank}(E_j)$\,. Each $A_i$ acts on the image of $E_j$ as a scalar, $A_iE_j=P_{ji}E_j$\,, and the matrix $P$ is the
\emph{eigenmatrix} of the scheme. Typical sources are group actions whose orbits on pairs of elements of $E$ are symmetric: for
instance the symmetric group acting on the edges of $K_n$ (Example~\ref{ex:tours} below), or an abelian group acting on itself by
translation, with two pairs in the same relation when their differences agree up to sign (\S\ref{sec:affine}).

\begin{definition}
A finite set $X\subseteq\R^E$ is \emph{concordant} with a scheme on $E$ if (a) $\sum_{x\in X}x_ex_f=\tau_i$ for all
$(e,f)$ in the $i$th relation (with $\tau_0\coloneqq\sum_{x\in X}x_e^2$), and (b) $\sum_{x\in X}x_e$ is independent of $e$.
\end{definition}

For $0/1$ vectors, $\tau_i$ is the number of members of $X$ containing both $e$ and $f$. Both conditions hold whenever a group
maps $X$ onto itself and its orbits on pairs are the relations, because the sums are then constant on orbits.

\begin{proposition}[spectral principle; see~\cite{L2}]\label{prop:principle}
Let $X$ be concordant with a scheme. Then $B$ lies in the Bose--Mesner algebra, $B=\sum_j\beta_jE_j$\,, with $\beta_j=\sum_i\tau_iP_{ji}$
for $j\ge1$, and
\begin{equation}\label{eq:principle}
\dim(\conv(X))=\dim(\mathbf 1^\top X)+\sum\bigl\{m_j ~:~ j\ge1,\ \beta_j\neq0\bigr\},
\end{equation}
where $\dim(\mathbf 1^\top X)\in\{0,1\}$ is $1$ exactly when the coordinate sum varies on $X$. Moreover, for $j\ge1$, $\beta_j\neq0$ if and only if
$E_jx\neq0$ for some $x\in X$; and $\aff(X)$ is the set of $y$ with $E_jy=0$ whenever $j\ge1$ and $\beta_j=0$, and with
$\mathbf 1^\top y$ equal to the common value of $\mathbf 1^\top x$ on $X$, if any.
\end{proposition}

\begin{proof}
By (b), $\bar x=c\mathbf 1$ for a constant $c$, and by (a), $\sum_xxx^\top=\sum_i\tau_iA_i$; hence
$B=\sum_i\tau_iA_i-|X|c^2J$ lies in the algebra. On the image of $E_j$\,, $j\ge1$, we have $J=0$ and $A_i=P_{ji}$\,, so
$\beta_j=\sum_i\tau_iP_{ji}$\,. On the image of $E_0$\,, the eigenvalue is $\beta_0=\frac1{|E|}\sum_x(\mathbf 1^\top x-\mathbf 1^\top\bar x)^2$, which is positive
exactly when $\mathbf 1^\top x$ varies. Because $B$ is positive semidefinite, its rank is the sum of $m_j$ over the $j$ with $\beta_j\neq0$, which
is~\eqref{eq:principle}. For $j\ge1$, $\operatorname{tr}(BE_j)=\beta_jm_j=\sum_x\|E_j(x-\bar x)\|^2=\sum_x\|E_jx\|^2$, which gives the
second statement. Finally, $\aff(X)$ lies in the affine space described, and the two have the same dimension.
\end{proof}

The statement in~\cite{L2} is for \emph{balanced} polytopes, on which $\mathbf 1^\top x$ is constant; centering at $\bar x$, as in the
definition of $B$, puts the variation of $\mathbf 1^\top x$ into $E_0$ and removes that restriction. The equations in the last statement, with
bases of the images of the vanishing $E_j$\,, are the canonical equation sets of~\cite{L3}. We use only symmetric schemes; the
extension of the method to coherent configurations (see~\cite{LeeLee}) allows nonsymmetric relations and several fibers.

\begin{example}[Hamiltonian tours]\label{ex:tours}
Let $E$ be the edge set of $K_n$\,, $n\ge4$, with two relations besides equality: edges sharing a vertex, and disjoint edges. This is
a scheme, with idempotents of ranks $1$, $n-1$, and $\binom n2-n$; the eigenvalues of the two relations on them are
$(2n-4,\,n-4,\,-2)$ and $(\binom{n-2}2,\,3-n,\,1)$ (see~\cite{L3}, Table~1). Let $X$ be the set of Hamiltonian tours, that is, of
the incidence vectors of the $N\coloneqq(n-1)!/2$ cycles through all $n$ vertices. The symmetric group maps $X$ onto itself, so $X$ is
concordant, and counting incidences gives $\tau_0=\frac{2N}{n-1}$, $\tau_1=\frac{2N}{(n-1)(n-2)}$, and $\tau_2=\frac{4N}{(n-1)(n-2)}$. With
$u\coloneqq\frac{2N}{(n-1)(n-2)}$ these are $(n-2)u,\,u,\,2u$, so
\begin{align*}
\beta_1&=\tau_0+(n-4)\tau_1+(3-n)\tau_2=\bigl[(n-2)+(n-4)-2(n-3)\bigr]u=0,\\
\beta_2&=\tau_0-2\tau_1+\tau_2=(n-2)u\neq0.
\end{align*}
As $\mathbf 1^\top x\equiv n$ is constant, $\dim(\conv(X))=m_2=\binom n2-n=\frac n2(n-3)$, the classical dimension of the tour polytope.
The vanishing of $\beta_1$\,, together with $\mathbf 1^\top x\equiv n$, is the family of the $n$ degree equations $x(\delta(v))=2$, one for each
vertex $v$: the images of $E_0$ and $E_1$ together are spanned by the indicator vectors of the vertex stars (see~\cite{L3}).
\end{example}

Witnesses sharpen a device of~\cite[Corollary~3.2]{L2}: there, because concordance restricts the dimension to a few values (for balanced polytopes on the edges of $K_n\,$, to $0$, $n-1$, $\binom n2-n$, or $\binom n2-1$), $n+1$ affinely independent points and two independent equations already pin the dimension to $\binom n2-n$. Here a single point $x$ with $E_jx\neq0$ certifies the whole eigenspace of $E_j$ (as $E_j\bar x=0$ for $j\ge1$), and this is what the witnesses of Theorem~\ref{thm:affine} and Proposition~\ref{prop:fourier} do.

\subsection{Amorphic schemes}\label{sec:amorphic}

A scheme is \emph{amorphic} if every fusion of its relations, that is, every way of partitioning the nontrivial relations into
groups and adding up the matrices in each group, is again a scheme. We need this notion because the spectra of amorphic schemes are
completely explicit. A graph on $n^2$ vertices is \emph{strongly regular of Latin-square type} if it is regular of valency
$t(n-1)$ and its \emph{restricted} eigenvalues (those with eigenvectors orthogonal to the all-ones vector) are $n-t$ and $-t$, for positive integers $n$ and $t$; it is of
\emph{negative} Latin-square type if $n$ and $t$ are negative integers. Examples are a disjoint union of $n$ cliques of size $n$ ($t=1$; a
\emph{square spread}), the lattice graph $L_2(n)$ on the cells of an $n\times n$ array, adjacent when they share a row or a column
($t=2$), and, for a prime power $n=q$, the graph on the points of the affine plane $\AG(2,q)$ in which two points are adjacent
when the line through them lies in one of $t$ fixed parallel classes. We use the following facts
(van Dam and Muzychuk~\cite{vDM}; see also~\cite{vDKX}, Section~2).
\begin{enumerate}
\item[(i)] If an amorphic scheme has $d\ge3$ classes, then either all its relations are strongly regular of Latin-square type or all
are of negative Latin-square type. (More generally, no scheme at all has both a relation of strictly Latin-square type and one of strictly negative
Latin-square type; see~\cite{vDKX}, Theorem~7.1. So only one eigenvalue pattern ever has to be considered.)
\item[(ii)] If relation $j$ has parameters $(n,t_j)$, then $\sum_jt_j=n+1$ (the valencies $t_j(n-1)$ add up to $n^2-1$), and the
nontrivial idempotents can be indexed $E_1\,,\dots,E_d$ so that $m_j=t_j(n-1)$ and
\[
A_iE_j=(n\delta_{ij}-t_i)E_j\qquad(i,j=1,\dots,d).
\]
In words: relation $i$ has eigenvalue $n-t_i$ on its own idempotent and $-t_i$ on all the others. (As a check,
$\sum_jm_j=(n-1)(n+1)=n^2-1$.)
\item[(iii)] The prototype is $\AG(2,q)$ with its $q+1$ parallel classes of lines: relation $j$ joins two points on a common line of
class $j$, so $n=q$ and $t_j=1$, and $E_j$ is the projection onto the mean-zero functions constant on the lines of class $j$, a space of dimension $q-1=t_j(n-1)$.
\end{enumerate}
By (ii), the eigenmatrix is known in closed form for any number of classes, and in the spectral principle the eigenmatrix is the only
input besides the counts $\tau_i$\,. A second reason for working with amorphic schemes is that every fusion of an amorphic scheme is
amorphic, so for a point set whose symmetry group is small, only the coarsest grouping of relations that the group can distinguish
need be checked for concordance.

\begin{theorem}[threshold formula]\label{thm:amorphic}
Let $\mathcal A\coloneqq\{I,A_1\,,\dots,A_d\}$ be an amorphic scheme on $E$, $|E|=n^2$, with parameters $(n,t_j)$ as in~(ii), and let
$X\subseteq\R^E$ be concordant with it, with counts $\tau_0\,,\tau_1\,,\dots,\tau_d$\,. Put
$\tau^*\coloneqq\frac1n\bigl(\sum_it_i\tau_i-\tau_0\bigr)$. Then
\[
\dim(\conv(X))=\dim(\mathbf 1^\top X)+(n-1)\sum\bigl\{t_j ~:~ \tau_j\neq\tau^*\bigr\}.
\]
Moreover, $\tau_j=\tau^*$ if and only if every $x\in X$ is orthogonal to the image of $E_j$\,.
\end{theorem}

\begin{proof}
By Proposition~\ref{prop:principle} and~(ii), $\beta_j=\tau_0+\sum_i\tau_i(n\delta_{ij}-t_i)=n(\tau_j-\tau^*)$, and $m_j=t_j(n-1)$.
\end{proof}

So each relation $j$ is either \emph{balanced}, $\tau_j=\tau^*$, in which case it forces $m_j$ equations on $X$, or it contributes
its full rank $m_j=(n-1)t_j$ to the dimension; there is no intermediate case. When the relations come from the parallel classes of a net
as in~(iii), the image of $E_j$ consists of the functions that are constant on the lines of the classes in relation $j$ and sum to zero,
so $\tau_j=\tau^*$ means that every member of $X$ has the same sum over every line of every such class.

\begin{figure}[t]
\centering
\begin{tikzpicture}[scale=0.5,every node/.style={font=\scriptsize}]
\foreach \k/\lab/\col in {0/{rows\\balanced}/black!60, 1/{columns\\balanced}/black!60, 2/{$x+y$\\witnessed}/blue!70, 3/{$x-y$\\witnessed}/red!70}{
 \begin{scope}[xshift=\k*4.6cm]
  \foreach \i in {0,1,2}\foreach \j in {0,1,2}\fill (\i,\j) circle (0.11);
  \ifnum\k=0 \foreach \c in {0,1,2}\draw[\col,thick] (0,\c)--(2,\c);\fi
  \ifnum\k=1 \foreach \c in {0,1,2}\draw[\col,thick] (\c,0)--(\c,2);\fi
  \ifnum\k=2 \draw[\col,thick] (0,0)--(1,2)--(2,1); \draw[\col,thick] (0,1)--(1,0)--(2,2); \draw[\col,thick] (0,2)--(1,1)--(2,0);\fi
  \ifnum\k=3 \draw[\col,thick] (0,0)--(1,1)--(2,2); \draw[\col,thick] (0,1)--(1,2)--(2,0); \draw[\col,thick] (0,2)--(1,0)--(2,1);\fi
  \node[align=center] at (1,-0.95) {\lab};
 \end{scope}}
\begin{scope}[xshift=18.4cm]
  \foreach \i in {0,1,2}\foreach \j in {0,1,2}\fill[black!30] (\i,\j) circle (0.11);
  \foreach \i/\j in {0/1,1/2,2/0}\fill[red!70] (\i,\j) circle (0.18);
  \node[align=center] at (1,-0.95) {a transversal\\of the $3$-net};
\end{scope}
\end{tikzpicture}
\caption{The parallel classes of $\AG(2,3)$ and a transversal of the $3$-net $\{\text{rows},\text{columns},x+y\}$, the Cayley table of $\mathbb Z_3\,$. Transversals meet every line of the net once, so its classes are balanced; the one shown is a line of the fourth class and witnesses that direction. The transversal polytope has dimension $2$ for this net and $4$ for the $2$-net $\{\text{rows},\text{columns}\}$ (permutation matrices).}
\label{fig:ag23}
\end{figure}
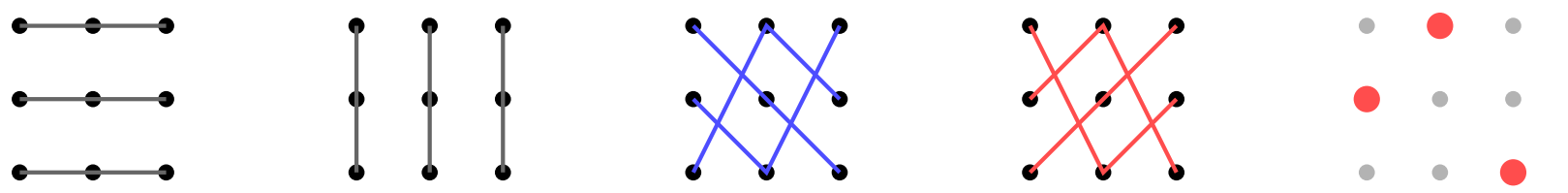
\begin{example}[the $3\times3$ permutation matrices]\label{ex:birkhoff}
Let $E\coloneqq\GF(3)^2$ be the nine cells of a $3\times3$ array, with its four parallel classes: rows, columns, and the two kinds of
diagonal (the scheme of (iii) with $q=3$, so $n=3$ and all $t_j=1$; Figure~\ref{fig:ag23}). Let $X$ be the six permutation matrices; each is a set of three
cells meeting every row and every column once, and the six are exactly the six diagonals. Translations of the array and the
map $z\mapsto-z$ preserve $X$ and have the four classes as the relations, so $X$ is concordant. Each cell lies in two members, so
$\tau_0=2$. Two cells in a common row or column lie in no common member: $\tau=0$ for these two classes. Two cells on a common diagonal lie
on exactly one common member (that diagonal): $\tau=1$ for the other two classes. Hence
$\tau^*=\frac13\bigl(0+0+1+1-2\bigr)=0$. Rows and columns are balanced, while the two diagonal classes are not and contribute
$(n-1)t_j=2$ each; because $\mathbf 1^\top x\equiv3$, the dimension is $2+2=4=(3-1)^2$, the dimension of the Birkhoff polytope of order $3$.
The balanced classes give the familiar equations: every row sum and every column sum equals $1$. A permutation matrix that is a diagonal of one of the two diagonal classes is not equidistributed along that class
(one line has sum $3$, the other two $0$); we will call such a member a \emph{witness} for that class. Replacing $3$ by an arbitrary prime power $q$ and adding more nets leads to
\S\ref{sec:affine}.
\end{example}

\section{Affine symmetry in $\GF(q)^d$}\label{sec:affine}

We now take the ground set to be a vector space over a finite field, $E=V=\GF(q)^d$, and let the affine maps of $V$ act as symmetries.
Point sets of finite geometry and design theory (lines, nets, codes, arrays) are of this kind. The scheme that results is the
translation scheme of $V$, but the proof needs no theory of schemes, only the Fourier transform, which we review first. The
structure is that of Example~\ref{ex:cyclic}, with $\mathbb Z_6$ replaced by $(V,+)$.

\subsection{Fourier analysis on $\GF(q)^d$}\label{sec:fourier}

Write $q=p^m$ with $p$ prime, let $\operatorname{Tr}\colon\GF(q)\to\GF(p)$ be the trace, and put
$\psi(a)\coloneqq e^{2\pi i\operatorname{Tr}(a)/p}$, a nontrivial homomorphism from $(\GF(q),+)$ to the unit circle. (For $q=p$, $\operatorname{Tr}$ is the
identity and $\psi(a)=e^{2\pi ia/p}$.) For $w\in V$ let $\psi_w(z)\coloneqq\psi(w\cdot z)$, where $w\cdot z\coloneqq\sum_iw_iz_i$\,. Vectors of $\mathbb C^V$
carry the inner product $\langle f,g\rangle\coloneqq\sum_zf(z)\overline{g(z)}$, and $\hat f(w)\coloneqq\langle f,\psi_w\rangle$.

\begin{enumerate}
\item[\textbf{F1}] \emph{Characters.} $\sum_{z\in V}\psi_w(z)$ is $|V|$ for $w=0$ and $0$ for $w\neq0$. Hence the $|V|$ vectors $\psi_w$ are
pairwise orthogonal, so they form a basis of $\mathbb C^V$, and $f=|V|^{-1}\sum_w\hat f(w)\psi_w$\,. (The trace is onto $\GF(p)$,
because a polynomial of degree $p^{m-1}<q$ cannot vanish on all of $\GF(q)$, so $\psi$ is nontrivial; if $w\neq0$, the map $z\mapsto w\cdot z$ is onto
$\GF(q)$ with fibers of size $q^{d-1}$, so $\sum_z\psi_w(z)=q^{d-1}\sum_{c\in\GF(q)}\psi(c)=0$.)
\item[\textbf{F2}] \emph{Translation-invariant matrices.} Call a matrix $B$ indexed by $V$ \emph{translation invariant} if
$B(z+b,z'+b)=B(z,z')$ for all $b$. Then every $\psi_w$ is an eigenvector of $B$:
\[
(B\psi_w)(z)=\sum_uB(z,z+u)\,\psi_w(z+u)=\psi_w(z)\sum_uB(0,u)\,\psi_w(u).
\]
This is the familiar diagonalization of circulant matrices by the characters of $\mathbb Z_n$\,, for the group $(V,+)$, and so the rank
of $B$ is the number of $w$ with a nonzero eigenvalue.
\item[\textbf{F3}] \emph{Hyperplane sums.} For $0\neq w\in V$ the sets $H_{w,c}\coloneqq\{z\in V ~:~ w\cdot z=c\}$, $c\in\GF(q)$, are the hyperplanes of a
\emph{parallel class}, which depends only on the \emph{direction} $[w]\coloneqq\{aw ~:~ a\in\GF(q)^*\}$, a point of the projective space
$\PG(d-1,q)$; there are $(q^d-1)/(q-1)$ directions. For $x\in\R^V$ let $S_c(x)\coloneqq\sum_{z\in H_{w,c}}x_z$\,. Then for $a\in\GF(q)^*$,
\[
\hat x(aw)=\sum_{c\in\GF(q)}S_c(x)\,\overline{\psi(ac)},
\]
the one-dimensional Fourier transform of $c\mapsto S_c(x)$. By F1 for $d=1$, this vanishes for every $a\neq0$ if and only if $S_c(x)$ does
not depend on $c$. We say that $x$ is \emph{equidistributed along $[w]$} in this case, and that $x$ \emph{witnesses} $[w]$ otherwise.
Equivalently, $x$ is orthogonal to the space $U_{[w]}$ of functions $z\mapsto g(w\cdot z)$ with $\sum_cg(c)=0$, which has dimension $q-1$
and is spanned over $\mathbb C$ by $\psi_{aw}$\,, $a\neq0$. By F1 again, $\R^V$ is the orthogonal sum of the constants and the spaces
$U_{[w]}$\,, one for each direction.
\end{enumerate}

\subsection{The dimension formula}\label{sec:affineformula}

Let $\mathcal H\coloneqq\{z\mapsto az+b ~:~ a\in\GF(q)^*,\ b\in V\}$ act on $\R^V$ by $(x_g)_z\coloneqq x_{g^{-1}z}$\,. We call a set $X\subseteq\R^V$
\emph{$\mathcal H$-invariant} if $x_g\in X$ for every $x\in X$ and $g\in\mathcal H$.

\begin{theorem}\label{thm:affine}
Let $X\subseteq\R^V$ be finite and $\mathcal H$-invariant, and let $W(X)$ be the set of directions witnessed by some member of $X$.
Then
\[
\dim(\conv(X))=\dim(\mathbf 1^\top X)+(q-1)\,|W(X)|,
\]
and $\aff(X)$ is the set of $y\in\R^V$ that are equidistributed along every direction not in $W(X)$ and whose coordinate sum is the
common value of $\mathbf 1^\top x$ on $X$, if any.
\end{theorem}

\begin{proof}
Let $u_x\coloneqq x-\bar x$, so that $B=\sum_{x\in X}u_xu_x^\top\,$.

\emph{Step 1: $B$ is translation invariant.} The average $\bar x$ is invariant under all of $\mathcal H$, in particular under translations,
so it is a constant vector. A translation permutes $X$ and hence permutes the vectors $u_x$\,, and therefore fixes $B=\sum_xu_xu_x^\top\,$.

\emph{Step 2: the eigenvalues.} By F2, each $\psi_w$ is an eigenvector of $B$, with some eigenvalue $\lambda_w$\,. Because $B$ is real
and $\|\psi_w\|^2=|V|$, and the $u_x$ are real,
\[
\lambda_w|V|=\langle B\psi_w\,,\psi_w\rangle=\sum_{x\in X}\bigl|\langle u_x\,,\psi_w\rangle\bigr|^2 .
\]
For $w=0$ this is $\sum_x(\mathbf 1^\top x-\mathbf 1^\top\bar x)^2$, which is nonzero exactly when $\mathbf 1^\top x$ varies on $X$. For $w\neq0$, $\langle\bar x,\psi_w\rangle=0$ by F1,
so $\langle u_x\,,\psi_w\rangle=\hat x(w)$, and $\lambda_w\neq0$ exactly when $\hat x(w)\neq0$ for some $x\in X$.

\emph{Step 3: homotheties.} For $g(z)=az$ we have $\widehat{x_g}(w)=\sum_yx_y\overline{\psi(aw\cdot y)}=\hat x(aw)$. Because $x\mapsto x_g$ permutes $X$,
$\hat x(aw)=0$ for all $x\in X$ if and only if $\hat x(w)=0$ for all $x\in X$. So $\lambda_{aw}\neq0$ if and only if $\lambda_w\neq0$: the
vanishing of $\lambda_w$ depends only on the direction of $w$.

\emph{Step 4: counting.} By F3, for $w\neq0$ we have $\hat x(aw)\neq0$ for some $a\neq0$ and some $x\in X$ exactly when some member witnesses
$[w]$. By Step 3 this makes all $q-1$ eigenvalues $\lambda_{aw}$\,, $a\neq0$, nonzero, and otherwise all vanish. Hence
$\dim(\conv(X))=\operatorname{rank}(B)=\dim(\mathbf 1^\top X)+(q-1)|W(X)|$.

\emph{Step 5: the affine hull.} Every $x\in X$ lies in the affine space $A$ described in the statement, so $\aff(X)\subseteq A$. Because $\R^V$ is
the orthogonal sum of the constants and the $(q^d-1)/(q-1)$ spaces $U_{[w]}$ (by F3), each of dimension $q-1$, the space $A$ has dimension
$\dim(\mathbf 1^\top X)+(q-1)|W(X)|$. The dimensions agree, so $\aff(X)=A$.
\end{proof}

In the language of \S\ref{sec:schemes}: the orbits of $\mathcal H$ on pairs of points are the diagonal and, for each direction
$[u]$, the relation joining $z$ and $z+\lambda u$, $\lambda\neq0$. These relations form the \emph{translation scheme} of $V$, whose
idempotents are the projections onto the constants and onto the spaces $U_{[w]}$\,, and Theorem~\ref{thm:affine} is the spectral principle
for it. For $d=2$ this scheme is the amorphic scheme of $\AG(2,q)$ in (iii) of \S\ref{sec:amorphic} (each relation a disjoint
union of $q$ cliques of size $q$, with $t_j=1$), and the theorem is then the threshold formula with concordance supplied by symmetry.
For $d\ge3$ the scheme is not amorphic, but its idempotents are just as explicit.

\begin{example}[the $3\times3$ permutation matrices again]\label{ex:birkhoff2}
Take $q=3$, $d=2$, and let $X$ be the six permutation matrices of Example~\ref{ex:birkhoff}, on the cells $z=(z_1\,,z_2)\in\GF(3)^2$.
There are four directions, represented by $w=(1,0)$, $(0,1)$, $(1,1)$, and $(1,2)$, whose parallel classes are the rows $z_1=c$, the columns
$z_2=c$, and the two kinds of diagonal $z_1+z_2=c$ and $z_1-z_2=c$. Every permutation matrix has exactly one cell in each row and in
each column, so $S_c(x)=1$ for all $c$ along the first two directions: no member witnesses them. The anti-diagonal $\{z_1+z_2=0\}$ is a member, with
$S_0=3$ and $S_1=S_2=0$, so it witnesses $[(1,1)]$; likewise the diagonal $\{z_1=z_2\}$ witnesses $[(1,2)]$. Hence $|W(X)|=2$, and because $\mathbf 1^\top x\equiv3$, the
theorem gives $\dim=2\cdot(3-1)=4$, as before. In terms of the nine characters, $B$ has eigenvalue $0$ at $\psi_0$ and at the four characters
$\psi_{\pm(1,0)}\,,\psi_{\pm(0,1)}$: these give the equations that the coordinate sum is $3$, that all row sums are equal, and that all column sums are equal,
which together say that every row sum and every column sum is $1$ (five independent equations). The remaining four characters have nonzero eigenvalues, so $\operatorname{rank}(B)=4=9-5$.
(For $q=3$ every permutation matrix happens to be a line; for larger $q$ they are not.)
\end{example}

\subsection{Transversals of Desarguesian nets}\label{sec:nets}

A \emph{$k$-net} of order $q$ consists of $q^2$ points and $k$ parallel classes of $q$ lines each, any two lines from different classes
meeting in exactly one point. The classes of $\AG(2,q)$ give \emph{Desarguesian} nets: a set $K$ of $k$ directions in $V=\GF(q)^2$,
$2\le k\le q$. A \emph{transversal} of the net is a $q$-subset of the plane meeting every line of every class in $K$ in exactly one point; let
$T_K$ be the set of transversals. (They are the common transversals of the $k-2$ mutually orthogonal Latin squares that the net encodes.)
For $k=2$, with the classes ``row'' and ``column'', the transversals are the permutation matrices of order $q$, so $\conv(T_K)$ is the Birkhoff polytope.
For $k=3$, with classes ``row'', ``column'', and ``symbol'' ($z_1+z_2=c$), a transversal is a permutation $\sigma$ of $\GF(q)$ such that
$a\mapsto a+\sigma(a)$ is also a permutation: a transversal of the Cayley table of $(\GF(q),+)$. We return to these in \S\ref{sec:latin}.

\begin{corollary}\label{cor:nets}
For $2\le k\le q$,
$\dim(\conv(T_K))=(q-1)(q+1-k)$, and
\[
\aff(T_K)=\bigl\{y\in\R^{\GF(q)^2} ~:~ y(\ell)=1 \text{ for every line } \ell \text{ of every class in } K\bigr\}.
\]
\end{corollary}

\begin{proof}
$T_K$ is $\mathcal H$-invariant, because every map in $\mathcal H$ sends lines of a class to lines of the same class. Every member has one point on
each line of each class in $K$, so it is equidistributed along those $k$ directions, and it has $q$ points, so $\mathbf 1^\top x\equiv q$. Now let $[w]\notin K$ be one of the other $q+1-k$
directions. Each line $\ell$ of the class $[w]$ meets every line of every other class in exactly one point, so $\ell\in T_K$; and $\ell$ has
$S_c=q$ for one value of $c$ and $S_{c'}=0$ for the others, so it witnesses $[w]$. Hence $W(T_K)$ consists of the $q+1-k$ directions outside
$K$, and Theorem~\ref{thm:affine} applies.
\end{proof}

The case $k=2$ is the Birkhoff dimension $(q-1)^2$, which is classical (see~\cite{BR}); for $k\ge3$ we do not know of the formula in the literature.
The condition $k\le q$ is needed for $T_K$ to be nonempty (any two points of a transversal lie on a common line, whose class is not in $K$); for $k=q$ the polytope is the simplex on the $q$ lines of the
last class, of dimension $q-1$.

\section{Assignment polytopes}\label{sec:assign}

The Birkhoff polytope has the permutation matrices as vertices; its multi-index generalizations, studied by Appa, Magos, and
Mourtos~\cite{AMM}, have as vertices the incidence vectors of Latin squares, of pairs of orthogonal Latin squares, and of their
higher analogues, the \emph{orthogonal arrays}. This section determines their dimensions when the alphabet size is a prime
power, finds exceptional cases in which they fall short of the dimension of the linear relaxation, and extends the result to
orders that are not prime powers. The tool is Theorem~\ref{thm:affine}, together with the spectral principle for the Hamming
scheme.

\subsection{Orthogonal arrays and their polytopes}\label{sec:oa}

An \emph{orthogonal array} $\mathrm{OA}(n^s,k,n,s)$ of index one is a set $C\subseteq[n]^k$ of $n^s$ \emph{rows} such that, for every
set of $s$ coordinates, the projection of $C$ onto those coordinates is a bijection onto $[n]^s$ (each $s$-tuple occurs exactly once).
The integer $s$ is its \emph{strength}. For $(k,s)=(2,1)$ the rows $(i,\sigma(i))$ are a permutation matrix; for $(3,1)$, a pair of
permutations $(i,\sigma(i),\tau(i))$ (the axial three-index assignments); for $(3,2)$, the rows $(i,j,L(i,j))$ of a Latin square $L$; for
$(4,2)$, a pair of orthogonal Latin squares. The \emph{$(k,s)$ assignment polytope} $P^{(k,s)}_n$ is the convex hull in $\R^{[n]^k}$ of the
incidence vectors $1_C$ of all such arrays. Write $X$ for the set of these vectors.

\emph{The equations.} Strength $s$ says that for every set $T$ of $s$ coordinates and every $t\in[n]^T$, the number of rows
of $C$ that agree with $t$ on $T$ is $1$. These are linear equations on $1_C$\,. The polytope of nonnegative solutions of
them, the \emph{linear relaxation}, contains the strictly positive point $n^{s-k}\mathbf 1$, so its dimension is that of the solution space,
which is
\[
L^{(k,s)}_n\coloneqq n^k-\sum_{j=0}^s\binom kj(n-1)^j
\]
(Proposition~\ref{prop:hamming} below explains this count). Appa, Magos, and Mourtos regard their results as strong evidence that
$\dim(P^{(k,s)}_n)=L^{(k,s)}_n$ whenever $P^{(k,s)}_n\neq\emptyset$, and Bulutoglu~\cite{Bul} established that the possible deficits are restricted. Our results
demonstrate that the equality can fail.

\emph{Levels.} The ground set $[n]^k$ carries the \emph{Hamming scheme}: two points are in the $i$th relation when they differ in
exactly $i$ coordinates. Its idempotents are the projections onto the \emph{ANOVA levels}, a decomposition that statisticians use for
$k$-way tables. For $u\subseteq[k]$ let $L_u\subseteq\R^{[n]^k}$ be the space of functions that depend only on the coordinates in
$u$ and have average zero over each of those coordinates when the others are fixed. Then $\dim(L_u)=(n-1)^{|u|}$,
$\R^{[n]^k}=\bigoplus_uL_u$ orthogonally, and the \emph{level} $U_r\coloneqq\bigoplus_{|u|=r}L_u$ has dimension $\binom kr(n-1)^r$. For $k=2$:
$U_0$ is the constants, $U_1$ consists of the functions $h_i+v_j$ with $\sum h=\sum v=0$ (row and column effects, as in
\S\ref{sec:almost}), and $U_2$ of the functions with all row and column sums zero (the interactions). In general $U_r$ is the eigenspace of the Hamming
graph $H(k,n)$ for the eigenvalue $k(n-1)-nr$, and the projections onto $U_0\,,\dots,U_k$ are the primitive idempotents of the scheme.

Strength $s$ says exactly that the marginal of $1_C$ on every set of $s$ coordinates is constant, that is, that $1_C$ is orthogonal to
$U_1\,,\dots,U_s$\,. Call a level $r$ \emph{witnessed} if some member of $X$ has a nonzero component in $U_r$\,.

\begin{proposition}\label{prop:hamming}
$\dim(P^{(k,s)}_n)=\sum\bigl\{\dim(U_r) ~:~ r>s \text{ is witnessed}\bigr\}$. In particular $\dim(P^{(k,s)}_n)\le\sum_{r>s}\dim(U_r)=L^{(k,s)}_n$\,,
and the deficit $L^{(k,s)}_n-\dim(P^{(k,s)}_n)$ is a sum of some of the $\dim(U_r)$.
\end{proposition}

\begin{proof}
The group $S_n\wr S_k$ (permutations of the symbols in each coordinate, and of the coordinates) maps orthogonal arrays to orthogonal arrays, so it
maps $X$ onto itself, and its orbits on pairs of points are exactly the Hamming relations. So $X$ is concordant with the Hamming scheme, and
because $\mathbf 1^\top x\equiv n^s$ is constant, Proposition~\ref{prop:principle} gives the dimension as the sum of $\dim(U_r)$ over the levels $r\ge1$ at
which some $x\in X$ has $E_rx\neq0$. Strength $s$ excludes $r\le s$, and $\sum_{r=0}^k\binom kr(n-1)^r=n^k$.
\end{proof}

This is Bulutoglu's decomposition~\cite{Bul}, obtained here from the spectral principle. When the alphabet is a field, the affine group refines
this: for invariant point sets the dimension is counted by directions (Theorem~\ref{thm:affine}), and
each level splits into direction pieces.

\subsection{Witnesses from linear codes}\label{sec:codes}

Let $n=q$ be a prime power, so that $V=[q]^k=\GF(q)^k$. The set $X$ is invariant under the group $\mathcal H$ of \S\ref{sec:affineformula}
(the map $z\mapsto az+b$ applies the affine map $z_i\mapsto az_i+b_i$ to the $i$th coordinate, a permutation of the symbols), and also under
permutations of the coordinates and under coordinatewise scalings $z_i\mapsto a_iz_i$\,. Moreover the level $U_r$ is the sum of the spaces $U_{[w]}$ of
\S\ref{sec:fourier} over the directions $[w]$ of \emph{weight} $r$, where the weight of $w$ is its number of nonzero coordinates: $L_u$ is spanned
by the characters $\psi_w$ with $\operatorname{supp}w=u$, and there are $(q-1)^{|u|}$ of these. By Theorem~\ref{thm:affine}, because $\mathbf 1^\top x\equiv q^s$,
\[
\dim(P^{(k,s)}_q)=(q-1)\cdot\#\{\text{witnessed directions}\},
\]
and strength $s$ rules out all directions of weight at most $s$. So the dimension is $L^{(k,s)}_q$ exactly when every direction of weight at least $s+1$ is
witnessed. To witness a direction we need a member of $X$ that is not equidistributed along it, and linear codes supply them.

\emph{Linear codes.} A \emph{linear code} of \emph{length} $k$ over $\GF(q)$ is a subspace $C\le\GF(q)^k$; its \emph{dual} is
$C^\perp\coloneqq\{w ~:~ w\cdot z=0 \text{ for all } z\in C\}$. The \emph{weight} of a word is its number of nonzero coordinates. A code of dimension $r$ and minimum weight $d$ is a
$[k,r,d]$ code, and $d\le k-r+1$ (the Singleton bound); codes with $d=k-r+1$ are \emph{MDS}. The dual of an MDS code is MDS. A linear code $C$ of dimension $s$
whose dual has minimum weight at least $s+1$ is an $\mathrm{OA}(q^s,k,q,s)$, and for $k\le q+1$ and every $1\le s<k$ such codes exist: the generalized
Reed--Solomon codes, which are the images of the polynomials of degree less than $s$ under evaluation at $k$ distinct points of $\GF(q)\cup\{\infty\}$, are MDS.
The two facts we use are these.

\begin{lemma}\label{lem:linear}
Let $C\le\GF(q)^k$ be linear. Then $\widehat{1_C}(w)=|C|$ if $w\in C^\perp$ and $0$ otherwise. Hence $C$ witnesses exactly the directions $[w]$ with
$0\neq w\in C^\perp$, and its component in the level $U_r$ is nonzero if and only if $C^\perp$ has a word of weight $r$.
\end{lemma}

\begin{proof}
$\widehat{1_C}(w)=\sum_{z\in C}\overline{\psi(w\cdot z)}$. If $w\in C^\perp$, every term is $1$. Otherwise $w\cdot z_0\neq0$ for some $z_0\in C$, and because $C$ is closed
under scalars, $w\cdot(az_0)$ takes every value of $\GF(q)$, so $z\mapsto\psi(w\cdot z)$ is a nontrivial character of the group $C$ and the sum vanishes. The rest follows from F3.
\end{proof}

\begin{lemma}\label{lem:mds}
Let $D$ be a linear $[k,k-s,s+1]$ MDS code over $\GF(q)$ with $k\le q+1$. Then $D$ has codewords of every
weight $j$ with $s+1\le j\le k$, except weight $k$ when $(k,s)=(q+1,q-1)$.
\end{lemma}

\begin{proof}
By the weight formula for MDS codes, the number of words of weight $j$ is $A_j=\binom kj(q-1)\sum_{i=0}^{m}(-1)^i\binom{j-1}{i}q^{m-i}$ with $m=j-s-1$. The ratio of the
$(i+1)$st to the $i$th term of the sum is $\frac{j-1-i}{(i+1)q}\le\frac{j-1}{q}\le1$. If $j\le q$, all ratios are less than $1$ and the
alternating sum is positive. If $j=q+1$, the first two terms cancel and the remaining terms again decrease strictly,
so the sum is positive if $m\ge2$ or $m=0$, and zero if $m=1$. Thus $A_j=0$ exactly when $j=k=q+1$ and $s=q-1$. (This agrees with Ezerman, Grassl, and Sol\'e~\cite{EGS}.)
\end{proof}

The exceptional case needs a second fact. A code $C\subseteq[q]^{q+1}$ with $q^{q-1}$ words and minimum distance $3$ is \emph{perfect}: the balls of radius $1$ around its words,
each of size $1+(q+1)(q-1)=q^2$, exactly tile $[q]^{q+1}$ (because $q^{q-1}\cdot q^2=q^{q+1}$). An $\mathrm{OA}(q^{q-1},q+1,q,q-1)$ is such a code.

\begin{lemma}\label{lem:perfect}
Every $\mathrm{OA}(q^{q-1},q+1,q,q-1)$, for any labeling of the symbols by $\GF(q)$, is equidistributed along
every direction of full weight $q+1$.
\end{lemma}

\begin{proof}
Such an array $C$ has minimum distance $3$ (two rows agreeing in $q-1$ coordinates would be equal), and it is perfect, as computed. Perfectness says that every word outside $C$ has exactly one neighbor in $C$
(a word at Hamming distance $1$), and no word of $C$ has a neighbor in $C$. So the adjacency matrix $A$ of the Hamming graph $H(q+1,q)$ satisfies $A1_C=\mathbf 1-1_C$\,,
and $A\mathbf 1=(q+1)(q-1)\mathbf 1$. Hence $A(1_C-q^{-2}\mathbf 1)=-(1_C-q^{-2}\mathbf 1)$: the centered indicator is an eigenvector for the eigenvalue $-1$. That
eigenspace is the level $U_r$ with $(q+1)(q-1)-qr=-1$, that is $U_q$\,, whatever the labeling. A function
$z\mapsto g(w\cdot z)-\bar g$ with $w$ of full weight lies in the top level $U_{q+1}=L_{[q+1]}$\,, which is orthogonal to $U_q$\,.
\end{proof}

\begin{theorem}\label{thm:assign}
Let $q$ be a prime power and $1\le s<k\le q+1$. Then
\[
\dim(P^{(k,s)}_q)=\begin{cases} L^{(k,s)}_q-(q-1)^{q+1}, & \text{if } (k,s)=(q+1,q-1);\\ L^{(k,s)}_q\,, & \text{otherwise}.\end{cases}
\]
\end{theorem}

\begin{proof}
Let $[w]$ be a direction of weight $j\ge s+1$. Take a generalized Reed--Solomon $[k,k-s,s+1]$ code $D$, which exists because $k\le q+1$, with a codeword of weight $j$ as in
Lemma~\ref{lem:mds}. A monomial map, permuting and scaling the coordinates, sends $D$ to an MDS code $D'$ that contains $w$. The dual $C=D'^\perp$ is a linear
$\mathrm{OA}(q^s,k,q,s)$, hence a member of $X$, and by Lemma~\ref{lem:linear} it witnesses $[w]$. So every direction of weight at least $s+1$ is witnessed, which gives
the dimension $L^{(k,s)}_q$\,. The only exception is $(k,s)=(q+1,q-1)$: there Lemma~\ref{lem:mds} gives no witness of full weight $q+1$ among linear arrays, and
Lemma~\ref{lem:perfect} establishes that no array at all witnesses one. There are $(q-1)^{q}$ directions of full weight, each contributing $q-1$, which accounts for the deficit $(q-1)^{q+1}$.
\end{proof}

\begin{example}[Latin squares of order $3$]\label{ex:ls3}
Let $q=3$ and $(k,s)=(3,2)$, so $X$ consists of the $12$ Latin squares of order $3$, and $L=27-(1+3\cdot2+3\cdot4)=8$. Strength $2$ excludes the directions of weight $1$ and $2$; the
directions of weight $3$ are the $2^3/2=4$ classes $[(\pm1,\pm1,1)]$. The cyclic Latin square $C=\{(a,b,a+b)\}$ is the linear code with $C^\perp$ spanned by $(1,1,-1)$, a word of weight $3$, so $C$
witnesses $[(1,1,-1)]$; its monomial images witness the other three. All four directions are witnessed, and $\dim=4\cdot(3-1)=8=L$. (Numerically, every one of the 12
squares has zero component in the levels $U_1$ and $U_2$ and a nonzero component in $U_3$\,, which has dimension $8$.)
\end{example}

\begin{figure}[t]
\centering
\begin{tikzpicture}[font=\tiny,tile/.style={draw=black!45,line width=0.3pt,minimum width=0.57cm,minimum height=0.38cm,inner sep=0pt,font=\tiny\ttfamily}]
\node[anchor=east,font=\scriptsize] at (-0.15,0.00) {weight 1};
\node[tile,fill=black!15] at (0.300,0.00) {1000};
\node[tile,fill=black!15] at (1.120,0.00) {0100};
\node[tile,fill=black!15] at (1.940,0.00) {0010};
\node[tile,fill=black!15] at (2.760,0.00) {0001};
\node[anchor=east,font=\scriptsize] at (-0.15,-0.50) {weight 2};
\node[tile,fill=black!15] at (0.300,-0.50) {1100};
\node[tile,fill=black!15] at (0.900,-0.50) {1200};
\node[tile,fill=black!15] at (1.720,-0.50) {1010};
\node[tile,fill=black!15] at (2.320,-0.50) {1020};
\node[tile,fill=black!15] at (3.140,-0.50) {1001};
\node[tile,fill=black!15] at (3.740,-0.50) {1002};
\node[tile,fill=black!15] at (4.560,-0.50) {0110};
\node[tile,fill=black!15] at (5.160,-0.50) {0120};
\node[tile,fill=black!15] at (5.980,-0.50) {0101};
\node[tile,fill=black!15] at (6.580,-0.50) {0102};
\node[tile,fill=black!15] at (7.400,-0.50) {0011};
\node[tile,fill=black!15] at (8.000,-0.50) {0012};
\node[anchor=east,font=\scriptsize] at (-0.15,-1.00) {weight 3};
\node[tile,fill=blue!30] at (0.300,-1.00) {1110};
\node[tile,fill=blue!30] at (0.900,-1.00) {1120};
\node[tile,fill=blue!30] at (1.500,-1.00) {1210};
\node[tile,fill=blue!30] at (2.100,-1.00) {1220};
\node[tile,fill=blue!30] at (2.920,-1.00) {1101};
\node[tile,fill=blue!30] at (3.520,-1.00) {1102};
\node[tile,fill=blue!30] at (4.120,-1.00) {1201};
\node[tile,fill=blue!30] at (4.720,-1.00) {1202};
\node[tile,fill=blue!30] at (5.540,-1.00) {1011};
\node[tile,fill=blue!30] at (6.140,-1.00) {1012};
\node[tile,fill=blue!30] at (6.740,-1.00) {1021};
\node[tile,fill=blue!30] at (7.340,-1.00) {1022};
\node[tile,fill=blue!30] at (8.160,-1.00) {0111};
\node[tile,fill=blue!30] at (8.760,-1.00) {0112};
\node[tile,fill=blue!30] at (9.360,-1.00) {0121};
\node[tile,fill=blue!30] at (9.960,-1.00) {0122};
\node[anchor=east,font=\scriptsize] at (-0.15,-1.50) {weight 4};
\node[tile,fill=red!40] at (0.300,-1.50) {1111};
\node[tile,fill=red!40] at (0.900,-1.50) {1112};
\node[tile,fill=red!40] at (1.500,-1.50) {1121};
\node[tile,fill=red!40] at (2.100,-1.50) {1122};
\node[tile,fill=red!40] at (2.700,-1.50) {1211};
\node[tile,fill=red!40] at (3.300,-1.50) {1212};
\node[tile,fill=red!40] at (3.900,-1.50) {1221};
\node[tile,fill=red!40] at (4.500,-1.50) {1222};
\node[tile,fill=black!15,minimum width=0.3cm] at (0.15,-2.15) {}; \node[anchor=west,font=\scriptsize] at (0.35,-2.15) {forced by strength $2$};
\node[tile,fill=blue!30,minimum width=0.3cm] at (3.55,-2.15) {}; \node[anchor=west,font=\scriptsize] at (3.75,-2.15) {witnessed by linear codes};
\node[tile,fill=red!40,minimum width=0.3cm] at (7.45,-2.15) {}; \node[anchor=west,font=\scriptsize] at (7.65,-2.15) {never witnessed};
\end{tikzpicture}
\caption{The $40$ directions of $\GF(3)^4$ for pairs of orthogonal Latin squares of order $3$, $(k,s,q)=(4,2,3)$, each labeled by the representative $w$ with first nonzero coordinate $1$ and grouped by support; each carries $q-1=2$ dimensions. The relaxation leaves the $24$ directions of weight $3$ and $4$ free ($L=48$); linear codes witness all $16$ of weight $3$, and no array witnesses the $8$ of full weight, so $\dim(P)=32$.}
\label{fig:dirs}
\end{figure}
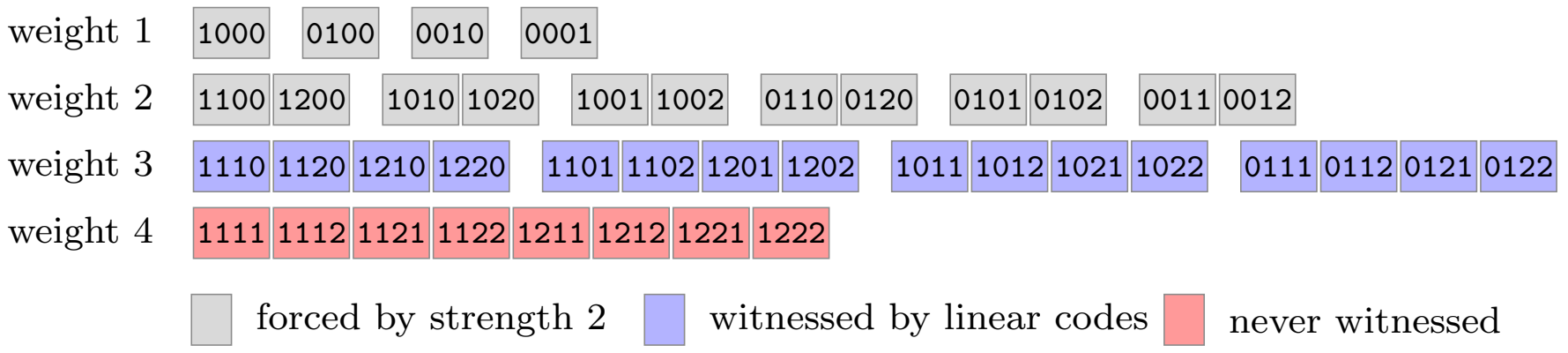
\begin{example}[the exceptional case $q=3$, $(k,s)=(4,2)$]\label{ex:tetra}
Now $X$ consists of the $72$ arrays $\mathrm{OA}(9,4,3,2)$ (pairs of orthogonal Latin squares of order $3$), and $L=81-(1+4\cdot2+6\cdot4)=48$. A linear member is
the \emph{tetracode} $C=\{(a,b,a+b,a+2b)\}$; it is self-dual, an MDS code $[4,2,3]$, whose $9$ words have weights $0$ and $3$ (eight of them) but none of weight $4$. So it witnesses the
directions of weight $3$ only. By Lemma~\ref{lem:perfect}, the same holds for every array in $X$ (Figure~\ref{fig:dirs}): the $2^4/2=8$ directions of weight $4$ are never witnessed, so
the dimension is $48-8\cdot2=32$, not $48$. The linear relaxation therefore has dimension $48$, but the polytope has dimension $32$: every pair of orthogonal Latin squares of order $3$ satisfies $16$ further
independent linear equations, namely, for each of the eight full-weight directions $w$, the three hyperplane sums $\sum_{w\cdot z=c}x_z$\,, $c\in\GF(3)$, are equal.
\end{example}

\begin{remark}
Theorem~\ref{thm:assign} recovers, for prime powers, the known values $\dim(P^{(2,1)})=(n-1)^2$ (see~\cite{BR}),
$\dim(P^{(3,1)})=n^3-3n+2$ for $n\ge3$ (see~\cite{BS}), $\dim(P^{(3,2)})=(n-1)^3$ for $n\ge4$ (see~\cite{EBG}), which the theorem also gives for $q=2,3$, and
$\dim(P^{(4,2)})=n^4-6n^2+8n-3$ for $n\ge4$ and $n\neq6$ (see~\cite{AMMJ}). It also explains the excluded small orders: $(3,1)$ at
$q=2$ and $(4,2)$ at $q=3$ are the first two exceptional cases, with dimensions $3$ and $32$. The next ones,
$(5,3)$ at $q=4$ (dimension $405$ instead of $648$), $(6,4)$ at $q=5$, and so on, demonstrate that the expected equality
fails for $(k,s,n)=(q+1,q-1,q)$ for every prime power $q$; for $q=2,3$ these orders lie outside the ranges of the
known theorems. The deficit is exactly the
one allowed by Bulutoglu's Theorem~6~\cite{Bul} (with $\ell=2$, because $\binom{s+1}{1}=q\equiv0\pmod q$).
\end{remark}

\subsection{Relation to the literature}\label{sec:lit}

\emph{Known cases.} Appa, Magos, and Mourtos~\cite{AMM} determine the dimension of the $(k,1)$ (``axial'')
polytopes for $n\ge3$ (their Theorem~9) and of the nonempty $(k,2)$ (``planar'') polytopes for $n\ge\max\{5,k\}$ (their
Theorem~17), for all orders, not only prime powers. For $s\le2$, Theorem~\ref{thm:assign} also covers the prime-power orders outside these
ranges, which include the case $k=n+1$ (for instance $(5,2)$ at $q=4$ and $(6,2)$ at $q=5$)
and the orders $n=2$ and $n=3$.

\emph{Bulutoglu's criterion.} Bulutoglu~\cite{Bul} demonstrates by representation theory that the dimension can fall short of
$L^{(k,s)}_n$ only by whole levels $U_{s+\ell}$ with $n\mid\binom{s+\ell-1}{\ell-1}$, so the polytope has dimension
$L^{(k,s)}_n$ whenever no such $\ell\le k-s$ exists (his Corollary~2). Over prime powers $q\le32$ with $k\le q+1$, this criterion
settles $2994$ of the $3217$ cases $(k,s,q)$; Theorem~\ref{thm:assign} settles all of them, and in the exceptional
family it demonstrates that the deficit his theorem allows, the whole top level $U_{q+1}$\,, is actually attained. His symmetry
group $S_n\wr S_k$ sees only whole levels (Proposition~\ref{prop:hamming}); the affine group sees the direction pieces inside them, which is what
makes the witness argument possible. As of his paper, no counterexample of index
one was known; he points to $(n,k,s)=(10,6,2)$ as the smallest open candidate.

\emph{Ingredients.} The weights of MDS codes (Lemma~\ref{lem:mds}) are classical; see~\cite{EGS}. That perfect codes
are completely regular, with the indicator in a single eigenspace of the Hamming graph, is due to Delsarte. The new
point is the use of both, through Theorem~\ref{thm:affine}, as witnesses and obstructions for the dimension.

\emph{The expectation in~\cite{AMM}.} It is not stated as a formal conjecture, and it carries no restriction on $n$; it
is Bulutoglu~\cite{Bul} who calls it a conjecture. Because it covers all $(k,s,n)$ with $P^{(k,s)}_n\neq\emptyset$, every
exceptional case is a counterexample, although for $q=2,3$ these orders are outside the ranges of Theorems~9 and~17
of~\cite{AMM}.

\subsection{Longer arrays and non-prime-power orders}\label{sec:longer}

By Proposition~\ref{prop:hamming}, $\dim(P^{(k,s)}_n)=\sum\{\dim(U_r) ~:~ r>s \text{ witnessed}\}$. Coordinatewise symbol permutations preserve
every level, so which levels are witnessed depends only on the equivalence classes of the arrays. For a linear code $C$ over $\GF(q)$, Lemma~\ref{lem:linear} says that
level $r$ is witnessed by $C$ exactly when $C^\perp$ has a word of weight $r$; this is the argument of Theorem~\ref{thm:assign} in the language of
levels. The two lemmas of the previous subsection also extend: the proof of Lemma~\ref{lem:perfect} uses only the Hamming graph, so it
holds for $\mathrm{OA}(n^{n-1},n+1,n,n-1)$ over any alphabet of size~$n$, and it has the following consequence one
step further out.

\begin{lemma}\label{lem:puncture}
The indicator of every $\mathrm{OA}(n^{n-1},n+2,n,n-1)$ has zero component at level $n+1$.
\end{lemma}

\begin{proof}
Deleting coordinate $i$ maps such an array $C$ bijectively onto an $\mathrm{OA}(n^{n-1},n+1,n,n-1)$, which we call $C_i$\,, whose
indicator lies in $U_0\oplus U_n$ by (the extension of) Lemma~\ref{lem:perfect}. For $u\subseteq[n+2]\setminus\{i\}$, the
$L_u$-component of $1_C$ depends only on the marginal of $1_C$ on the coordinates other than $i$, which is $1_{C_i}$; so it
vanishes unless $|u|\in\{0,n\}$. Every $u$ with $|u|=n+1$ omits exactly one coordinate.
\end{proof}

\begin{theorem}\label{thm:qplus2}
Let $q\ge4$ be a power of~$2$, and $k\coloneqq q+2$.
\begin{enumerate}
\item[(i)] If $q\ge8$, then $\dim(P^{(q+2,3)}_q)=L^{(q+2,3)}_q$\,.
\item[(ii)] $\dim(P^{(q+2,q-1)}_q)=L^{(q+2,q-1)}_q-(q+2)(q-1)^{q+1}$.
\end{enumerate}
For $q=4$ the two cases coincide, and the dimension is $3402-1458=1944$.
\end{theorem}

\begin{proof}
A hyperoval of $\PG(2,q)$ (a set of $q+2$ points, no three collinear, which exists for even $q$) gives a linear $[q+2,3,q]$ MDS code $H$, whose dual is a $[q+2,q-1,4]$ MDS code; so
both polytopes are nonempty. As in the proof of Theorem~\ref{thm:assign}, a level $r$ is witnessed as soon as some
linear MDS code of the dual parameters has a word of weight $r$, because monomial images of MDS codes are MDS.
(i) The dual code is $H^\perp$, of dimension $q-1$. For $4\le j\le q+1$ the argument of Lemma~\ref{lem:mds} gives
$A_j>0$. For $j=q+2$ the alternating sum is the binomial expansion of $(q-1)^{q+1}$ without its last three terms,
which gives $A_{q+2}=(q-1)q^{-3}\bigl((q-1)^{q+1}+\binom{q+1}2q^2-(q+1)q+1\bigr)>0$. So every level $r\ge4$ is witnessed.
(ii) The dual code is $H$. A linear form vanishes at $0$ or $2$ points of the hyperoval, so $H$ has words of weights
$q+2$ and $q$, and none of weight $q+1$. So the levels $q$ and $q+2$ are witnessed, and the level $q+1$ is not, by
Lemma~\ref{lem:puncture}; it has dimension $\binom{q+2}{q+1}(q-1)^{q+1}$.
\end{proof}

So $(q+2,q-1,q)$, for every power $q\ge4$ of $2$, is a second family of counterexamples to the expectation
of~\cite{AMM}, alongside $(q+1,q-1,q)$.

\emph{Orders that are not prime powers.} If $C_1\in P^{(k,s)}_{n_1}$ and $C_2\in P^{(k,s)}_{n_2}$\,, the rows obtained by combining every row of $C_1$ with
every row of $C_2$ coordinatewise form an array $C$ in $P^{(k,s)}_{n_1n_2}$\,, with $1_C=1_{C_1}\otimes1_{C_2}$ under $[n_1n_2]^k=[n_1]^k\times[n_2]^k$. For $u\subseteq[k]$, the space
$L_u(n_1)\otimes L_\emptyset(n_2)$ lies in $L_u(n_1n_2)$, and the component of $1_C$ there is the $L_u$-component of $1_{C_1}$ times the
nonzero constant $|C_2|/n_2^k\,$. So every level witnessed in $P^{(k,s)}_{n_1}$ (or in $P^{(k,s)}_{n_2}$) is witnessed in $P^{(k,s)}_{n_1n_2}$\,.

\begin{proposition}\label{prop:products}
Let $1\le s<k$, and let $n\coloneqq q_1\cdots q_m$ with $m\ge2$ and the $q_i$ pairwise coprime prime powers. Suppose that for each $i$ either
$k\le q_i+1$, or $q_i$ is even, $k=q_i+2$, and $s\in\{3,q_i-1\}$. Then $\dim(P^{(k,s)}_n)=L^{(k,s)}_n$\,.
\end{proposition}

\begin{proof}
Each $P^{(k,s)}_{q_i}$ is nonempty, so products give members of $P^{(k,s)}_n$\,. By Theorems~\ref{thm:assign} and~\ref{thm:qplus2}
and their proofs, $P^{(k,s)}_{q_i}$ witnesses every level $r$ with $s<r\le k$ unless $q_i=s+1$, which holds for at most one $i$. So
every such level is witnessed in some factor, hence in $P^{(k,s)}_n$\,, and by Proposition~\ref{prop:hamming} the dimension is
$\sum_{r>s}\dim(U_r)=L^{(k,s)}_n$\,.
\end{proof}

With Theorem~17 of~\cite{AMM} and Bulutoglu's criterion, this leaves only orders and parameters at which
orthogonal arrays exist that are not products of linear MDS codes, a question about existence rather than about
dimension.

\section{Transversals of Latin squares}\label{sec:latin}

\begin{figure}[t]
\centering
\begin{tikzpicture}[scale=0.42,every node/.style={font=\scriptsize}]
\fill[red!35] (0,7) rectangle ++(1,1);\node at (0.5,7.5) {00};
\node at (1.5,7.5) {01};
\node at (2.5,7.5) {02};
\node at (3.5,7.5) {03};
\node at (4.5,7.5) {10};
\node at (5.5,7.5) {11};
\node at (6.5,7.5) {12};
\node at (7.5,7.5) {13};
\node at (0.5,6.5) {01};
\fill[yellow!80!black!35] (1,6) rectangle ++(1,1);\node at (1.5,6.5) {02};
\node at (2.5,6.5) {03};
\node at (3.5,6.5) {00};
\node at (4.5,6.5) {11};
\node at (5.5,6.5) {12};
\node at (6.5,6.5) {13};
\node at (7.5,6.5) {10};
\node at (0.5,5.5) {02};
\node at (1.5,5.5) {03};
\node at (2.5,5.5) {00};
\node at (3.5,5.5) {01};
\fill[violet!35] (4,5) rectangle ++(1,1);\node at (4.5,5.5) {12};
\node at (5.5,5.5) {13};
\node at (6.5,5.5) {10};
\node at (7.5,5.5) {11};
\node at (0.5,4.5) {03};
\node at (1.5,4.5) {00};
\node at (2.5,4.5) {01};
\node at (3.5,4.5) {02};
\node at (4.5,4.5) {13};
\fill[cyan!80!black!35] (5,4) rectangle ++(1,1);\node at (5.5,4.5) {10};
\node at (6.5,4.5) {11};
\node at (7.5,4.5) {12};
\node at (0.5,3.5) {10};
\node at (1.5,3.5) {11};
\node at (2.5,3.5) {12};
\fill[magenta!35] (3,3) rectangle ++(1,1);\node at (3.5,3.5) {13};
\node at (4.5,3.5) {00};
\node at (5.5,3.5) {01};
\node at (6.5,3.5) {02};
\node at (7.5,3.5) {03};
\node at (0.5,2.5) {11};
\node at (1.5,2.5) {12};
\node at (2.5,2.5) {13};
\node at (3.5,2.5) {10};
\node at (4.5,2.5) {01};
\node at (5.5,2.5) {02};
\fill[green!70!black!35] (6,2) rectangle ++(1,1);\node at (6.5,2.5) {03};
\node at (7.5,2.5) {00};
\node at (0.5,1.5) {12};
\node at (1.5,1.5) {13};
\node at (2.5,1.5) {10};
\node at (3.5,1.5) {11};
\node at (4.5,1.5) {02};
\node at (5.5,1.5) {03};
\node at (6.5,1.5) {00};
\fill[orange!35] (7,1) rectangle ++(1,1);\node at (7.5,1.5) {01};
\node at (0.5,0.5) {13};
\node at (1.5,0.5) {10};
\fill[blue!35] (2,0) rectangle ++(1,1);\node at (2.5,0.5) {11};
\node at (3.5,0.5) {12};
\node at (4.5,0.5) {03};
\node at (5.5,0.5) {00};
\node at (6.5,0.5) {01};
\node at (7.5,0.5) {02};
\draw[black!50] (0,0) grid (8,8);
\node at (-0.6,7.5) {00}; \node at (0.5,8.5) {00};
\node at (-0.6,6.5) {01}; \node at (1.5,8.5) {01};
\node at (-0.6,5.5) {02}; \node at (2.5,8.5) {02};
\node at (-0.6,4.5) {03}; \node at (3.5,8.5) {03};
\node at (-0.6,3.5) {10}; \node at (4.5,8.5) {10};
\node at (-0.6,2.5) {11}; \node at (5.5,8.5) {11};
\node at (-0.6,1.5) {12}; \node at (6.5,8.5) {12};
\node at (-0.6,0.5) {13}; \node at (7.5,8.5) {13};
\draw[-{Stealth},thick] (9.2,4) -- (10.4,4);
\begin{scope}[xshift=11.5cm]
\node[circle,fill,inner sep=1.3pt,label=left:{00}] (r0) at (0,7.50) {}; \node[circle,fill,inner sep=1.3pt,label=right:{00}] (c0) at (3.2,7.50) {};
\node[circle,fill,inner sep=1.3pt,label=left:{01}] (r1) at (0,6.50) {}; \node[circle,fill,inner sep=1.3pt,label=right:{01}] (c1) at (3.2,6.50) {};
\node[circle,fill,inner sep=1.3pt,label=left:{02}] (r2) at (0,5.50) {}; \node[circle,fill,inner sep=1.3pt,label=right:{02}] (c2) at (3.2,5.50) {};
\node[circle,fill,inner sep=1.3pt,label=left:{03}] (r3) at (0,4.50) {}; \node[circle,fill,inner sep=1.3pt,label=right:{03}] (c3) at (3.2,4.50) {};
\node[circle,fill,inner sep=1.3pt,label=left:{10}] (r4) at (0,3.50) {}; \node[circle,fill,inner sep=1.3pt,label=right:{10}] (c4) at (3.2,3.50) {};
\node[circle,fill,inner sep=1.3pt,label=left:{11}] (r5) at (0,2.50) {}; \node[circle,fill,inner sep=1.3pt,label=right:{11}] (c5) at (3.2,2.50) {};
\node[circle,fill,inner sep=1.3pt,label=left:{12}] (r6) at (0,1.50) {}; \node[circle,fill,inner sep=1.3pt,label=right:{12}] (c6) at (3.2,1.50) {};
\node[circle,fill,inner sep=1.3pt,label=left:{13}] (r7) at (0,0.50) {}; \node[circle,fill,inner sep=1.3pt,label=right:{13}] (c7) at (3.2,0.50) {};
\draw[red,line width=1.1pt] (r0)--(c0);
\draw[yellow!80!black,line width=1.1pt] (r1)--(c1);
\draw[violet,line width=1.1pt] (r2)--(c4);
\draw[cyan!80!black,line width=1.1pt] (r3)--(c5);
\draw[magenta,line width=1.1pt] (r4)--(c3);
\draw[green!70!black,line width=1.1pt] (r5)--(c6);
\draw[orange,line width=1.1pt] (r6)--(c7);
\draw[blue,line width=1.1pt] (r7)--(c2);
\end{scope}
\end{tikzpicture}
\caption{The Cayley table of $\mathbb Z_2\times\mathbb Z_4$ (elements written $gh$) with a transversal, shaded by symbol: the graph of a complete mapping $\sigma$, the cell $(a,\sigma(a))$ carrying the symbol $a+\sigma(a)$; equivalently, a rainbow perfect matching of $K_{8,8}$\,.}
\label{fig:cayley}
\end{figure}
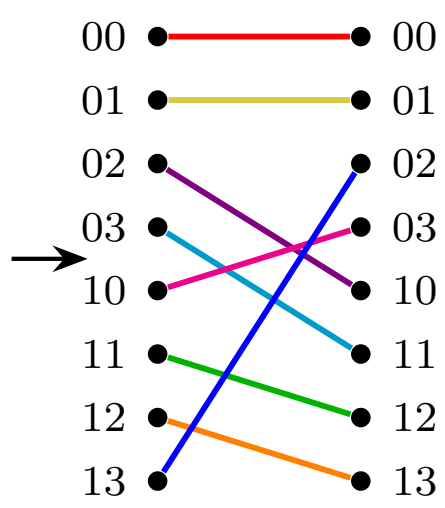

A \emph{Latin square} of order $n$ is an $n\times n$ array of symbols from $[n]$ in which every symbol occurs once in each row and once in
each column. A \emph{transversal} is a set of $n$ cells, one in each row, one in each column, and carrying each symbol once (Figure~\ref{fig:cayley}). Let
$\mathcal T(L)\subseteq\R^{[n]\times[n]}$ be the \emph{transversal polytope}, the convex hull of the incidence vectors of the transversals of $L$.
The $3n$ line-sum equations (every row, every column, and every symbol class has sum $1$) have rank $3n-2$, because the three families have the same
total; so $\dim(\mathcal T(L))\le n^2-(3n-2)=(n-1)(n-2)$. \S\ref{sec:nets} established that for the Cayley table of $(\GF(q),+)$ the bound is attained
(Corollary~\ref{cor:nets} with $k=3$). Here we ask for which Latin squares it is, beginning with the Cayley tables of all finite abelian groups, where the
Fourier analysis of \S\ref{sec:fourier} still applies.

\subsection{Group tables and a Fourier criterion}\label{sec:fourier2}

Let $G$ be a finite abelian group of order $n$, written additively, and let $L(a,b)\coloneqq a+b$ be its Cayley table. A \emph{complete mapping} of $G$ is a
bijection $\sigma\colon G\to G$ such that $a\mapsto a+\sigma(a)$ is also a bijection. The transversals of $L$ are exactly the sets
$\{(a,\sigma(a)) ~:~ a\in G\}$ for complete mappings $\sigma$: the cell $(a,b)$ has symbol $a+b$. If $|G|$ is odd, the identity is a complete
mapping. By a theorem of Hall and Paige~\cite{HP}, an abelian group has a complete mapping if and only if its Sylow $2$-subgroup, the subgroup of elements whose order is a power of $2$, is trivial or not cyclic
(for example, $\mathbb Z_{12}$ has Sylow $2$-subgroup $\mathbb Z_4$ and no complete mapping, while $\mathbb Z_2\times\mathbb Z_4$ has one). We assume
this from now on.

\emph{Characters.} A \emph{character} of $G$ is a homomorphism $u$ from $G$ to the unit circle. The characters form a group under pointwise multiplication, of order $n$; we write it
additively, so $u+v$ is the character $a\mapsto u(a)v(a)$ and $0$ is the trivial character, and we identify it with $G$ when convenient. As in F1, $\sum_{a\in G}u(a)$ is
$n$ if $u=0$ and $0$ otherwise, so characters are orthogonal. For $G=\mathbb Z_n$ the characters are $a\mapsto\omega^{ua}$ with $\omega\coloneqq e^{2\pi i/n}$ and $u\in\mathbb Z_n$; the characters of $G\times G$
are the products $(a,b)\mapsto u(a)v(b)$, written $(u,v)$. For a complete mapping $\sigma$ put $\hat\sigma(u,v)\coloneqq\sum_{a\in G}u(a)v(\sigma(a))$.

\emph{Some terminology.} Beyond the fundamental theorem of finite abelian groups we use only the following. For a prime $p$, a \emph{$p$-group} is a group of $p$-power order, and the \emph{$p$-part} $G_p$ of $G$, also called its Sylow $p$-subgroup, is the subgroup of its elements of $p$-power order; $G$ is the direct sum of its $p$-parts. By the fundamental theorem a $p$-group is $P\cong\mathbb Z_{p^{e_1}}\oplus\dots\oplus\mathbb Z_{p^{e_r}}$ with all $e_i\ge1$, and the number $r$ of summands is the \emph{rank} of $P$; so $P$ is cyclic exactly when $r=1$. We call $P$ \emph{homocyclic} if all $e_i$ are equal and \emph{elementary abelian} if all are $1$, so that $P\cong\mathbb Z_p^r$ is the additive group of $\GF(p^r)$. For a $2$-group $P$ we write $2P\coloneqq\{2a:a\in P\}$, so that $P/2P\cong\mathbb Z_2^r\,$. An \emph{automorphism} is a group isomorphism $G\to G$, such as multiplication by an integer prime to $n$. A character $u$ of $G$ \emph{restricts} to a character $u_N$ of a subgroup $N$; if $u$ is trivial on $N$, it \emph{induces} the character $\bar u(a+N)\coloneqq u(a)$ of the quotient $G/N$, and every character of $G/N$ arises in this way. A \emph{section} of $G\to G/N$ is a choice of one element $s(x)$ in each coset $x$.

\begin{proposition}[Fourier criterion]\label{prop:fourier}
Suppose $G$ has a complete mapping. Call a pair $(u,v)$ of characters \emph{allowed} if $u\neq0$, $v\neq0$, and $u\neq v$; there are ${(n-1)(n-2)}$ allowed pairs. Then $\dim(\mathcal T(L))$ is the number of allowed pairs $(u,v)$ for which some complete mapping $\sigma$, a \emph{witness}, has $\hat\sigma(u,v)\neq0$. So $\dim(\mathcal T(L))=(n-1)(n-2)$ exactly when every allowed pair has a witness.
\end{proposition}

(This is the notion of \S\ref{sec:fourier}: a transversal whose Fourier coefficient at the character $(u,v)$ does not vanish, that is, a member of $X$ that is not orthogonal to the eigenspace of $(u,v)$.)

\begin{proof}
Let $X$ be the set of incidence vectors of transversals, indexed by the cells $(a,b)\in G\times G$. Translations $(a,b)\mapsto(a+x,b+y)$ shift the rows, the columns, and
the symbols, so they permute $X$. Steps~1 and~2 of the proof of Theorem~\ref{thm:affine} apply verbatim to the abelian group $G\times G$ in place of $(V,+)$: the matrix
$B$ is translation invariant, the characters $(u,v)$ are its eigenvectors, and the eigenvalue at $(u,v)$ is $\frac1{n^2}\sum_{x\in X}|\hat x(u,v)|^2$ for $(u,v)\neq(0,0)$, which is nonzero exactly when some
$x\in X$ has $\hat x(u,v)\neq0$. At $(0,0)$ the eigenvalue is $0$, because the coordinate sum is $n$ on all of $X$. For $x=1_T$ with $T=\{(a,\sigma(a))\}$, $\hat x(u,v)$ is the complex conjugate of $\hat\sigma(u,v)$.
Now the pairs that are not allowed give $\hat\sigma(u,v)=0$ for every complete mapping $\sigma$: if $v=0$ and $u\neq0$, then $\hat\sigma=\sum_au(a)=0$; if $u=0$ and $v\neq0$, then $\hat\sigma=\sum_av(\sigma(a))=\sum_bv(b)=0$ because
$\sigma$ is a bijection; if $u=v\neq0$, then $\hat\sigma=\sum_au(a+\sigma(a))=0$ likewise. (These are the row, column, and symbol sums.) Hence $\operatorname{rank}(B)$ counts the allowed pairs with a witness.
\end{proof}

\begin{example}[$G=\mathbb Z_5$]\label{ex:z5}
Let $\omega\coloneqq e^{2\pi i/5}$. The map $\sigma_k(a)=ka$ is a complete mapping when $k$ and $k+1$ are units, that is, for $k\in\{1,2,3\}$, and
\[
\hat\sigma_k(u,v)=\sum_a\omega^{ua+vka}=\begin{cases}5,&u+kv=0;\\0,&\text{otherwise.}\end{cases}
\]
For an allowed pair put $k\coloneqq-u/v$: it is nonzero because $u\neq0$, and different from $-1$ because $u\neq v$, so $k\in\{1,2,3\}$. Hence every one of the $12$ allowed
pairs has a (linear) witness, and $\dim(\mathcal T(L))=12=(5-1)(5-2)$. The same argument works for every prime $n$.
\end{example}

For $\mathbb Z_9$\,, linear maps are no longer enough. For the pair $(u,v)=(3,1)$, $\hat\sigma_k(u,v)\neq0$ requires $3+k\equiv0$, that is $k=6$, which is not a unit modulo $9$, so $\sigma_k$ is not a
complete mapping. We will nevertheless find witnesses for all allowed pairs of every abelian group that has complete mappings (Theorem~\ref{thm:main}),
by building nonlinear complete mappings from a few constructions.

\subsection{Overview of the proof of Theorem~\ref{thm:main}}\label{sec:overview}
\begin{figure}[t]
\centering
\thickmuskip=5mu\relax\medmuskip=4mu\relax
\begin{tikzpicture}[font=\scriptsize, box/.style={draw,rounded corners,align=center,inner sep=2.5pt}, >={Stealth}]
\node[box,fill=blue!6,text width=5.7cm,align=center] (odd) at (0,0) {\textbf{odd $p$-groups} (Lemmas~\ref{lem:bigp}, \ref{lem:three}, \ref{lem:pgroups})\\ $p\ge5$: switches with multipliers of order $p-1$\\ $p=3$: explicit nonlinear family\\ then liftings from cyclic subgroups};
\node[box,fill=red!6,text width=5.85cm,align=center,right=0.3cm of odd] (two) {\textbf{noncyclic $2$-groups} (Theorem~\ref{thm:twogroups})\\ $(\mathbb Z_{2^e})^m$: switches over Galois rings\\ rank $2$: explicit family, three-cell modifications,\\ cyclic-quotient design\\ rank $\ge3$: liftings, $\GF(4)$ design};
\node[box,fill=black!5,text width=11.85cm,align=center,below=0.3cm of $(odd.south)!0.5!(two.south)$] (comb) {\textbf{combine over the Sylow components} \mbox{$G=\bigoplus_pG_p$}\\ restrictions of $(u,v)$ all allowed or zero: product of witnesses\\ some excluded, some allowed: twisted product\\ two excluded of different types: layered coupling ($r\ge5$, quadratic residues) or the $\mathbb Z_3$ design};
\node[box,fill=green!8,below=0.3cm of comb] (res) {a witness for every allowed pair $(u,v)$ $\Longrightarrow$ $\dim(\mathcal T(L))=(n-1)(n-2)$};
\draw[->] (odd.south) -- (odd.south|-comb.north);
\draw[->] (two.south) -- (two.south|-comb.north);
\draw[->] (comb) -- (res);
\end{tikzpicture}
\caption{Structure of the proof of Theorem~\ref{thm:main}: witnesses for each component type, assembled by the symmetries and three combination rules.}
\label{fig:arch}
\end{figure}
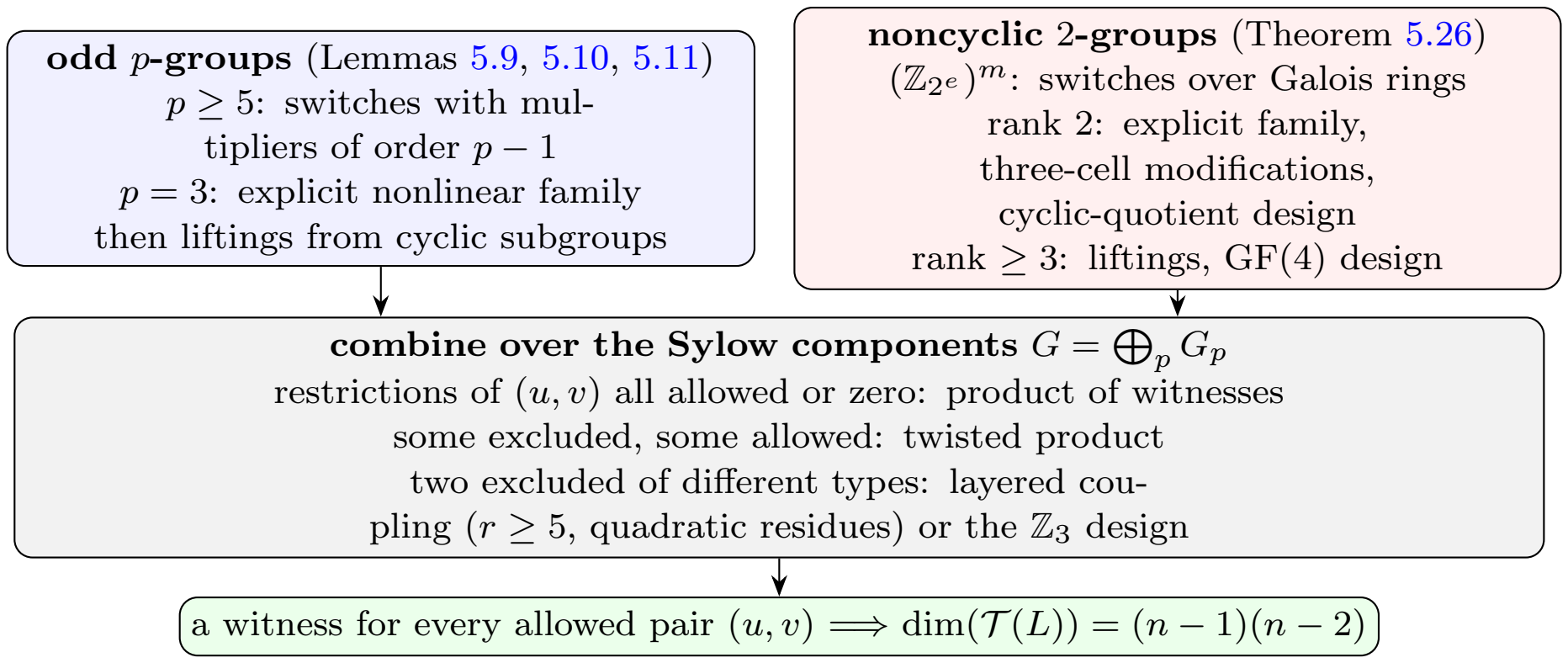
\paragraph{The constructions.}
Three \emph{symmetries} of the set of complete mappings turn witnesses into witnesses: conjugation by automorphisms of $G$, the inverse $\sigma^{-1}$ (swapping $u$ and $v$), and $x\mapsto-x+\sigma(-x)$ (mapping $(u,v)$ to $(v-u,v)$); together they permute the three types of excluded pairs. A \emph{lifting} assembles a witness on $G$ coset by coset from one on a subgroup or on a quotient (Lemma~\ref{lem:lift}). A \emph{switch} modifies the linear map $\sigma_k$ on a single orbit of a unit acting by multiplication, and is a witness as soon as that orbit separates two auxiliary characters (Lemma~\ref{lem:switch}); with units of order $p-1$ from $\mathbb Z_{p^f}$ this covers the cyclic $p$-groups for $p\ge5$, and with units of a Galois ring~\cite{Wan} the $2$-groups $(\mathbb Z_{2^e})^m$, $m\ge2$ (Lemmas~\ref{lem:bigp} and~\ref{lem:homocyclic}). Cyclic $3$-groups have no suitable units and need an explicit nonlinear family (Lemma~\ref{lem:three}). Finally, two \emph{couplings} combine complete mappings on two components: a twisted product (Lemma~\ref{lem:twist}) when one of the two restrictions is allowed, and a layered coupling built on quadratic residues when neither is (Lemma~\ref{lem:two}, and Proposition~\ref{prop:coupling3} for a separate design over $\mathbb Z_3$).

\paragraph{The proof (Figure~\ref{fig:arch}).}
Write $G=\bigoplus_pG_p$ with $G_p$ its Sylow $p$-subgroups. Every allowed pair on an odd $p$-group has a witness, by liftings from a cyclic subgroup on which the pair restricts to an allowed one (Lemma~\ref{lem:pgroups}); every allowed pair on a noncyclic $2$-group has one (Theorem~\ref{thm:twogroups}, below). For $(u,v)$ allowed on $G$, each restriction to a component is zero, allowed, or excluded of one of three types. If none is excluded, a product of witnesses and complete mappings is a witness, because $\hat\sigma$ is multiplicative; if some are excluded and some allowed, a twisted product is; otherwise every nonzero restriction is excluded, not all of the same type (else $(u,v)$ would be excluded), and a layered coupling on two components of different types makes the coefficient nonzero.

\paragraph{The $2$-groups.}
For $2$-groups the odd tools fail, and the hardest cases are the groups $\mathbb Z_{2^a}\times\mathbb Z_{2^b}$ of rank $2$, where liftings cover only some of the pairs. There we use an explicit family of complete mappings whose Fourier coefficients factor into three elementary sums, so that the pairs it witnesses are known exactly (Lemma~\ref{lem:family}), modified at three cells when needed (Lemmas~\ref{lem:cycle} and~\ref{lem:shape}). One class of pairs, the \emph{cyclic-quotient pairs} $(u,2^{j-1}u)$ with $u$ of order $2^j\ge4$, is invisible to all such constructions and needs a design of its own (Lemma~\ref{lem:cq}); a second design over $\GF(4)$ handles two residual cases of rank $3$ (Lemma~\ref{lem:gf4d}). Groups of rank at least $3$ reduce to these by liftings, and an induction on the order completes the proof of Theorem~\ref{thm:twogroups} (Propositions~\ref{prop:z2zb}--\ref{prop:rank2}).

\subsection{Tools for building witnesses}\label{sec:tools}

\emph{Symmetries.} Three operations turn witnesses into witnesses.
(S1) If $\alpha$ is an automorphism of $G$ (for example $a\mapsto ka$ for a multiplier $k$) and $\sigma$ is a witness for $(u,v)$, then $\alpha\sigma\alpha^{-1}$ is a witness for $(u\circ\alpha^{-1},v\circ\alpha^{-1})$.
(S2) The inverse $\sigma^{-1}$ is a complete mapping, and $\widehat{\sigma^{-1}}(v,u)=\hat\sigma(u,v)$, so a witness for $(u,v)$ gives one for $(v,u)$.
(S3) Exchanging the roles of columns and symbols: the map $x\mapsto-x+\sigma(-x)$ is a complete mapping $\sigma'$ with $\hat\sigma'(u',v')=\hat\sigma(v'-u',v')$, so a witness for $(u,v)$ gives one for $(v-u,v)$.
Call a pair \emph{of type rows, columns, or symbols} if it is not allowed because $v=0$, $u=0$, or $u=v$, respectively. (S2) and (S3) permute the three types arbitrarily.

\emph{Multipliers and linear maps.} A \emph{multiplier} is an integer $k$ coprime to $n$; more generally, when $G$ is the additive group of a finite commutative ring $R$ (such as $\mathbb Z_n$\,, or a Galois ring in \S\ref{sec:even}), a multiplier is a unit of $R$, and $ka$ denotes the ring product.
A multiplier acts on $G$ by $a\mapsto ka$ and on characters by $(kw)(a)\coloneqq w(ka)$. If $k$ and $k+1$ are multipliers, then $\sigma_k(a)=ka$ is a complete mapping, and $\hat\sigma_k(u,v)=|G|$ if $u+kv=0$ and $0$ otherwise (Example~\ref{ex:z5}).

\emph{Lifting.} Witnesses on a subgroup or a quotient lift to $G$, coset by coset.

\begin{lemma}[lifting]\label{lem:lift}
Let $N\le G$ be a subgroup such that $G/N$ has a complete mapping $\varphi$.
(i) If $N$ also has a complete mapping (automatic when $|N|$ is odd), $u$ and $v$ are trivial on $N$, and $\varphi$ is a witness for the induced pair on $G/N$, then $(u,v)$ has a witness.
(ii) If $(u_N\,,v_N)$, the restriction to $N$, is allowed and has a witness on $N$, then $(u,v)$ has a witness.
\end{lemma}

\begin{proof}
Fix a section $s\colon G/N\to G$, a complete mapping $\theta$ of $N$, and elements $c_x\in N$ for $x\in G/N$, and let $\sigma(s(x)+m)\coloneqq s(\varphi(x))+\theta(m)+c_x$ for $m\in N$.
This is a complete mapping of $G$: it maps the coset $s(x)+N$ bijectively onto $s(\varphi(x))+N$, and $a\mapsto a+\sigma(a)$ maps it bijectively onto $s(x+\varphi(x))+N$, a different coset
for each $x$. Its coefficient is
\[
\hat\sigma(u,v)=\hat\theta(u_N\,,v_N)\sum_{x}u(s(x))\,v(s(\varphi(x)))\,v(c_x).
\]
In case (i), the map $\theta$ is a complete mapping of $N$ (here the hypothesis is used), we have $v(c_x)=1$ and $\hat\theta=|N|$, and the sum is the coefficient of $\varphi$ for the induced pair. In case (ii), take $\theta$ a witness on $N$ (so $\hat\theta\neq0$) and
$c_x=0$; if the sum vanishes, change one $c_{x_0}$ to an element with $v(c_{x_0})\neq1$, which exists because $v_N\neq0$; the sum then equals $u(s(x_0))v(s(\varphi(x_0)))(v(c_{x_0})-1)\neq0$.
\end{proof}

\begin{example}[lifting in $\mathbb Z_{15}$]\label{ex:lift}
Let $G\coloneqq\mathbb Z_{15}$\,, $N\coloneqq\langle3\rangle=\{0,3,6,9,12\}\cong\mathbb Z_5$\,, and $(u,v)\coloneqq(5,10)$. Both characters depend only on $a$ modulo $3$, so they are trivial on $N$, and the induced pair on
$G/N\cong\mathbb Z_3$ is $(1,2)$, for which the identity $\varphi$ is a witness ($u+v=0$). Take $s(x)\coloneqq x$, $\theta(3j)\coloneqq6j$ (the map $y\mapsto2y$ on $N\cong\mathbb Z_5$), and $c_x\coloneqq0$: then
$\sigma(x+3j)=x+6j$ for $x\in\{0,1,2\}$, $j\in\mathbb Z_5$\,. One checks that $\sigma$ and $a\mapsto a+\sigma(a)=2x+9j$ are bijections of $\mathbb Z_{15}$\,, and
$\hat\sigma(5,10)=\sum_a\omega_3^{x+2x}=15$.
\end{example}

\emph{Two more constructions.} Let $\mu$ be a multiplier such that $\mu+1$ is also a multiplier, and put $k\coloneqq-\mu^{-1}-1$. Then $k$ and $k+1$ are multipliers, and $1+k\mu=(1+k)\mu^2$.

\begin{lemma}[switches]\label{lem:switch}
Let $\sigma_k(a)\coloneqq ka$, let $c\in G$, let $O\coloneqq\{\mu^ic\}$ be its orbit under multiplication by $\mu$, and define $\sigma(a)\coloneqq k\mu a$ for $a\in O$ and $\sigma(a)\coloneqq ka$ otherwise. Then $\sigma$ is a complete
mapping, and $\hat\sigma(u,v)=\hat\sigma_k(u,v)+\sum_{a\in O}(\chi_1(a)-\chi_2(a))$ with $\chi_1\coloneqq u+k\mu v$ and $\chi_2\coloneqq u+kv$. If $\chi_1$ and $\chi_2$ lie in different orbits of
multiplication by $\mu$ on characters, then, for some $c$, one of $\sigma_k$ and $\sigma$ is a witness for the pair $(u,v)$.
\end{lemma}

\begin{proof}
On $O$, $\sigma$ maps $O$ onto $k\mu O=kO$, as $\sigma_k$ does, and $a+\sigma(a)=a+k\mu a=(1+k)\mu^2a$ maps $O$ onto $(1+k)O$, as $a+\sigma_k(a)$ does; so $\sigma$ is a complete mapping. The coefficient formula is the
definition. Let $d$ be the order of $\mu$ acting on $G$. The function $c\mapsto\sum_{i<d}(\chi_1-\chi_2)(\mu^ic)$ is the combination $\sum_i(\mu^i\chi_1-\mu^i\chi_2)$ of characters, with $(\mu\chi)(a)=\chi(\mu a)$; it is not identically zero, because characters are linearly independent and the $\mu$-orbits of $\chi_1$ and $\chi_2$ are
disjoint. If the orbit of $c$ has $m<d$ elements, each of its elements is counted $d/m$ times in that sum, so for a $c$ with a nonzero value, $\sum_{a\in O}(\chi_1-\chi_2)(a)\neq0$. So $\sigma$ or $\sigma_k$ is a witness.
\end{proof}

\begin{example}[a switch in $\mathbb Z_{25}$]\label{ex:z25}
For the pair $(u,v)=(5,1)$ in $\mathbb Z_{25}$\,, a linear map needs $5+k\equiv0$, that is $k=20$, not a unit. Take $\mu=7$, which has order $4$ ($7^2=49\equiv-1$), so $\mu$ and $\mu+1=8$ are units, and
$k=-\mu^{-1}-1=-18-1\equiv6$ (as $7\cdot18\equiv1$); then $1+k\mu=43\equiv18=(1+k)\mu^2$. Here $\hat\sigma_6(5,1)=0$ because $5+6\neq0$. Now $\chi_1=u+k\mu v=47\equiv22$ and $\chi_2=u+kv=11$, whose orbit under
multiplication by $7$ is $\{11,2,14,23\}$ and does not contain $22$. Take $c=1$, with orbit $O=\{1,7,24,18\}$: the switched map is $\sigma(a)=6a$ except $\sigma(1)=17$, $\sigma(7)=19$, $\sigma(24)=8$, $\sigma(18)=6$ (where $\sigma_6$ has $6,17,19,8$);
it is a complete mapping, with $|\hat\sigma(5,1)|\approx2.64\neq0$. (For $c=5$ the coefficient is $0$: some choice of $c$ is needed.)
\end{example}

\begin{lemma}[twisted products]\label{lem:twist}
Let $G\coloneqq H\oplus K$, and suppose $\sigma_H$ is a complete mapping of $H$ with $\hat\sigma_H(u_H\,,v_H)=0$ and $\sigma_K$ is a witness for $(u_K\,,v_K)$ on $K$, with $v_K\neq0$. Fix $h_0\in H$ and $t\in K$ with $v_K(t)\neq1$, and put
$\sigma(h,x)\coloneqq(\sigma_H(h),\sigma_K(x)+t)$ if $h=h_0$\,, and $\sigma(h,x)\coloneqq(\sigma_H(h),\sigma_K(x))$ otherwise. Then $\sigma$ is a complete mapping and
\[
\hat\sigma(u,v)=\bigl(v_K(t)-1\bigr)\,\hat\sigma_K(u_K\,,v_K)\,u_H(h_0)\,v_H(\sigma_H(h_0))\neq0 .
\]
\end{lemma}

\begin{proof}
For each $h$ the maps $x\mapsto\sigma_K(x)+[h=h_0]t$ and $x\mapsto x+\sigma_K(x)+[h=h_0]t$ are bijections of $K$, and $h\mapsto\sigma_H(h)$, $h\mapsto h+\sigma_H(h)$ are bijections of $H$; so $\sigma$ is a complete mapping.
The coefficient is $\hat\sigma_H\hat\sigma_K+u_H(h_0)v_H(\sigma_H(h_0))\hat\sigma_K\,(v_K(t)-1)$, and $\hat\sigma_H=0$.
\end{proof}

\begin{example}[a twisted product in $\mathbb Z_3\times\mathbb Z_5$]\label{ex:twist}
Let $G\coloneqq\mathbb Z_3\times\mathbb Z_5\cong\mathbb Z_{15}$ and $(u,v)\coloneqq((1,1),(1,2))$, which is allowed. On $\mathbb Z_3$ the restriction is $(1,1)$, of type symbols, so every complete mapping has coefficient $0$; take $\sigma_H\coloneqq\mathrm{id}$.
On $\mathbb Z_5$ it is $(1,2)$, with the linear witness $\sigma_K(x)=2x$ (because $1+2\cdot2=5\equiv0$). With $h_0=0$ and $t=1$, the twisted product is $\sigma(h,x)=(h,2x)$ for $h\neq0$ and $\sigma(0,x)=(0,2x+1)$, and
$\hat\sigma(u,v)=5(\omega_5^{2}-1)\neq0$, where $\omega_5\coloneqq e^{2\pi i/5}$.
\end{example}

\subsection{Odd $p$-groups and couplings}\label{sec:odd}

Let $p$ be an odd prime. We first find witnesses for all allowed pairs on abelian $p$-groups, by treating cyclic ones and then reducing to them. For $p\ge5$ linear maps and switches with a suitable multiplier suffice. For $p=3$ multipliers are useless: if $\mu$ and $\mu+1$ are units modulo $3$, then $\mu\equiv1$, so a switch moves nothing modulo $3$, and we need a nonlinear family. At the end of the subsection we construct the couplings that combine groups of coprime orders.

\emph{Teichm\"uller multipliers.} Let $p\ge5$ be prime and $f\ge1$. The group of units of $\mathbb Z_{p^f}$ is cyclic of order $p^{f-1}(p-1)$, so it contains an element $\mu_p$ of order $d_p\coloneqq p-1$, and because the kernel of reduction modulo $p$ has order $p^{f-1}$, prime to $d_p$\,, the reduction of $\mu_p$ also has order $d_p\ge4$. Then $\mu_p+1$ is a unit (because $\mu_p\not\equiv-1$ modulo $p$), and $1-\mu_p^i$ is $0$ if $d_p\mid i$ and a unit otherwise.

\begin{lemma}[components with $p\ge5$]\label{lem:bigp}
Every allowed pair on $\mathbb Z_{p^f}$\,, $p\ge5$, has a witness.
\end{lemma}

\begin{proof}
Let $\mu\coloneqq\mu_p$ and $k\coloneqq-\mu^{-1}-1$. By Lemma~\ref{lem:switch}, either $\sigma_k$ or a switch is a witness, unless $\chi_1=\mu^j\chi_2$ for some $j$, that is, $u+k\mu v=\mu^j(u+kv)$; because $k\mu=-1-\mu$ and $k=-(1+\mu)\mu^{-1}$, this reads
$(1-\mu^j)u=(1+\mu)(1-\mu^{j-1})v$. If $d_p\mid j$, the left side is $0$ and $(1+\mu)(1-\mu^{j-1})$ is a unit, so $v=0$; if $d_p\mid j-1$, then $u=0$; both are excluded. Otherwise both coefficients are units, so $v=cu$ for a unit $c$.
If $c\not\equiv1$ modulo $p$, the linear map with $k'\coloneqq-c^{-1}$ is a witness: $u+k'v=0$, and $k'$ and $k'+1=(c-1)/c$ are units. If $c\equiv1$ modulo $p$, the equation $\bigl[(1-\mu^j)-c(1+\mu)(1-\mu^{j-1})\bigr]u=0$ with $u\neq0$
forces the bracket to be a non-unit; modulo $p$ the bracket is $\mu(\mu^{j-2}-1)$, so $j\equiv2$ modulo $d_p$\,, and then the original equation reads $(1-\mu)(1+\mu)(u-v)=0$, that is $u=v$, which is excluded.
\end{proof}

\begin{lemma}[cyclic $3$-groups]\label{lem:three}
Every allowed pair on $\mathbb Z_{3^e}$ has a witness.
\end{lemma}

Before the proof, an example of the family used. For $\mathbb Z_9$ take $s\coloneqq5$ (so $s^2=25\equiv-2$ and $s\equiv2$ modulo $3$) and $k\coloneqq s-1=4$, which is $\equiv1$ modulo $3$. Let $F\coloneqq\{0,1,5,6,7,8\}$, let $\pi$ be the permutation of $F$ with cycles $(0\,1\,6\,5)$ and $(7\,8)$,
and let $\sigma(a)\coloneqq4\pi(a)$ on $F$ and $\sigma(a)\coloneqq4a$ elsewhere. Then
\begin{gather*}
\bigl(\sigma(a)\bigr)_{a=0,\dots,8}=(4,6,8,3,7,0,2,5,1),\\
\bigl(a+\sigma(a)\bigr)_{a=0,\dots,8}=(4,7,1,6,2,5,8,3,0),
\end{gather*}
both permutations of $\mathbb Z_9$\,, so $\sigma$ is a complete mapping. For the pair $(u,v)=(6,1)$, that is $(-3,1)$, the pair for which linear maps fail, $|\hat\sigma(6,1)|\approx7.60\neq0$. (The pair $(3,1)$ then follows by the symmetries, as in the proof.)

\begin{proof}[Proof of Lemma~\ref{lem:three}]
\emph{Step 1: reduction to $(\pm3^{e-1},1)$.} By induction on $e$. On $\mathbb Z_3$ the allowed pairs are $(1,2)$ and $(2,1)$, for which $u+v=0$ and the identity is a witness.
Let $e\ge2$. If $3\mid u$ and $3\mid v$, both characters are trivial on the subgroup of order $3$ and the induced pair on the quotient $\mathbb Z_{3^{e-1}}$ has a witness by induction, so Lemma~\ref{lem:lift}(i) applies. If
$3\nmid u,v,u-v$, the restriction to the subgroup of order $3$ is allowed, and Lemma~\ref{lem:lift}(ii) applies. Otherwise exactly one of $u$, $v$, $u-v$ is divisible by $3$, and by the symmetries (S2), (S3) we may assume that it is $u$:
$3^m\,\|\,u$ with $1\le m\le e-1$ and $3\nmid v$. If $m\le e-2$, the restriction of $(u,v)$ to the subgroup $N\cong\mathbb Z_{3^{m+1}}$ is $(u\bmod3^{m+1},v\bmod3^{m+1})=(3^mu',v)$ with $u'$ and $v$ units, hence allowed, with a witness by induction, so
Lemma~\ref{lem:lift}(ii) applies. So let $m\coloneqq e-1$. Multiplying by $v^{-1}$ (S1), we may take $(u,v)=(\varepsilon3^{e-1},1)$ with $\varepsilon=\pm1$; and the two signs are exchanged by the symmetries:
$(3^{e-1},1)\to(1,3^{e-1})\to(3^{e-1}-1,3^{e-1})\to(3^{e-1},3^{e-1}-1)$ by (S2), (S3), (S2), and multiplying by $(3^{e-1}-1)^{-1}=-1-3^{e-1}$ gives $(-3^{e-1},1)$. So it suffices to find a witness for $(-3^{e-1},1)$, say, and, for $e\ge3$, either sign works for the argument below.

\emph{Step 2: the family.} Let $n\coloneqq3^e$. Let $s$ satisfy $s^2\equiv-2\pmod n$ and $s\equiv2\pmod3$; it exists because the solution $s=2$ modulo $3$ (as $-2\equiv1\equiv2^2$) lifts, one power of $3$ at a time, to a solution modulo $n$, the derivative $2s$ being a unit (Hensel's lemma). Put $k\coloneqq s-1$ (so $k\equiv1$ modulo $3$ and $1+k=s$ are units) and $y\coloneqq-s^{-1}$,
let $F\coloneqq\{0,1,s,s+1,y,y+1\}$, let $\pi$ be the permutation of $F$ with cycles $(0,1,s+1,s)$ and $(y,y+1)$, and let $\sigma(a)\coloneqq k\pi(a)$ on $F$ and $\sigma(a)\coloneqq ka$ elsewhere. Then $\sigma$ permutes $kR$ as $\sigma_k$ does, and, using $s^2=-2$ and $sy=-1$, the symbols $a+k\pi(a)$ for
$a=0,1,s+1,s,y,y+1$ are $s-1,\,-2,\,-1,\,s,\,s-2,\,0$, which are exactly the symbols $sa=(1+k)a$ of $a\in F$. So $\sigma$ is a complete mapping.

\emph{Step 3: its coefficient.} Let $\zeta\coloneqq e^{2\pi i/n}$ and $\omega\coloneqq\zeta^{3^{e-1}}$. For $(u,v)=(\varepsilon3^{e-1},1)$ we have $\hat\sigma_k(u,v)=0$ (as $k+\varepsilon3^{e-1}$ is a unit), so
\[
\hat\sigma(u,v)=\sum_{a\in F}\omega^{\varepsilon a}\bigl(\zeta^{k\pi(a)}-\zeta^{ka}\bigr).
\]
A vanishing integer combination of $n$th roots of unity has coefficients constant on the cosets of the subgroup of order $3$, because the relations among them are generated by $\zeta^j+\zeta^{j+3^{e-1}}+\zeta^{j+2\cdot3^{e-1}}=0$.
Modulo $3^{e-1}$ the twelve exponents are $k\pi(a)$ and $ka$. For $e\ge3$ the elements of $F$ are distinct modulo $3^{e-1}$: their residues modulo $3$ are $0,1,2,0,1,2$, and a coincidence $s+1\equiv0$, $y\equiv1$ or $y+1\equiv s$ would force $s\equiv-1$, hence
$s^2\equiv1\equiv-2$ and $3\equiv0$ modulo $3^{e-1}$. So the coset of $kr$, for $r\in F$, carries exactly the two terms $+\zeta^{\varepsilon3^{e-1}\pi^{-1}(r)+kr}$ and $-\zeta^{\varepsilon3^{e-1}r+kr}$, which occupy different exponents of the coset unless $\pi^{-1}(r)\equiv r$ modulo $3$, in which case they cancel.
Hence the sum vanishes only if $\pi^{-1}(r)\equiv r$ modulo $3$ for all $r\in F$. But $\pi^{-1}(1)=0\not\equiv1$. So $\hat\sigma(u,v)\neq0$. For $e=2$, where $s=5$ and $k=4$, the example above exhibits a nonzero coefficient for $(-3,1)$.
\end{proof}

\begin{lemma}[odd $p$-groups]\label{lem:pgroups}
Every allowed pair on an abelian $p$-group $G$, $p$ odd, has a witness.
\end{lemma}

\begin{proof}
Let $K_u\coloneqq\ker(u)$, $K_v\coloneqq\ker(v)$, and $K_{u-v}\coloneqq\ker(u-v)$. Each is a proper subgroup, of index at least $p\ge3$. Any two of $u$, $v$, $u-v$ determine the third, so the pairwise intersections of these
all equal $K_0\coloneqq K_u\cap K_v$\,. By inclusion and exclusion, $|K_u\cup K_v\cup K_{u-v}|=|K_u|+|K_v|+|K_{u-v}|-2|K_0|\le|G|-2|K_0|<|G|$. So there is $g\in G$ outside all three, that is, with $u(g),v(g),(u-v)(g)\neq1$; then $(u,v)$ restricts to an
allowed pair on the cyclic $p$-group $N\coloneqq\langle g\rangle$, which has a witness by Lemma~\ref{lem:bigp} or Lemma~\ref{lem:three}, and Lemma~\ref{lem:lift}(ii) lifts it to $G$, the identity being a complete mapping of $G/N$. (For instance, in $\mathbb Z_3^2$ with $u=(1,0)$ and $v=(0,1)$, the three kernels have $3+3+3-2=7$ elements in total, and $g=(1,2)$ lies outside them; $\langle g\rangle\cong\mathbb Z_3$ carries the allowed pair $(1,2)$.)
\end{proof}

\begin{lemma}[layered coupling]\label{lem:two}
Let $G\coloneqq H\oplus C'$, where $H$ is an abelian group with a complete mapping $\rho$ and $C'$ is cyclic of order $r^f$ with $r\ge5$ prime. Let $u$ be a nontrivial character of $H$ and $v$ a character of $C'$ of order $r$, both extended trivially to $G$. Then $(u,v)$ has a witness.
\end{lemma}

\emph{Idea.} The pair has type rows on $H$ and columns on $C'$, so no product of witnesses is available. We split $C'$ into three \emph{layers} according to quadratic residuosity modulo $r$, apply $\rho$ on two layers and $\rho$ shifted by some $t$ on the third, and multiply the $C'$-coordinate by integers whose quadratic characters keep the map bijective; the Gauss periods of $\mathbb Z_r$ then make the coefficient a nonzero multiple of $u(i_0)-u(\pi(i_0))$.

\begin{proof}
Split $C'=\mathbb Z_{r^f}$ into the layers $Q_0\coloneqq rC'$ (the non-units), $Q_1$ (the units whose residue modulo $r$ is a square), and $Q_2$ (the units with a nonsquare residue), and write $\ell(h)\in\{0,1,2\}$ for the layer of $h$. Write $\nu(k)\in\{0,1\}$ for the quadratic character of $k$ modulo $r$ ($0$ for squares); multiplication by $k$ maps $Q_1$ onto $Q_{1+\nu(k)}$ and $Q_2$ onto $Q_{2-\nu(k)}$\,, and $k$ and $k+1$ are prime to $r$ for $k\in[1,r-2]$. Fix $t\in H$, $t\neq0$, and let $\rho_0\coloneqq\rho_1\coloneqq\rho$ and $\rho_2\coloneqq\rho-t$, which are complete mappings of $H$. For integers $k_{i,1}\,,k_{i,2}\in[1,r-2]$, $i\in H$, and $k_{i,0}\coloneqq1$, put
\[
\sigma(i,h)\coloneqq\bigl(\rho_{\ell(h)}(i),\ k_{i,\ell(h)}\,h\bigr).
\]

\emph{When $\sigma$ is a complete mapping.} Let $\pi(i)\coloneqq\rho^{-1}(\rho(i)+t)$, $\theta\coloneqq\mathrm{id}+\rho$, and $\pi'(i)\coloneqq\theta^{-1}(\theta(i)+t)$. The preimages of a point $(\rho(i),h')$ with $h'\in Q_1$ come from layer-$1$ sources $(i,h)$ if $\nu(k_{i,1})=0$, and from layer-$2$ sources $(\pi(i),h)$ if $\nu(k_{\pi(i),2})=1$; there must be exactly one, which holds exactly when $\nu(k_{i,1})=\nu(k_{\pi(i),2})$. The same condition arises for $Q_2$\,, and layer $0$ is mapped by $\rho\times\mathrm{id}$. So $\sigma$ is a bijection if and only if $\nu(k_{i,1})=\nu(k_{\pi(i),2})$ for all $i$. For the symbols, $(i,h)\mapsto(\theta(i)-t[\ell(h)=2],(1+k_{i,\ell(h)})h)$, and the same count establishes that $a\mapsto a+\sigma(a)$ is a bijection if and only if $\nu(1+k_{i,1})=\nu(1+k_{\pi'(i),2})$ for all $i$.

\emph{Its coefficient.} Because $v$ has order $r$, it is trivial on $Q_0$\,, and its sums over the layers are Gauss periods: let $\eta_0$ and $\eta_1$ be the sums of $\psi(a)$ over the nonzero squares, respectively nonsquares, modulo $r$, where $\psi$ is the character of $\mathbb Z_r$ through which $v$ factors. They satisfy $\eta_0+\eta_1=-1$ and $\eta_0\neq\eta_1$\,, so $1+2\eta_0\neq0$ and $1+2\eta_1=-(1+2\eta_0)$. Summing over $h$, $\hat\sigma(u,v)=\sum_iu(i)Y_i$ with $Y_i\coloneqq r^{f-1}(1+\eta_{\nu_i}+\eta_{1-\nu'_i})$, where $\nu_i\coloneqq\nu(k_{i,1})$ and $\nu'_i\coloneqq\nu(k_{i,2})$ (indices of $\eta$ taken in $\{0,1\}$). So $Y_i=0$ if $\nu_i=\nu'_i$\,, $Y_i=r^{f-1}(1+2\eta_0)$ if $(\nu_i\,,\nu'_i)=(0,1)$, and $Y_i=-r^{f-1}(1+2\eta_0)$ if $(\nu_i\,,\nu'_i)=(1,0)$.

\emph{Choice of parameters.} Choose every $\nu(1+k)=1$, which satisfies the second condition; choose $\nu_{i_0}=0$ for one $i_0$ and $\nu_i=1$ for $i\neq i_0$\,, and put $\nu'_i\coloneqq \nu_{\pi^{-1}(i)}$\,, which satisfies the first. Integers $k$ with these values exist: with $\chi(k)\coloneqq(-1)^{\nu(k)}$ for $k$ prime to $r$ (the Legendre symbol), the number of $k\in[1,r-2]$ with $(\nu(k),\nu(k+1))=(a,b)$ is
\[
\tfrac14\bigl(r-2-(-1)^a\chi(-1)-(-1)^b-(-1)^{a+b}\bigr),
\]
which is $\frac14(r-\chi(-1))>0$ for $(a,b)=(0,1)$ and $\frac14(r-2+\chi(-1))>0$ for $(1,1)$, whenever $r\ge5$. Because $\pi$ has no fixed point ($\pi(i)=i$ would mean $t=0$), $\nu_i\neq \nu'_i$ exactly for $i=i_0$ and $i=\pi(i_0)$, with $(\nu,\nu')=(0,1)$ and $(1,0)$ there, and the coefficient is $r^{f-1}(1+2\eta_0)\bigl(u(i_0)-u(\pi(i_0))\bigr)$. If $u(\pi(i))=u(i)$ for all $i$ and all $t\neq0$, then $\rho^{-1}(x+t)-\rho^{-1}(x)\in\ker(u)$ for all $x$ and $t$, so the bijection $\rho^{-1}$ would take all its values in one coset of the proper subgroup $\ker(u)$, which is absurd. So for some $t$ and $i_0$ the coefficient is nonzero.
\end{proof}

\begin{example}[coupling in $\mathbb Z_3\times\mathbb Z_5$]\label{ex:qr}
Let $G\coloneqq\mathbb Z_3\times\mathbb Z_5\cong\mathbb Z_{15}$\,, with $u\coloneqq(1,0)$ and $v\coloneqq(0,1)$: rows on $\mathbb Z_3$\,, columns on $\mathbb Z_5$; no product of witnesses applies. Take $\rho$ the identity of $\mathbb Z_3$\,, $t\coloneqq1$, and $i_0\coloneqq0$, so $\pi(i)=i+1$. The squares modulo $5$ are $\{1,4\}$, so the layers are $Q_0=\{0\}$, $Q_1=\{1,4\}$, $Q_2=\{2,3\}$, and the $k\in[1,3]$ with $\nu(1+k)=1$ are $k=1$, with $\nu(k)=0$, and $k=2$, with $\nu(k)=1$. With $\nu=(0,1,1)$ and $\nu'_i=\nu_{i-1}$\,, that is $\nu'=(1,0,1)$, we get $(k_{i,1})=(1,2,2)$ and $(k_{i,2})=(2,1,2)$. The map $\sigma(i,h)=(i-[\ell(h)=2],k_{i,\ell(h)}h)$ is
\[
\begin{array}{c|ccccc}
 & h=0 & h=1 & h=2 & h=3 & h=4\\\hline
i=0 & (0,0) & (0,1) & (2,4) & (2,1) & (0,4)\\
i=1 & (1,0) & (1,2) & (0,2) & (0,3) & (1,3)\\
i=2 & (2,0) & (2,2) & (1,4) & (1,1) & (2,3)
\end{array}
\]
and one checks that $\sigma$ and $a\mapsto a+\sigma(a)$ are bijections of $G$, and $|\hat\sigma(u,v)|=|1+2\eta_0|\,|1-\omega_3|=\sqrt5\cdot\sqrt3=\sqrt{15}$, where $\omega_3\coloneqq e^{2\pi i/3}$.
\end{example}

\subsection{Abelian groups of even order}\label{sec:even}

When $|G|$ is even, its Sylow $2$-subgroup $P$ is not cyclic, by the theorem of Hall and Paige recalled in \S\ref{sec:fourier2}. Several tools used for odd $p$-groups break down. The reduction to cyclic subgroups in Lemma~\ref{lem:pgroups} fails: $\mathbb Z_2^2$ is
the union of the kernels of its three nonzero characters, so there is no cyclic subgroup carrying an allowed restriction. Integer multipliers fail: there is no integer $k$ with $k$ and $k+1$ both odd. And
Lemma~\ref{lem:lift}(i) needs a complete mapping of the subgroup $N$, which for a $2$-group means that $N$ is noncyclic; this is a serious restriction (part (ii) needs one on $N$ as well). In $P=\mathbb Z_2\times\mathbb Z_{2^b}$\,, for instance, no nontrivial proper subgroup $N$ has complete mappings on both $N$ and $P/N$, because one of the two is always cyclic, so no lifting is available at all. The structure that remains is that of a \emph{homocyclic} $2$-group $(\mathbb Z_{2^e})^m$, $m\ge2$, which is the additive group of a Galois ring; we treat these first, then all noncyclic $2$-groups, and finally the coupling with an odd part.

\emph{Galois rings.} For $\mathbb Z_2^m$ the field $\GF(2^m)$ supplies multipliers, units $k$ with $k+1$ also a unit, and hence the linear complete mappings $a\mapsto ka$; Galois rings do the same for $(\mathbb Z_{2^e})^m$. For $e\ge1$ and $m\ge2$, let $R\coloneqq\mathrm{GR}(2^e,m)=\mathbb Z_{2^e}[x]/(f)$, where $f$ is a monic polynomial of degree $m$ over $\mathbb Z_{2^e}$ whose reduction modulo $2$ is irreducible.
For $e=1$ this is the field $\GF(2^m)$; in general the additive group of $R$ is $(\mathbb Z_{2^e})^m$ (with basis $1,x,\dots,x^{m-1}$). The following facts are standard (see Wan~\cite{Wan}).
An element of $R$ is a unit if and only if its reduction modulo $2$, an element of the residue field $R/2R\cong\GF(2^m)$, is nonzero. The unit group contains a cyclic subgroup of order $2^m-1$, whose elements are the \emph{Teichm\"uller elements}; it maps isomorphically onto $\GF(2^m)^*$. With $\operatorname{Tr}(a)\in\mathbb Z_{2^e}$ the trace of the $\mathbb Z_{2^e}$-linear map $b\mapsto ab$ of $R$, the \emph{trace form} $(a,b)\mapsto\operatorname{Tr}(ab)$ is nondegenerate, so the characters of $(R,+)$ are $a\mapsto e^{2\pi i\operatorname{Tr}(wa)/2^e}$, $w\in R$; we identify
characters with elements of $R$, and then $(kw)(a)=w(ka)$ is the ring product $kw$. Everything in \S\ref{sec:tools} about multipliers therefore applies to $G=(R,+)$ with units of $R$ as multipliers.

\begin{example}[$P=\mathbb Z_2^2=\GF(4)$]\label{ex:gf4}
Let $\GF(4)=\{0,1,\mu,\mu^2\}$ with $\mu^2=\mu+1$, so $1+\mu=\mu^2$ and $1+\mu^2=\mu$, and let $\operatorname{Tr}(a)\coloneqq a+a^2$, which is $0$ for $a\in\{0,1\}$ and $1$ for $a\in\{\mu,\mu^2\}$. The characters are $\chi_w(a)\coloneqq(-1)^{\operatorname{Tr}(wa)}$.
For an allowed pair $(u,v)$ we have $u,v\in\GF(4)^*$ and $u\neq v$, so $k\coloneqq u/v\in\{\mu,\mu^2\}$; then $k$ and $k+1=k^2$ are units, $\sigma_k(a)=ka$ is a complete mapping, and $u+kv=u+u=0$ in characteristic $2$, so
$\hat\sigma_k(u,v)=4\neq0$. All $6$ allowed pairs have linear witnesses, which gives $\dim=(4-1)(4-2)=6$, as for the $8$ transversals of the Klein four-group.
\end{example}

\begin{lemma}[homocyclic $2$-groups]\label{lem:homocyclic}
Let $P\cong(\mathbb Z_{2^e})^m$ with $m\ge2$. Every allowed pair has a witness.
\end{lemma}

\begin{proof}
Let $P\coloneqq(R,+)$ with $R=\mathrm{GR}(2^e,m)$, and let $\mu\in R$ be a Teichm\"uller element of order $d\coloneqq2^m-1\ge3$. Its residue has order $d$ as well and differs from $0$ and $1$, so $\mu$ and $\mu+1$ are units, and $1-\mu^i$ is a unit exactly when $d\nmid i$. These are the only properties of $\mu_p$ used in the proof of Lemma~\ref{lem:bigp}, which therefore applies verbatim, with $\mathbb Z_{p^f}$ replaced by $R$, reduction modulo $p$ by reduction modulo $2$, and $d_p$ by $d$; Lemma~\ref{lem:switch} is valid for units of $R$.
\end{proof}

\emph{Beyond homocyclic $2$-groups.} Next we treat $2$-groups that are not homocyclic, with three further tools: a general design built over a cyclic subgroup, an explicit family of complete mappings of $\mathbb Z_{2^a}\times\mathbb Z_{2^b}$ with a product formula for its coefficients, and local modifications by $3$-cycles. Throughout, $A\coloneqq2^a$, $B\coloneqq2^b$, $\zeta_N\coloneqq e^{2\pi i/N}$, and characters of $\mathbb Z_A\times\mathbb Z_B$ are written $(u_1\,,u_2)$ with $u_1\in\mathbb Z_A$\,, $u_2\in\mathbb Z_B$\,, so that $(u_1\,,u_2)(g,h)=\zeta_A^{u_1g}\zeta_B^{u_2h}\,$.

\begin{lemma}[extension design]\label{lem:ext}
Let $K\coloneqq\langle\kappa\rangle\le G$ be cyclic of even order $m$, let $Q\coloneqq G/K$, let $s\colon Q\to G$ be a section with $s(0)=0$, and let $c\colon Q\times Q\to\mathbb Z_m$\,, the \emph{cocycle} of $s$, be defined by $s(x)+s(y)=s(x+y)+c(x,y)\kappa$; it records the carries that occur when representatives are added. Let $\varphi_0\,,\varphi_1$ be permutations of $Q$, let $t\colon Q\to\mathbb Z_m$\,, and define, for $x\in Q$ and $0\le i<m$ with $h(i)\coloneqq[i\ge m/2]$,
\[
\sigma\bigl(s(x)+i\kappa\bigr)\coloneqq s\bigl(\varphi_{h(i)}(x)\bigr)+(i+t_x)\kappa .
\]
Suppose that (A) $t_{\varphi_0^{-1}(y)}=t_{\varphi_1^{-1}(y)}$ for all $y\in Q$, and (B) the map $(x,h)\mapsto\bigl(x+\varphi_h(x),\,t_x+c(x,\varphi_h(x))\bmod2\bigr)$ is a bijection from $Q\times\{0,1\}$ onto $Q\times\mathbb Z_2$\,. Then $\sigma$ is a complete mapping of $G$, and for characters $u,v$ of $G$,
\[
\hat\sigma(u,v)=S_0\sum_{x\in Q}u(s(x))\,v(t_x\kappa)\bigl[v(s(\varphi_0(x)))+\varepsilon\,v(s(\varphi_1(x)))\bigr],
\]
where $S_0\coloneqq\sum_{i<m/2}(uv)(i\kappa)$ and $\varepsilon\coloneqq(uv)(\tfrac m2\kappa)$.
In particular $S_0\neq0$ if and only if the restriction of $uv$ to $K$ is trivial ($\varepsilon=1$) or of order $m$ ($\varepsilon=-1$).
\end{lemma}

\begin{proof}
For a target $s(y)+i'\kappa$, the preimages under $\sigma$ are the $s(x)+i\kappa$ with $\varphi_{h(i)}(x)=y$ and $i=i'-t_x$; for each $h\in\{0,1\}$ there is one candidate $x_h\coloneqq\varphi_h^{-1}(y)$, and by (A) the two values $t_{x_h}$ agree, so exactly one of $i'-t_{x_0}$\,, $i'-t_{x_1}$ lies in the half $h$ it needs. So $\sigma$ is a bijection. The symbol of $s(x)+i\kappa$ is $s(x+\varphi_h(x))+(c(x,\varphi_h(x))+t_x+2i)\kappa$ with $h=h(i)$; as $i$ runs over a half, $2i$ runs over $2K$ once, so the symbols from $(x,h)$ form the coset $s(x+\varphi_h(x))+(c(x,\varphi_h(x))+t_x)\kappa+2K$, and (B) says that these cosets partition $G$. For the coefficient, $u(s(x)+i\kappa)=u(s(x))u(i\kappa)$, so
$\hat\sigma(u,v)=\sum_xu(s(x))v(t_x\kappa)\sum_i(uv)(i\kappa)\,v(s(\varphi_{h(i)}(x)))$, and the inner sum is $S_0v(s(\varphi_0x))+S_1v(s(\varphi_1x))$ with $S_1=\varepsilon S_0$\,. Finally $S_0$ is a geometric sum, zero exactly when the restriction of $uv$ to $K$ has order strictly between $1$ and $m$.
\end{proof}

When $u$ and $v$ are trivial on $K$ and $\varphi_1=\varphi_0+d$ for some $d\in Q$, the formula reads $\hat\sigma(u,v)=\frac m2\bigl(1+v(s(d))\bigr)\hat\varphi_0(\bar u,\bar v)$, where $\bar u,\bar v$ are the induced characters of $Q$: a witness on $Q$ lifts to a witness on $G$ whenever (A) and (B) can be met, even though $K$ has no complete mapping. We use the lemma in two concrete forms.

\begin{lemma}[cyclic-quotient pairs]\label{lem:cq}
Let $K\le G$ be cyclic of even order $m$ such that $Q\coloneqq G/K$ is cyclic of order $2^k$, with a section $s$ whose cocycle takes values in $2K$. Let $u$ be a character of $G$ of order $2^j$, $2\le j\le k$, that is trivial on $K$. Then $(u,2^{j-1}u)$ has a witness, and so has every pair obtained from it by the symmetries; for $j=2$ these are all allowed pairs whose two characters lie in $\{u,-u,2u\}$. In particular this applies to $G\coloneqq\mathbb Z_A\times\mathbb Z_B$ with $1\le a<b$, $u\coloneqq(A/2,B/4)$, and $K\coloneqq\langle(1,2)\rangle$.
\end{lemma}

\begin{proof}
Write the elements of $Q$ as $y\in\{0,\dots,2^k-1\}$ with binary digits $y_0\,,\dots,y_{k-1}$\,, and define $\varphi_0(y)$ to be the number whose digit $0$ is $y_{j-1}$\,, whose digits $1,\dots,j-2$ are $y_0\,,\dots,y_{j-3}$\,, whose digits $j-1,\dots,k-2$ are $y_j\,,\dots,y_{k-1}$\,, and whose digit $k-1$ is $y_{j-2}$; let $\varphi_1\coloneqq\varphi_0+2^{k-1}$ and $t_y\coloneqq y_{j-1}$\,. We apply Lemma~\ref{lem:ext}. Condition (A) holds because $\varphi_0$ and $\varphi_1$ differ in digit $k-1$ only, so $\varphi_1^{-1}(z)$ and $\varphi_0^{-1}(z)$ differ in digit $j-2$ only, which leaves $t$ unchanged. For (B), write $y=\ell+2^{j-2}y_{j-2}+2^{j-1}y_{j-1}+2^jH$ with $0\le\ell<2^{j-2}$ and $0\le H<2^{k-j}$, so that $\varphi_h(y)=y_{j-1}+2\ell+2^{j-1}H+2^{k-1}(y_{j-2}+h\bmod2)$ and
\[
y+\varphi_h(y)=3\ell+(2^{j-1}+1)y_{j-1}+2^{j-2}y_{j-2}+3\cdot2^{j-1}H+2^{k-1}(y_{j-2}+h\bmod2).
\]
Because the cocycle takes values in $2K$, (B) says that for each $p\in\{0,1\}$ the map sending $(\ell,y_{j-2}\,,h,H)$ to $y+\varphi_h(y)$, with $y_{j-1}=p$, is a bijection onto $Q$. Modulo $2^{j-2}$ the symbol is $3\ell+p$, and $\ell\mapsto3\ell+p$ is a bijection onto $\mathbb Z_{2^{j-2}}$ (for $j=2$ there is no $\ell$ and nothing to check); so it suffices that for fixed $\ell$ the map $(y_{j-2}\,,h,H)\mapsto(y+\varphi_h(y)-3\ell-p)/2^{j-2}$ is a bijection onto $\mathbb Z_{2^{m'}}$\,, $m'\coloneqq k-j+2$. Up to the constant $2p$, this map is $y_{j-2}+2^{m'-1}(y_{j-2}+h\bmod2)+6H$, whose values are $2\mathcal S\cup(2\mathcal S+2^{m'-1})$ for $y_{j-2}=0$ and $1+2\mathcal S\cup(1+2\mathcal S+2^{m'-1})$ for $y_{j-2}=1$, where $\mathcal S\coloneqq\{3H\bmod2^{m'-1} ~:~ H<2^{m'-2}\}$ is a transversal of the cosets of $2^{m'-2}$ in $\mathbb Z_{2^{m'-1}}$ because $3$ is a unit; so the values are all of $\mathbb Z_{2^{m'}}$\,.

The coefficient: $u$ and $v\coloneqq2^{j-1}u$ are trivial on $K$, so $S_0=m/2$, $\varepsilon=1$, and $v(t_y\kappa)=1$; the induced characters are $\bar u(y)=\zeta_{2^j}^{oy}$ with $o$ odd and $\bar v(y)=(-1)^y$; the parity of $\varphi_h(y)$ is $y_{j-1}$ for both $h$; and $\bar u(y)$ depends only on $y$ modulo $2^j$. Hence
\begin{align*}
\hat\sigma(u,v)&=\frac m2\cdot2\cdot2^{k-j}\sum_{y<2^j}\zeta_{2^j}^{oy}(-1)^{y_{j-1}}\\
&=m\,2^{k-j}\bigl(1-\zeta_{2^j}^{o2^{j-1}}\bigr)\sum_{y<2^{j-1}}\zeta_{2^j}^{oy}=m\,2^{k-j+1}\cdot\frac{-2}{\zeta_{2^j}^o-1}\neq0.
\end{align*}
For $j=2$ the symmetries give the other pairs: (S2) gives $(2u,u)$, (S3) maps $(2u,u)$ to $(-u,u)$, (S2) then gives $(u,-u)$, and the automorphism $a\mapsto-a$ maps $(u,2u)$ and $(2u,u)$ to $(-u,2u)$ and $(2u,-u)$. Finally, for $G=\mathbb Z_A\times\mathbb Z_B$ and $K=\langle(1,2)\rangle$, of order $B/2$, the quotient is cyclic of order $2A$, generated by the image of $(0,1)$, with $s(x)=(0,x)$ and cocycle $A[x+y\ge2A]\kappa$ because $(0,2A)=A(1,2)$; $A$ is even, and $(A/2,B/4)$ is trivial on $K$ and of order $4$.
\end{proof}

The lemma handles, for a character $u$ of $2$-power order $2^j\ge4$, the pair $(u,2^{j-1}u)$, which factors through the cyclic quotient $G/\ker(u)$ and is invisible to the designs built from $2$-adic fibers, such as the family below.

\begin{lemma}[the $2^a$-adic family]\label{lem:family}
Let $1\le a<b$, $L\coloneqq B/A$, and $c\coloneqq(A+1)^{-1}$ modulo $B$. For a permutation $\pi$ of $\mathbb Z_A$ and a function $\beta\colon\mathbb Z_L\to\mathbb Z_A$\,, put $\alpha(h)\coloneqq\pi(\lfloor h/L\rfloor)+\beta(h\bmod L)$ for $0\le h<B$ and
\[
\sigma(g,h)\coloneqq\bigl(\alpha(h+cg),\ Ah+g\bigr)\qquad(0\le g<A,\ h\in\mathbb Z_B).
\]
Then $\sigma$ is a complete mapping of $\mathbb Z_A\times\mathbb Z_B$\,, and for characters $u=(u_1\,,u_2)$ and $v=(v_1\,,v_2)$, with $w\coloneqq u_2+Av_2$\,,
\begin{gather*}
\hat\sigma(u,v)=\Bigl(\sum_{g<A}\xi^g\Bigr)\,\hat\pi\,\hat\beta,\qquad
\xi\coloneqq\zeta_A^{u_1}\zeta_B^{v_2-wc}\,,\\
\hat\pi\coloneqq\sum_{j\in\mathbb Z_A}\zeta_A^{wj+v_1\pi(j)}\,,\qquad
\hat\beta\coloneqq\sum_{r<L}\zeta_B^{wr}\zeta_A^{v_1\beta(r)}\,.
\end{gather*}
\end{lemma}

\begin{proof}
The map $(g,h)\mapsto Ah+g$ is $A$-to-one: the preimages of $s$ are the $(g,h)$ with $g=s\bmod A$ and $h\equiv(s-g)/A$ modulo $L$, that is, $h=h_0+jL$ for $0\le j<A$. On such a fiber, $\alpha(h+cg)=\pi(\lfloor h_0'/L\rfloor+j)+\beta(h_0'\bmod L)$ with $h_0'\coloneqq h_0+cg$ takes every value once; so $\sigma$ is a bijection. The symbol of $(g,h)$ is $(g+\alpha(h+cg),(A+1)h+g)$; the second coordinate determines $h=c(s'-g)$ for each $g$, and then $g+\alpha(h+cg)=g+\alpha(cs')$, a translate of $g$, so the symbols are distinct. For the coefficient, substitute $h'\coloneqq h+cg$: the sum factors as $\sum_g\zeta_A^{u_1g}\zeta_B^{(v_2-wc)g}$ times $\sum_{h'}\zeta_B^{wh'}\zeta_A^{v_1\alpha(h')}\,$, and writing $h'=r+jL$ the second factor is $\hat\beta\sum_j\zeta_A^{wj+v_1\pi(j)}\,$.
\end{proof}

\begin{corollary}\label{cor:family}
Let $u=(u_1\,,u_2)$ and $v=(v_1\,,v_2)$ be characters of $\mathbb Z_A\times\mathbb Z_B$\,, $1\le a<b$, and put $w\coloneqq u_2+Av_2$ and $\xi\coloneqq\zeta_A^{u_1}\zeta_B^{v_2-wc}\,$. Some member of the family witnesses $(u,v)$ if and only if $\sum_{g<A}\xi^g\neq0$, $w\equiv0$ modulo $A$ exactly when $v_1=0$, and $w\equiv0$ modulo $B$ when $v_1=0$; otherwise every member has coefficient $0$ there.
\end{corollary}

\begin{proof}
The first factor in Lemma~\ref{lem:family} does not depend on $\pi$ or $\beta$. The factor $\hat\pi$ is the coefficient of the permutation $\pi$ of $\mathbb Z_A$ at the characters $w\bmod A$ and $v_1$; it vanishes for every $\pi$ if exactly one of them is $0$, it equals $A$ if both are $0$, and if both are nonzero it is nonzero for some $\pi$, because the centered incidence vectors of the permutations of $\mathbb Z_A$ span the space of dimension $(A-1)^2$ orthogonal to the line sums (the Birkhoff polytope; see~\cite{BR}). The factor $\hat\beta$ can be made nonzero by the choice of $\beta$ when $v_1\neq0$ (changing $\beta(0)$ changes it by a nonzero amount), and for $v_1=0$ it is $\sum_{r<L}\zeta_B^{wr}\,$, which for $A\mid w$ is nonzero exactly when $B\mid w$.
\end{proof}

\begin{lemma}[$3$-cycles]\label{lem:cycle}
Let $\sigma$ be a complete mapping of $G$ and let $a_1\,,a_2\,,a_3\in G$ be distinct with $\{a_1+\sigma(a_2),\,a_2+\sigma(a_3),\,a_3+\sigma(a_1)\}=\{\theta(a_1),\theta(a_2),\theta(a_3)\}$, where $\theta\coloneqq\mathrm{id}+\sigma$. Then $\sigma'$, defined by $\sigma'(a_1)\coloneqq\sigma(a_2)$, $\sigma'(a_2)\coloneqq\sigma(a_3)$, $\sigma'(a_3)\coloneqq\sigma(a_1)$, and $\sigma'=\sigma$ elsewhere, is a complete mapping, and
\begin{multline*}
\hat\sigma'(u,v)-\hat\sigma(u,v)=u(a_1)\bigl(v(\sigma a_2)-v(\sigma a_1)\bigr)\\
+u(a_2)\bigl(v(\sigma a_3)-v(\sigma a_2)\bigr)+u(a_3)\bigl(v(\sigma a_1)-v(\sigma a_3)\bigr).
\end{multline*}
\end{lemma}

\begin{proof}
The images are permuted and the symbols are permuted among themselves; the formula is the definition.
\end{proof}

\begin{lemma}[$3$-cycles in the family]\label{lem:shape}
Let $\sigma$ be a member of the family of Lemma~\ref{lem:family}, with $D\coloneqq(A^2+A+1)^{-1}$ modulo $B$, and let $g_1\,,g_2\,,g_3\in\{0,\dots,A-1\}$ and $h_1\in\mathbb Z_B$\,. Put
\[
h_3-h_2\coloneqq D\bigl((A+1)g_2-Ag_3-g_1\bigr),\qquad h_1-h_3\coloneqq A(h_3-h_2)+g_3-g_2\,,
\]
and $a_i\coloneqq(g_i\,,h_i)$, assumed distinct. Then $a_1\,,a_2\,,a_3$ satisfy the hypothesis of Lemma~\ref{lem:cycle} if $g_1+\alpha(h_2')=g_3+\alpha(h_3')$ and $g_2+\alpha(h_3')=g_1+\alpha(h_1')$ modulo $A$, where $h_i'\coloneqq h_i+cg_i$\,.
\end{lemma}

\begin{proof}
The symbols are $\theta(g,h)=(g+\alpha(h'),(A+1)h+g)$ and $a_1+\sigma(a_2)=(g_1+\alpha(h_2'),h_1+Ah_2+g_2)$, with the second coordinates computed with $g_i$ as integers. The two displayed relations are the unique solution of $h_1+Ah_2+g_2=(A+1)h_3+g_3$ and $h_2+Ah_3+g_3=(A+1)h_1+g_1$\,, which with their consequence $h_3+Ah_1+g_1=(A+1)h_2+g_2$ say that the second coordinates of $a_1+\sigma(a_2)$, $a_2+\sigma(a_3)$, $a_3+\sigma(a_1)$ are those of $\theta(a_3)$, $\theta(a_1)$, $\theta(a_2)$; the two conditions on $\alpha$, with their consequence $g_3+\alpha(h_1')=g_2+\alpha(h_2')$, do the same for the first coordinates.
\end{proof}

\begin{proposition}\label{prop:z2zb}
Every allowed pair on $\mathbb Z_2\times\mathbb Z_{2^b}$\,, $b\ge2$, has a witness.
\end{proposition}

\begin{proof}
For $b=2$ the family of Lemma~\ref{lem:family} witnesses $17$ of the $42$ allowed pairs, and their images under the symmetries (S1)--(S3), with the automorphisms $(g,h)\mapsto(g+h,h)$ and $(g,h)\mapsto(g,h+2g)$ of $\mathbb Z_2\times\mathbb Z_4$\,, are all $42$; we omit this routine enumeration. Let $b\ge3$, $B\coloneqq2^b$, $K\coloneqq\langle(1,2)\rangle$, and let $(u,v)$ be allowed. The symmetries (S2) and (S3) act on $(u_1\,,v_1)\in\mathbb Z_2^2$ as on the pair itself, so they permute its three nonzero values transitively, and they preserve the pairs trivial on $K$; so we may assume $(u_1\,,v_1)\in\{(0,0),(1,1)\}$, and that $(u,v)$ is not witnessed by the family (Corollary~\ref{cor:family}).

\emph{Pairs trivial on $K$.} The quotient by $K$ is cyclic of order $4$, with characters $0$, $v_0\coloneqq(1,B/4)$, $2v_0=(0,B/2)$, and $-v_0$\,; so a pair trivial on $K$ is one of the six of Lemma~\ref{lem:cq}.

\emph{The other pairs: a $3$-cycle.} Let $d\coloneqq7^{-1}$ and $c\coloneqq3^{-1}$ modulo $B$, let $h\in\mathbb Z_B$\,, and consider $a_1\coloneqq(0,h-d)$, $a_2\coloneqq(1,h-3d)$, $a_3\coloneqq(0,h)$ and a member $\sigma$ of the family, with $\pi$ the identity, whose $\alpha$ satisfies
\begin{equation}\label{eq:alphacond}
\alpha(h-d)=\alpha(h)+1,\qquad\alpha(h-3d+c)=\alpha(h).
\end{equation}
This is the shape of Lemma~\ref{lem:shape} for $A=2$ and the pattern $(0,1,0)$ (there $D=7^{-1}=d$, $h_3-h_2=3d$, and $h_1-h_3=-d$), whose two conditions are \eqref{eq:alphacond}. Such a $\beta$ exists: the three arguments are distinct modulo $B/2$, because their pairwise differences $-d$, $c-2d=21^{-1}$, and $c-3d=-2\cdot21^{-1}$ are nonzero modulo $B/2$ for $b\ge3$, and $\alpha$ is determined modulo $B/2$ by $\beta$ up to the known term $\lfloor h/(B/2)\rfloor$. So Lemma~\ref{lem:cycle} applies. The coefficient of $\sigma$ is $0$ by Corollary~\ref{cor:family}, so the coefficient of $\sigma'$ is the difference in Lemma~\ref{lem:cycle}, which, with $Z\coloneqq\zeta_B^{d}$ (so $\zeta_B=Z^7$), $x\coloneqq Z^{u_2}$, $y\coloneqq Z^{v_2}$, $\varepsilon_1\coloneqq(-1)^{u_1}$, and $\varepsilon_2\coloneqq(-1)^{v_1}$, evaluates to
\[
\hat\sigma'(u,v)=\frac{(-1)^{v_1\alpha(h)}\zeta_B^{(u_2+2v_2)h}}{x^3y^2}\,\Delta,\qquad
\Delta\coloneqq x^2y^2(y-x)+\varepsilon_2x^2(x-1)-\varepsilon_1y^2(y-1).
\]
(Each of $\sigma(a_1)=(\alpha(h)+1,2h-2d)$, $\sigma(a_2)=(\alpha(h),2h-6d+1)$, and $\sigma(a_3)=(\alpha(h),2h)$ is evaluated at $v$, each $a_i$ at $u$, and $\zeta_B^{-6d+1}=Z$.) It remains to find when $\Delta=0$.

\emph{Vanishing of $\Delta$.} If $(u_1\,,v_1)=(0,0)$, then $\Delta=-(x-1)(y-1)(x-y)(xy+x+y)$, and $xy+x+y\neq0$, because $|x+1|=|x|=1$ would make $x$ a primitive cube root of unity; so $\Delta=0$ forces $u=0$, $v=0$, or $u=v$. If $(u_1\,,v_1)=(1,1)$, then $\Delta=-(x-y)H$ with $H\coloneqq x^2y^2+x^2+y^2+xy-x-y$. A vanishing sum of roots of unity of order dividing $B$ splits into pairs of opposite terms: if $\zeta$ has order $B$, a relation $\sum_jc_j\zeta^j=0$ with integers $c_j$ says that $\sum_jc_jX^j$ is divisible by the minimal polynomial $X^{B/2}+1$, which, after reducing exponents modulo $B/2$, means $c_j=c_{j+B/2}$: the terms $\zeta^j$ and $\zeta^{j+B/2}=-\zeta^j$ occur equally often. Because $H(1,y)=2y^2$ and $H(x,1)=2x^2$, we have $x,y\neq1$. If $-x$ is paired with a positive term of $H$, that term is $x^2y^2$ or $y^2$ (the others would force $x=1$ or $y=1$), and then $-y$ is paired with $x^2y^2$ or $x^2$; the three possible combinations give $x^5=1$, $y^5=1$, or $y^3=1$, which is impossible. So $-x$ is paired with $-y$, that is $y=-x$, and the remaining terms $x^4,x^2,x^2,-x^2$ pair up only if $x^2=-1$. Then $x=\pm i$ and $y=\mp i$, so $u_2\,,v_2\in\{B/4,3B/4\}$ with $v_2=u_2+B/2$, and $\{u,v\}=\{v_0\,,-v_0\}$ is trivial on $K$. So every other pair is witnessed by $\sigma'$.
\end{proof}

\begin{lemma}[$\GF(4)$ designs]\label{lem:gf4d}
Let $P\coloneqq C\oplus K$ with $C\cong\mathbb Z_2^2\,$, identified with $\GF(4)$, and $K=\langle\kappa\rangle$ cyclic of even order $m$. Let $k\in\GF(4)\setminus\{0,1\}$ and $0\neq d\in\GF(4)$, let $t$ be the functional with $t(k^{-1}d)=0$ and $t((1+k)^{-1}d)=1$, and let $\ell$ be any functional. With $s(x)\coloneqq x+\ell(x)\kappa$, $\varphi_0(x)\coloneqq kx$, and $\varphi_1(x)\coloneqq kx+d$, the map $\sigma$ of Lemma~\ref{lem:ext} is a complete mapping of $P$.
\end{lemma}

\begin{proof}
The elements $k^{-1}d$ and $(1+k)^{-1}d$ are distinct and nonzero, hence a basis of $\GF(4)$ over $\GF(2)$, so $t$ exists. The cocycle of $s$ is $c(x,y)=\ell(x)+\ell(y)-\ell(x+y)\in\{0,2\}$, with $\ell$ read in $\{0,1\}$. Condition (A) says $t(x)=t(x+k^{-1}d)$, because $\varphi_1^{-1}(y)=\varphi_0^{-1}(y)+k^{-1}d$; and as the cocycle is even, (B) says that for each $y$ the two solutions of $x+\varphi_h(x)=y$, namely $x=(1+k)^{-1}y$ for $h=0$ and $x=(1+k)^{-1}(y+d)$ for $h=1$, have different values of $t$, that is, $t((1+k)^{-1}d)=1$; both hold by the choice of $t$.
\end{proof}

\begin{proposition}[rank at least $3$]\label{prop:rank3}
Let $P$ be an abelian $2$-group of rank $r\ge3$. For $e\in P\setminus2P$ put $N_e\coloneqq\langle e\rangle+2P$, so that $P/N_e\cong\mathbb Z_2^{r-1}\,$. Suppose that every allowed pair on every noncyclic $N_e$ has a witness. Then every allowed pair on $P$ has a witness.
\end{proposition}

\begin{proof}
Characters of order at most $2$ are those trivial on $2P$; we identify them with linear functionals on $P/2P\cong\mathbb Z_2^r\,$. A noncyclic $N_e$ has a complete mapping, and so does $P/N_e$\,. If $2P$ is noncyclic, every $N_e$ is noncyclic. If $2P$ is cyclic, then $P=\mathbb Z_2^{r-1}\times\mathbb Z_{2^\gamma}$ with $\gamma\ge2$, we write characters as $\chi=(\chi',\chi_\gamma)$ with $\chi(x,j)=(-1)^{\chi'(x)}\zeta_{2^\gamma}^{\chi_\gamma j}\,$, and $N_e$ for $e=(x,j)$ is noncyclic exactly when $j$ is even, in which case $N_e=\langle(x,0)\rangle\oplus2P$. Let $(u,v)$ be allowed.

\emph{(I) None of $u$, $v$, $u-v$ has order at most $2$.} All three are nontrivial on $2P\le N_e$\,, so the restriction to any noncyclic $N_e$ is allowed, and Lemma~\ref{lem:lift}(ii) applies.

\emph{(II) Exactly one has order $2$,} say $v$ after (S2) and (S3), so that $u$ has order at least $4$. If $2P$ is noncyclic, choose $e\notin2P$ with $v(e)\neq1$; the restriction to $N_e$ is allowed, and Lemma~\ref{lem:lift}(ii) applies. Let $2P$ be cyclic. If $v'\neq0$, take $e=(x,0)$ with $v'(x)\neq0$ and argue as before. If $v'=0$, then $v=(0,2^{\gamma-1})$ is trivial on every noncyclic $N_e$\,, and no automorphism of $P$ moves it, because its kernel $\mathbb Z_2^{r-1}\times2\mathbb Z_{2^\gamma}$ is the set of $a\in P$ with $2^{\gamma-1}a=0$, which every automorphism maps onto itself. We distinguish three cases.
If $r\ge4$, let $N\coloneqq\langle(x,0),(0,1)\rangle\cong\mathbb Z_2\times\mathbb Z_{2^\gamma}$ for some $x\neq0$. Then $v(0,1)=-1$, and $u$ and $u-v$ have order at least $4$ on $\langle(0,1)\rangle$, so the restriction to $N$ is allowed; it has a witness by Proposition~\ref{prop:z2zb}, $P/N\cong\mathbb Z_2^{r-2}$ is noncyclic, and Lemma~\ref{lem:lift}(ii) applies.
If $r=3$ and $u_\gamma$ is even, then $\gamma\ge3$; choose $x\neq0$ with $u'(x)=0$. Then $u$ and $v$ are trivial on $N\coloneqq\langle(x,0),(0,2^{\gamma-1})\rangle\cong\mathbb Z_2^2\,$, the induced pair on $P/N\cong\mathbb Z_2\times\mathbb Z_{2^{\gamma-1}}$ is allowed and has a witness by Proposition~\ref{prop:z2zb}, and Lemma~\ref{lem:lift}(i) applies.
If $r=3$ and $u_\gamma$ is odd, we use Lemma~\ref{lem:gf4d} with $C\coloneqq\mathbb Z_2^2\times\{0\}=\GF(4)$, $\kappa\coloneqq(0,1)$, $m=2^\gamma$, and a functional $\ell$ with $\ell(d)=1$, together with the coefficient formula of Lemma~\ref{lem:ext}. The value of $uv$ at $\kappa$ has order $2^\gamma$, so $S_0\neq0$ and $\varepsilon=-1$; moreover $v(t(x)\kappa)=(-1)^{t(x)}$, $v(s(y))=(-1)^{\ell(y)}$, and $u(s(x))=(-1)^{u'(x)}z^{\ell(x)}$ with $z\coloneqq\zeta_{2^\gamma}^{u_\gamma}$\,, which is not real. As $\ell(d)=1$, the bracket in Lemma~\ref{lem:ext} is $2(-1)^{\ell(kx)}$, so the coefficient is $2S_0\sum_x(-1)^{w(x)}z^{\ell(x)}$ with $w\coloneqq u'+t+\ell\circ k$; writing $\ker(\ell)=\{0,f\}$ and taking $a\notin\ker(\ell)$, this is $2S_0\bigl(1+(-1)^{w(f)}\bigr)\bigl(1+(-1)^{w(a)}z\bigr)$, which is nonzero when $w\in\{0,\ell\}$. Identify the functionals on $\GF(4)$ with its elements through the trace form, $x\mapsto\operatorname{Tr}(\lambda x)$ corresponding to $\lambda$, as in Example~\ref{ex:gf4}, with $u'$ corresponding to $\lambda$, and take $k$ with $k^2=k+1$, so that $k^{-1}=k^2$ and $(1+k)^{-1}=k$. Then $t$ corresponds to $kd^{-1}$, because $\operatorname{Tr}(1)=0$ and $\operatorname{Tr}(k^2)=1$, and $\ell$ corresponding to $bd^{-1}$ with $b\in\{k,k^2\}$ has $\ell(d)=\operatorname{Tr}(b)=1$. The condition $w\in\{0,\ell\}$ then reads $\lambda\in d^{-1}\{k(1+b),k+bk^2\}$, which is $\{d^{-1},k^2d^{-1}\}$ for $b=k$ and $\{k^2d^{-1},0\}$ for $b=k^2$. So $b\coloneqq k^2$ serves for $\lambda=0$, and $b\coloneqq k$ with $d\coloneqq\lambda^{-1}$ for $\lambda\neq0$.

\emph{(III) All three have order at most $2$,} that is, $u$ and $v$ are trivial on $2P$. If $2P$ is noncyclic, Lemma~\ref{lem:lift}(i) with $N=2P$ applies, the induced pair on $\mathbb Z_2^r$ being allowed with a witness by Lemma~\ref{lem:homocyclic}. Let $2P$ be cyclic, so $u_\gamma\,,v_\gamma\in\{0,2^{\gamma-1}\}$. If some nonzero $x$ lies in $\ker(u')\cap\ker(v')$, which happens when $r\ge4$ or when $u'$ and $v'$ are linearly dependent, then $u$ and $v$ are trivial on $N_e$ for $e=(x,0)$, and Lemma~\ref{lem:lift}(i) applies. Otherwise $r=3$ and $u'$ and $v'$ are independent. Then $K\coloneqq\ker(u)\cap\ker(v)$ consists of the $(x,j)$ with $u'(x)=j[u_\gamma\neq0]$ and $v'(x)=j[v_\gamma\neq0]$ modulo $2$; for even $j$ this forces $x=0$, and for odd $j$ it determines $x=x_1$\,, so $K=\langle\kappa\rangle$ with $\kappa\coloneqq(x_1\,,1)$ is cyclic of order $m=2^\gamma$, and the subgroup $\{(x,0)\}$ meets $K$ only in $0$, so that $P=\{(x,0)\}\oplus K$. We use Lemma~\ref{lem:gf4d} with $C\coloneqq\{(x,0)\}=\GF(4)$, $\ell\coloneqq0$, $k\in\GF(4)$ such that $u'(x)=v'(kx)$ for all $x$ (it exists, and $k\neq0,1$, because $u'$ and $v'$ are independent), and $d$ spanning $\ker(v')$. Because $u$ and $v$ are trivial on $K$, we have $S_0=m/2$, $\varepsilon=1$, and $v(t(x)\kappa)=1$, so the coefficient of Lemma~\ref{lem:ext} is $\frac m2\sum_x(-1)^{u'(x)}\bigl((-1)^{v'(kx)}+(-1)^{v'(kx+d)}\bigr)=\frac m2\cdot4\cdot2\neq0$, because $u'(x)=v'(kx)$ and $v'(d)=0$.
\end{proof}

\begin{proposition}\label{prop:rank2}
Every allowed pair on $\mathbb Z_{2^a}\times\mathbb Z_{2^b}$\,, $2\le a<b$, has a witness.
\end{proposition}

\begin{proof}
By induction on $a$. Let $A\coloneqq2^a$, $B\coloneqq2^b$, $L\coloneqq B/A\ge2$, $P\coloneqq\mathbb Z_A\times\mathbb Z_B$\,, $c\coloneqq(A+1)^{-1}$, and $D\coloneqq(A^2+A+1)^{-1}$ modulo $B$. The subgroup $2P\cong\mathbb Z_{A/2}\times\mathbb Z_{B/2}$ is noncyclic, every allowed pair on it has a witness (by Proposition~\ref{prop:z2zb} if $a=2$, and by induction if $a\ge3$), and $P/2P\cong\mathbb Z_2^2\,$. Let $(u,v)$ be allowed. If none of $u$, $v$, $u-v$ has order at most $2$, all three are nontrivial on $2P$, the restriction to $2P$ is allowed, and Lemma~\ref{lem:lift}(ii) applies. If all three have order at most $2$, then $u$ and $v$ are trivial on $2P$, and Lemma~\ref{lem:lift}(i) applies, the induced pair on $\mathbb Z_2^2$ having a witness by Lemma~\ref{lem:homocyclic}. In the remaining case exactly one of them has order $2$, and by (S2) and (S3) we may assume that it is $v$; so $v\in\{(A/2,0),(0,B/2),(A/2,B/2)\}$, and $u=(u_1\,,u_2)$ has order at least $4$.

\emph{The family.} By Corollary~\ref{cor:family}, either some member of the family of Lemma~\ref{lem:family} witnesses $(u,v)$, or every member has coefficient $0$ at $(u,v)$; in the latter case we call the pair uncovered.

\emph{The $3$-cycle.} Put $p_2\coloneqq D(A^2+A)/2-1$ and $p_3\coloneqq DA^2/2-1$, which are odd because $A(A+1)/2$ and $A^2/2$ are even, and let $h\in\mathbb Z_B$\,, $a_1\coloneqq(0,h)$, $a_2\coloneqq(A/2,h+p_2)$, $a_3\coloneqq(A/2+1,h+p_3)$; this is the shape of Lemma~\ref{lem:shape} for the pattern $(0,A/2,A/2+1)$, whose conditions read
\[
\alpha(y_2)=\alpha(y_1)+1\qquad\text{and}\qquad\alpha(y_3)=\alpha(y_1)-A/2,
\]
where $y_1\coloneqq h$, $y_2\coloneqq h+p_2+cA/2$, and $y_3\coloneqq h+p_3+c(A/2+1)$.
A computation gives $y_2-y_1=N_2\,cD$, $y_3-y_2=\frac{A^2+A+2}2\,cD$, and $y_3-y_1=-\frac A2\,cD$, with $N_2\coloneqq\frac A2(A+1)^2-(\frac A2+1)(A^2+A+1)$. Both $N_2$ and $(A^2+A+2)/2$ are odd, so $y_2-y_1$ and $y_3-y_2$ are units and the residues $r_i\coloneqq y_i\bmod L$ satisfy $r_2\neq r_1$ and $r_2\neq r_3$\,. Write $\ell_i\coloneqq\lfloor y_i/L\rfloor$, so that $\alpha(y_i)=\pi(\ell_i)+\beta(r_i)$. If $r_3\neq r_1$\,, any $\pi$ will do, and $\beta(r_2)$ and $\beta(r_3)$ are chosen to meet the two conditions. If $r_3=r_1$\,, then $\ell_3\neq \ell_1$ because $y_3\neq y_1$ modulo $B$; choose $\pi$ with $\pi(\ell_3)=\pi(\ell_1)-A/2$ and then $\beta(r_2)$. In both cases Lemma~\ref{lem:cycle} applies. For an uncovered pair the coefficient of the modified mapping is the difference in that lemma. Here $\sigma(a_1)=(\alpha_0\,,Ah)$, $\sigma(a_2)=(\alpha_0+1,A(h+p_2)+A/2)$, and $\sigma(a_3)=(\alpha_0-A/2,A(h+p_3)+A/2+1)$ with $\alpha_0\coloneqq\alpha(h)$; because $v_1\in\{0,A/2\}$ and $4\mid A$, $v$ takes the values $\zeta_A^{v_1\alpha_0}\,$, $\zeta_A^{v_1\alpha_0}\delta$, $\zeta_A^{v_1\alpha_0}\varepsilon$ at these, where $\delta\coloneqq(-1)^{[v_1\neq0]}$ and $\varepsilon\coloneqq(-1)^{[v_2\neq0]}$, while $u$ takes the values $\zeta_B^{u_2h}$ times $1$, $(-1)^{u_1}x^{p_2}$, $(-1)^{u_1}\eta x^{p_3}$ at $a_1$\,, $a_2$\,, $a_3$\,, where $x\coloneqq\zeta_B^{u_2}$ and $\eta\coloneqq\zeta_A^{u_1}\,$. So the coefficient is $\zeta_A^{v_1\alpha_0}\zeta_B^{u_2h}\Delta$ with
\[
\Delta\coloneqq(\delta-1)+(-1)^{u_1}x^{p_2}(\varepsilon-\delta)+(-1)^{u_1}\eta\,x^{p_3}(1-\varepsilon),
\]
and every uncovered pair with $\Delta\neq0$ has a witness.

\emph{The zeros of $\Delta$.} For $v=(A/2,0)$, $\Delta=2\bigl((-1)^{u_1}x^{p_2}-1\bigr)=2\bigl(u(g_v)-1\bigr)$ with $g_v\coloneqq(A/2,p_2)$. For $v=(A/2,B/2)$, $\Delta=2\bigl((-1)^{u_1}\eta x^{p_3}-1\bigr)=2\bigl(u(g_v)-1\bigr)$ with $g_v\coloneqq(A/2+1,p_3)$. For $v=(0,B/2)$, $\Delta=2(-1)^{u_1}x^{p_2}\bigl(\eta x^{p_3-p_2}-1\bigr)=2(-1)^{u_1}x^{p_2}\bigl(u(g_v)-1\bigr)$ with $g_v\coloneqq(1,p_3-p_2)$. So $\Delta=0$ exactly when $u$ is trivial on $g_v$\,. In each case $v(g_v)=1$, because $\zeta_A^{A^2/4}=1$, $p_2$ and $p_3$ are odd, and $p_3-p_2$ is even. Moreover $g_v\notin2P$, because it has an odd coordinate, so $P/\langle g_v\rangle$ is cyclic, as $P$ has rank $2$, and the characters trivial on $g_v$ form a cyclic group, which contains $u$, of order $2^j\ge4$, and $v$, of order $2$; hence $v=2^{j-1}u$. Now Lemma~\ref{lem:cq} applies with $K\coloneqq\langle g_v\rangle$. This group is cyclic of even order, because $g_v$ has an odd coordinate; $P/K$ is cyclic of order at least $2^j\ge4$, because $u$ factors through it; and for $e\coloneqq(1,0)$ if the second coordinate of $g_v$ is odd and $e\coloneqq(0,1)$ otherwise, the image of $e$ generates $P/K$, and the section $n\mapsto ne$ ($0\le n<|P/K|$) has cocycle values in $2K$. Indeed, in the first case $\langle e\rangle\cap K=0$, and in the second $\langle e\rangle\cap K=\langle Ag_v\rangle$, so in both cases $|P/K|\,e$ is an even multiple of $g_v$\,.
\end{proof}

\begin{theorem}\label{thm:twogroups}
Every allowed pair on every noncyclic abelian $2$-group has a witness.
\end{theorem}

\begin{proof}
By induction on the order. A noncyclic abelian $2$-group of rank $2$ is homocyclic, or $\mathbb Z_2\times\mathbb Z_{2^b}$\,, or $\mathbb Z_{2^a}\times\mathbb Z_{2^b}$ with $2\le a<b$, and Lemma~\ref{lem:homocyclic}, Proposition~\ref{prop:z2zb}, and Proposition~\ref{prop:rank2} apply. For rank at least $3$, Proposition~\ref{prop:rank3} and induction apply.
\end{proof}

\emph{Coupling with $3$.} Lemma~\ref{lem:two} couples any group with a complete mapping to a cyclic $r$-group with $r\ge5$. For $r=3$ quadratic residues are useless, and we need a different design, which we give for $2$-groups.

\begin{lemma}[three layers over $\mathbb Z_3$]\label{lem:layers}
Let $H$ be an abelian group with a complete mapping $\rho$, let $G\coloneqq H\oplus\mathbb Z_3$\,, and let $\omega\coloneqq e^{2\pi i/3}$. Let $u$ be a nontrivial character of $H$ and $v(i,h)\coloneqq\omega^h$, both extended to $G$. For $t\in H$ put $\theta\coloneqq\mathrm{id}+\rho$, $\pi(i)\coloneqq\rho^{-1}(\rho(i)+t)$, $\pi'(i)\coloneqq\theta^{-1}(\theta(i)+t)$, $\tau\coloneqq\pi'\pi^{-1}$, and for an orbit $O$ of $\tau$ on $H$ put $S_O\coloneqq\sum_{i\in O}\bigl(u(i)-u(\pi^{-1}(i))\bigr)$. Then for every orbit $O$ there is a complete mapping $\sigma$ of $G$ with $\hat\sigma(u,v)=(1-\omega)S_O$\,. In particular $(u,v)$ has a witness whenever $S_O\neq0$ for some $t$ and some $O$.
\end{lemma}

\begin{proof}
For functions $x,y,z\colon H\to\mathbb Z_3$ put
\[
\sigma(i,0)\coloneqq(\rho(i),x_i),\qquad\sigma(i,1)\coloneqq(\rho(i),y_i),\qquad\sigma(i,2)\coloneqq(\rho(i)-t,z_i).
\]
Counting preimages, $\sigma$ is a bijection if and only if $\{x_i\,,y_i\,,z_{\pi(i)}\}=\mathbb Z_3$ for all $i$, and $a\mapsto a+\sigma(a)$ is a bijection if and only if $\{x_i\,,1+y_i\,,2+z_{\pi'(i)}\}=\mathbb Z_3$ for all $i$. Both hold if $y_i=x_i+1$ and $z_{\pi(i)}=z_{\pi'(i)}=x_i+2$ for all $i$, that is, if $z$ is constant on the orbits of $\tau$ and $x_i\coloneqq z_{\pi(i)}-2$, $y_i\coloneqq z_{\pi(i)}-1$. With $\omega^{-2}+\omega^{-1}=-1$ the coefficient is then
\[
\hat\sigma(u,v)=\sum_iu(i)\bigl(\omega^{x_i}+\omega^{y_i}+\omega^{z_i}\bigr)=\sum_i\omega^{z_i}\bigl(u(i)-u(\pi^{-1}(i))\bigr)=\sum_{O'}\omega^{z_{O'}}S_{O'}\,,
\]
where $z_{O'}$ is the value of $z$ on the orbit $O'$. Because $\sum_{O'}S_{O'}=\sum_iu(i)-\sum_iu(\pi^{-1}(i))=0$, taking $z\coloneqq0$ on $O$ and $z\coloneqq1$ on the other orbits gives $\hat\sigma(u,v)=S_O-\omega S_O$\,.
\end{proof}

\begin{lemma}[coupling $2$ with $3$]\label{lem:twothree}
Let $G\coloneqq P\oplus\mathbb Z_3$ with $P\cong\mathbb Z_2^m\,$, $m\ge2$, let $u$ and $v$ be nontrivial characters of $P$ and $\mathbb Z_3$\,, extended trivially to $G$. Then $(u,v)$ has a witness.
\end{lemma}

\begin{proof}
Applying an automorphism of $\mathbb Z_3$ we may take $v(h)=\omega^h$, and we apply Lemma~\ref{lem:layers}. Identify $P$ with the additive group of $\GF(2^m)$, and for $k\in\GF(2^m)\setminus\{0,1\}$ let $\rho(i)\coloneqq ki$, a complete mapping because $k$ and $k+1$ are nonzero. Then $\pi(i)=i+k^{-1}t$ and $\pi'(i)=i+(1+k)^{-1}t$, so $\tau$ is the translation by $t''\coloneqq(k(1+k))^{-1}t$, its orbits are the cosets $\{i,i+t''\}$, and $S_O=\bigl(1-u(k^{-1}t)\bigr)\bigl(u(i)+u(i+t'')\bigr)$ for $O=\{i,i+t''\}$. So $S_O\neq0$ if $u(t'')=1$ and $u(k^{-1}t)\neq1$, that is, if $t''\in\ker(u)$ and $(1+k)t''\notin\ker(u)$, because $k^{-1}t=(1+k)t''$. Choose $t''\neq0$ in $\ker(u)$, which has order $2^{m-1}\ge2$, and $y\notin\ker(u)$; then $y=k't''$ for some $k'\neq0$, and $k'\neq1$, so $k\coloneqq k'+1\notin\{0,1\}$ satisfies both conditions, and $t\coloneqq k(1+k)t''$ completes the construction.
\end{proof}

\begin{lemma}[coupling $\mathbb Z_2\times\mathbb Z_{2^b}$ with $3$]\label{lem:twothreeB}
Let $G\coloneqq P\oplus\mathbb Z_3$ with $P\coloneqq\mathbb Z_2\times\mathbb Z_B$\,, $B\coloneqq2^b\ge4$, let $u$ and $v$ be nontrivial characters of $P$ and $\mathbb Z_3$\,, extended trivially to $G$. Then $(u,v)$ has a witness.
\end{lemma}

\begin{proof}
As in Lemma~\ref{lem:twothree} we may take $v(h)=\omega^h$, and we apply Lemma~\ref{lem:layers} with the complete mapping $\rho$ of Lemma~\ref{lem:family} for $a=1$, $\pi=\mathrm{id}$, and $\beta=0$. Let $H\coloneqq B/2$, $c\coloneqq3^{-1}$ modulo $B$, $\alpha(x)\coloneqq[x\ge H]$ for $x\in\{0,\dots,B-1\}$, and $\zeta\coloneqq\zeta_B$\,. We write the elements of $P$ as $(g,x)$ with $x\coloneqq h+cg$, so that $(g,x)$ stands for $(g,x-cg)$; in these coordinates, because $1-2c\equiv c$,
\[
\rho(g,x)=\bigl(\alpha(x),\,2x+cg\bigr),\qquad\theta(g,x)=\bigl(g+\alpha(x),\,3x\bigr),\qquad u(g,x)=\lambda^g\zeta^{u_2x},
\]
where the second coordinates of $\rho$ and $\theta$ are in the original coordinate $h$, where $u=(u_1\,,u_2)$ with $u_2\neq0$ or $u_1=1$, and where $\lambda\coloneqq(-1)^{u_1}\zeta^{-cu_2}$. Consequently $\theta^{-1}(e,h')=(e-\alpha(ch'),ch')$, and $\rho^{-1}(e,h')$ is the $(g,x)$ with $g\equiv h'$ modulo $2$, $2x\equiv h'-cg$ modulo $B$, and $\alpha(x)=e$.

\emph{Case $u_1=1$ and $u_2$ even: $t\coloneqq(0,1)$.} Here $\pi(0,y)=(1,x)$ with $x\equiv y+(1-c)/2$ modulo $H$ and $\alpha(x)=\alpha(y)$; as $(3c-1)/2\equiv0$ modulo $H$, this says $\pi^{-1}(1,x)=(0,\hat x)$ with $\hat x\equiv x-c$ modulo $H$ and $\alpha(\hat x)=\alpha(x)$. Also $\pi'(g,x)=(g+\alpha(x)-\alpha(x+c),x+c)$. Choose $x$ with $\alpha(x-c)\neq\alpha(x)$, which exists because translation by $c$ does not preserve $\{0,\dots,H-1\}$. Then $\hat x=x-c+H$, so $\tau(1,x)=\pi'(0,\hat x)=(1,x+H)$, and likewise $\tau(1,x+H)=(1,x)$: $O\coloneqq\{(1,x),(1,x+H)\}$ is an orbit of $\tau$. Because $u_2$ is even, $\zeta^{u_2H}=1$, so $u(1,x)=u(1,x+H)=-\zeta^{u_2(x-c)}$ and $u(0,\hat x)=u(0,\hat x+H)=\zeta^{u_2(x-c)}$, and $S_O=-4\zeta^{u_2(x-c)}\neq0$.

\emph{The other cases: $t\coloneqq(0,2)$.} Now $\pi(g,x)=(g,r(x))$, where $r$ rotates each of $\{0,\dots,H-1\}$ and $\{H,\dots,B-1\}$ cyclically by one step, and $\pi'(g,x)=(g+\alpha(x)-\alpha(x+2c),x+2c)$. Because $2c-1\equiv-c$, $\tau(g,x)=(g+\alpha(x)-\alpha(\tau_1(x)),\tau_1(x))$ with $\tau_1(x)\coloneqq x-c+H[x\in\{0,H\}]$. So $g+\alpha(x)$ is invariant under $\tau$, and the orbits of $\tau$ are the sets $O_X^\epsilon\coloneqq\{(\epsilon-\alpha(x),x) ~:~ x\in X\}$ for the orbits $X$ of $\tau_1$ and $\epsilon\in\{0,1\}$. The orbit of $0$ under $\tau_1$ is $X_0\coloneqq\{0\}\cup\{H-jc ~:~ 1\le j<H\}$, because $H-jc\in\{0,H\}$ forces $j\equiv0$ or $j\equiv H$ modulo $B$ ($c$ is a unit and $3H\equiv H$), and $H-Hc\equiv0$ as $c$ is odd; $X_0$ meets each residue class modulo $H$ exactly once. Write $f(g,x)\coloneqq u(g,x)-u(\pi^{-1}(g,x))=\lambda^g\zeta^{u_2x}(1-\zeta^{-u_2})$ for $x\notin\{0,H\}$ and $f(g,x)=\lambda^g\zeta^{u_2x}(1-(-1)^{u_2}\zeta^{-u_2})$ for $x\in\{0,H\}$.

If $u_2$ is odd, take $z$ equal to $0$ on $O_{X_0}^0\cup O_{X_0}^1$ and $1$ elsewhere in the proof of Lemma~\ref{lem:layers}; the coefficient is then $(1-\omega)(S_{O^0}+S_{O^1})$, and with $a\coloneqq\zeta^{-cu_2}$ and $a'\coloneqq\zeta^{-u_2}$,
\begin{align*}
S_{O^0}+S_{O^1}&=\sum_{x\in X_0}\sum_{g}f(g,x)=(1+\lambda)\Bigl[(1-a')\sum_{x\in X_0}\zeta^{u_2x}+2a'\Bigr]\\
&=(1+\lambda)\Bigl[(1-a')\Bigl(2+\frac2{a-1}\Bigr)+2a'\Bigr]=\frac{2(1+\lambda)(a-a')}{a-1}\,,
\end{align*}
using $\sum_{x\in X_0}\zeta^{u_2x}=1-\sum_{j=1}^{H-1}a^j=2-\sum_{j=0}^{H-1}a^j=2+\frac2{a-1}$ (as $a^H=-1$). This is nonzero: $\lambda=\pm\zeta^{-cu_2}$ is not real, and $a/a'=\zeta^{2cu_2}\neq1$.

If $u_1=0$ and $u_2$ is even, then $\lambda=\zeta^{-cu_2}\neq1$, $\zeta^{u_2H}=1$, and the second formula for $f$ agrees with the first; as the elements of $O_{X_0}^1$ are the $(1-\alpha(x),x)$, $x\in X_0$\,,
\begin{align*}
S_{O_{X_0}^1}&=\lambda(1-\zeta^{-u_2})\sum_{x\in X_0}\lambda^{-\alpha(x)}\zeta^{u_2x}\\
&=\lambda(1-\zeta^{-u_2})\Bigl[\sum_{y\in\mathbb Z_H}\zeta^{u_2y}+(\lambda^{-1}-1)\sum_{y\in D}\zeta^{u_2y}\Bigr],
\end{align*}
where $D$ is the set of classes $y$ modulo $H$ whose representative in $X_0$ lies in $\{H,\dots,B-1\}$; here $\zeta^{u_2y}$ depends only on $y$ modulo $H$. The first sum vanishes because $\zeta^{u_2}$ is a nontrivial $H$-th root of unity. For $1\le y<H$ the representative is $H-jc$ with $jc\equiv-y$ modulo $H$ and $1\le j<H$, so $jc$ modulo $B$ is the element $w\in\{H-y,B-y\}$ with $3w$ modulo $B$ less than $H$, and the representative lies in $\{H,\dots,B-1\}$ exactly when $w=B-y$, that is, when $3y$ modulo $B$ exceeds $H$. As $3y<B+H$, this says $H<3y<B$: $D$ is the set of integers strictly between $H/3$ and $2H/3$, which has $(H-1)/3$ or $(H+1)/3$ elements, an odd number because $H$ is even. So $\sum_{y\in D}\zeta^{u_2y}=\zeta^{u_2y_0}(\xi^{|D|}-1)/(\xi-1)$ with $\xi\coloneqq\zeta^{u_2}\neq1$ of $2$-power order, which does not divide $|D|$; the sum is nonzero, and so is $S_{O_{X_0}^1}$\,.
\end{proof}

\begin{proposition}\label{prop:coupling3}
Let $P$ be a noncyclic abelian $2$-group, $C'$ cyclic of order $3^f$, $u$ a nontrivial character of $P$, and $v$ a character of $C'$ of order $3$. Then $(u,v)$ has a witness on $P\oplus C'$.
\end{proposition}

\begin{proof}
The subgroup $N\coloneqq3C'$ has odd order, $u$ and $v$ are trivial on it, $(P\oplus C')/N\cong P\oplus\mathbb Z_3$ has a complete mapping, and the induced pair is allowed; so by Lemma~\ref{lem:lift}(i) it suffices to treat $C'=\mathbb Z_3$\,. Then we argue by induction on $|P|$. If $P$ is elementary abelian, this is Lemma~\ref{lem:twothree}. If $2P$ is noncyclic and $u$ has order at least $4$, then $u$ is nontrivial on $2P$, the restriction of $(u,v)$ to $N\coloneqq2P\oplus C'$ is allowed and has a witness by induction, and $(P\oplus C')/N\cong P/2P$ has a complete mapping, so Lemma~\ref{lem:lift}(ii) applies. If $2P$ is noncyclic and $u$ has order $2$, then $u$ and $v$ are trivial on $N\coloneqq2P$, which has a complete mapping, and Lemma~\ref{lem:lift}(i) applies, the induced pair on $(P/2P)\oplus C'$ having a witness by Lemma~\ref{lem:twothree}. If $2P$ is cyclic and nontrivial, then $P\cong\mathbb Z_2^{r-1}\times\mathbb Z_{2^\gamma}$ with $r\ge2$ and $\gamma\ge2$; for $r=2$ this is Lemma~\ref{lem:twothreeB}, and for $r\ge3$ we use the noncyclic subgroups $N_e=\langle(x,0)\rangle\oplus2P$ of Proposition~\ref{prop:rank3}. If $u$ has order at least $4$, it is nontrivial on $2P\le N_e$\,, and Lemma~\ref{lem:lift}(ii) applies with $N\coloneqq N_e\oplus C'$, by induction. If $u$ has order $2$, choose $x\neq0$ with $u(x,0)=1$, possible because $r-1\ge2$; then $u$ and $v$ are trivial on $N_e$\,, and Lemma~\ref{lem:lift}(i) applies, the induced pair on $\mathbb Z_2^{r-1}\oplus C'$ having a witness by Lemma~\ref{lem:twothree}.
\end{proof}

\subsection{All abelian groups}\label{sec:main}

\begin{theorem}\label{thm:main}
If $G$ is an abelian group with a transversal, that is, whose Sylow $2$-subgroup is trivial or not cyclic, then the transversal polytope of its Cayley table has dimension ${(n-1)(n-2)}$.
\end{theorem}

\begin{proof}
Write $G=\bigoplus_pG_p$ as the sum of its $p$-parts, the \emph{components}, and let $(u,v)$ be allowed. Every component has a complete mapping, the identity if $p$ is odd and by Theorem~\ref{thm:twogroups} if $p=2$, and every allowed pair on a component has a witness, by Lemma~\ref{lem:pgroups} and Theorem~\ref{thm:twogroups}. The restriction $(u_p\,,v_p)$ of $(u,v)$ to $G_p$ is zero, allowed, or excluded of one of the three types, and $\hat\sigma$ is multiplicative over the components.
\emph{(A) No restriction is excluded.} The product of a witness on each component with an allowed restriction and any complete mapping on the others is a witness: each factor of its coefficient is nonzero or equal to $|G_p|$.
\emph{(B) Some restriction is excluded and some is allowed.} Let $H$ be the sum of the components with an excluded restriction and $K$ the sum of the others. A product of complete mappings on the components of $H$ has coefficient $0$, because it has a factor $\hat\sigma_p(u_p\,,v_p)$ of an excluded pair; the product on $K$ of (A) is a witness, with $v_K\neq0$; so a twisted product (Lemma~\ref{lem:twist}) is a witness.
\emph{(C) Every restriction is excluded or zero.} Because $(u,v)$ is allowed, the excluded restrictions do not all have the same type. If some $G_r$ with $r\ge5$ has an excluded restriction, choose another component $G_p$ whose excluded restriction has a different type; otherwise the components with excluded restrictions are $G_2$ and $G_3$\,, and we put $p\coloneqq2$, $r\coloneqq3$. The symmetries permute the types on all components simultaneously, as the full symmetric group of the three types, so we may assume rows on $G_p$ ($u_p\neq0$, $v_p=0$) and columns on $G_r$ ($u_r=0$, $v_r\neq0$). Choose $g\in G_r$ with $v(g)\neq1$, let $r^j$ be the order of $v(g)$, and let $C'\coloneqq\langle r^{j-1}g\rangle$, a cyclic $r$-group on which $v$ has order $r$. Let $\rho$ be a product of complete mappings of the components that fix $0$ (a translate of a complete mapping is one), the identity on $G_r$\,, and $S\coloneqq G_p\oplus C'$. Replace $\rho$ on $S$ by the complete mapping of Lemma~\ref{lem:two} (with the restriction of $\rho$ to $G_p$) if $r\ge5$, and of Proposition~\ref{prop:coupling3} if $r=3$; both map $S$ onto $S$ with symbols $S$, so the result is a complete mapping of $G$. Its coefficient is that of $\rho$, which is $0$ because it has the factor $\sum_{i\in G_p}u_p(i)=0$, minus the contribution $\sum_{(i,h)\in S}u_p(i)v(h)=0$ of $S$, plus the nonzero coefficient of the replacement.
\end{proof}

\begin{example}[the pairs of $\mathbb Z_{15}$]\label{ex:z15}
The cases of the proof all occur in $\mathbb Z_3\times\mathbb Z_5\cong\mathbb Z_{15}$\,.
(A) The pair $((1,1),(2,2))$ is allowed on both components. A witness on $\mathbb Z_3$ is the identity ($1+2\equiv0$), and one on $\mathbb Z_5$ is $x\mapsto2x$ ($1+2\cdot2\equiv0$); their product is a witness for the pair, with $|\hat\sigma|=3\cdot5=15$. The pair $(5,10)$ in the cyclic notation of Example~\ref{ex:lift} restricts to the allowed pair $(2,1)$ on $\mathbb Z_3$ and to $(0,0)$ on $\mathbb Z_5$\,, so the identity is a witness.
(B) The pair $((1,1),(1,2))$ is excluded of type symbols on $\mathbb Z_3$ and allowed on $\mathbb Z_5$: Example~\ref{ex:twist}.
(C) The pair $((1,0),(0,1))$ has type rows on $\mathbb Z_3$ and columns on $\mathbb Z_5$: Example~\ref{ex:qr}.
\end{example}

\subsection{Nonabelian groups and general Latin squares}\label{sec:general}

Nonabelian groups would need a coherent-configuration version of the method (see~\cite{LeeLee}), because the orbits of the symmetries $(a,b)\mapsto(xa,by)$ on pairs of cells no longer form a commutative association scheme (the associated permutation representation is not multiplicity-free). By computation, every
group of order at most $21$ whose Cayley table has a transversal (by Hall and Paige~\cite{HP}, and the completion of
their conjecture recorded in~\cite{EMM}, those with trivial or noncyclic Sylow $2$-subgroups; $35$ groups, including
all $13$ of order $16$), and also $\mathbb Z_5^2\,$, attains ${(n-1)(n-2)}$: for each group, the incidence vectors of a few
transversals found by randomized search, together with their images under $(a,b)\mapsto(xa,by)$, already span an affine
space of that dimension.

For general Latin squares there are two simple obstructions: the dimension is at most $|T|-1$, where $T$ is the set of
transversals, and every cell in no transversal adds the equation $x_e=0$. In random samples of orders $5$ to $7$,
about four fifths of the squares that fall short attain exactly the lesser of these two bounds; the rest have further
equations. From order $8$ on, all sampled squares were full-dimensional ($25$, $20$, $15$, and $10$ squares of
orders $8$, $9$, $10$, and $11$). The following statement
explains why many transversals with evenly spread pair counts suffice.

\begin{proposition}\label{prop:stability}
Let $L$ be a Latin square of order $n$ with transversal set $T\neq\emptyset$. For a cell $e$ let $\tau_e$ be the number of
transversals through $e$, and for cells $e,f$ in different rows, columns, and symbols let $\tau_{ef}$ be the number through
both; let $\bar\tau$ be the mean of the $\tau_{ef}$\,. If $|\tau_{ef}-\bar\tau|\le\varepsilon\bar\tau$ and $|\tau_e-(n-2)\bar\tau|\le\varepsilon(n-2)\bar\tau$ for all such
$e,f$, with $\varepsilon(1+\varepsilon)(n-2)<1$, then the transversal polytope has dimension ${(n-1)(n-2)}$.
\end{proposition}

\begin{proof}
Let $W$ be the orthogonal complement of the line indicators, so $\dim(W)=(n-1)(n-2)$, and differences of transversals
lie in $W$. It suffices that no nonzero $w\in W$ has $x\cdot w$ constant on $T$, that is, $w^\top Bw>0$ for
$B\coloneqq\sum_{x\in T}(x-\bar x)(x-\bar x)^\top=M-|T|\bar x\bar x^\top$, with $M\coloneqq\sum_xxx^\top$. Counting pairs gives $\overline{\tau_e}=(n-2)\bar\tau$ and
$|T|=\sum_e\tau_e/n$. The averaged matrix $\bar M\coloneqq(n-2)\bar\tau I+\bar\tau O$, with $O$ the adjacency matrix of the relation
``different row, column, and symbol'', acts on $W$ as $n\bar\tau$ (Theorem~\ref{thm:amorphic} for the $3$-net scheme).
By the row-sum bound, $\|M-\bar M\|\le\varepsilon n(n-2)\bar\tau$, and $|T|(\bar x\cdot w)^2=(\sum_e(\tau_e-(n-2)\bar\tau)w_e)^2/|T|\le
\varepsilon^2n(n-2)\bar\tau\,\|w\|^2$. So
\[
w^\top Bw\ge n\bar\tau\bigl(1-\varepsilon(1+\varepsilon)(n-2)\bigr)\|w\|^2>0.
\]
\end{proof}

The hypothesis is strong and far from necessary: in random squares of orders $8$ to $10$ some pairs of cells lie in no
common transversal, so $\varepsilon\ge1$, although these polytopes are full-dimensional. The same argument gives a stability
version of~\eqref{eq:principle} for any scheme: a point set whose pair counts are
close enough to constant on each relation has at least the dimension predicted by the averaged counts.

\section{Almost amorphic schemes and symmetric matrices}\label{sec:almost}

The paper of van Dam, Koolen, and Xiong~\cite{vDKX} asks how far from amorphic a scheme can be while most of its relations remain
of Latin-square type. Earlier, the same authors established that a scheme with at most one relation not of (negative) Latin-square type is
amorphic; the paper constructs non-amorphic schemes with exactly two such relations, and proves several structural statements about
the eigenvalues of strongly regular relations inside a scheme. Two of its results enter here, as applications of \S\ref{sec:schemes}; \S\S\ref{sec:affine}--\ref{sec:latin}
do not depend on them. \S\ref{sec:symm} then treats symmetric matrices directly, without that paper: a count of points and lines on quadrics gives the dimension of the polytopes of large sets of symmetric $3\times3$ and $4\times4$ matrices with pairwise differences of rank at least $3$, for every $q$.

\subsection{Marginals}

The eigenspace of the lattice graph $L_2(n)$ for the eigenvalue $n-2$ consists of the functions
$(i,j)\mapsto h_i+v_j$ with $\sum_ih_i=\sum_jv_j=0$, of dimension $2n-2$ (these are the \emph{marginal} functions of an $n\times n$
array). A lemma of~\cite{vDKX} (Lemma~7.3) says that, if $L_2(n)$ is one relation of a scheme on the cells, this eigenspace is a single
primitive idempotent of the scheme: it cannot be split, however fine the scheme is.

\begin{proposition}\label{prop:lattice}
Let $X\subseteq\R^{[n]\times[n]}$ be concordant with an association scheme having the lattice graph $L_2(n)$ as one
of its relations. Then either every $x\in X$ has all row sums equal and all column sums equal, or $\aff(X)$ contains all
$2n-2$ dimensions of marginal functions.
\end{proposition}

\begin{proof}
By Lemma~7.3 of~\cite{vDKX}, the eigenspace of $L_2(n)$ for $n-2$ is the image of a single idempotent $E$ of the scheme, and the
statement is the dichotomy $\beta=0$ or $\beta\neq0$ for $E$ in Proposition~\ref{prop:principle}.
\end{proof}

Without the lattice graph as a single relation (for instance when rows and columns are different relations) the statement
fails, which is why it is phrased this way.

\subsection{Common fissions}

A scheme $\mathcal C$ is a \emph{fission} of $\mathcal A$ if every relation of $\mathcal A$ is a union
of relations of $\mathcal C$; then every idempotent of $\mathcal A$ is a sum of idempotents of $\mathcal C$. Let
$\mathcal A\coloneqq\{I,A_1\,,\dots,A_e\,,B\}$ be a scheme on $n^2$ points with idempotents $E_1\,,\dots,E_e\,,F$, where $B$ is strongly regular of
(negative) Latin-square type, of valency other than $(n^2-1)/2$, with spectral idempotent $F$ (the one for the eigenvalue $n-t$).
Van Dam, Koolen, and Xiong establish in~\cite{vDKX} (Lemma~4.1 and Theorem~4.2) that if $B=B_1+\dots+B_f$ with each $B_j$ strongly regular of the same type, with valency
$t_j(n-1)$, then $\mathcal C=\{I,A_1\,,\dots,A_e\,,B_1\,,\dots,B_f\}$ is again a scheme, with idempotents $E_1\,,\dots,E_e\,,F_1\,,\dots,F_f$\,,
where $F=\sum_jF_j$ and $F_j$ has rank $t_j(n-1)$. Put $t\coloneqq\sum_jt_j$\,.

\begin{proposition}\label{prop:fission}
Let $X$ be concordant with $\mathcal C$, with counts $\tau_0$\,, $\tau_{A_a}$\,, and $\tau_{B_j}$\,. Let $\alpha_{ia}$ and $\alpha_{Fa}$ be the
eigenvalues of $A_a$ on $E_i$ and on $F$ in $\mathcal A$, and $\bar\tau_B\coloneqq\frac1t\sum_jt_j\tau_{B_j}$\,. Then
\[
\dim(\conv(X))=\dim(\mathbf 1^\top X)+\sum\bigl\{\operatorname{rk}(E_i) ~:~ \lambda_i\neq0\bigr\}+(n-1)\sum\bigl\{t_j ~:~ \tau_{B_j}\neq\tau_B^*\bigr\},
\]
where $\lambda_i\coloneqq\tau_0+\sum_a\tau_{A_a}\alpha_{ia}-t\bar\tau_B$ and $\tau_B^*\coloneqq\frac1n\bigl(t\bar\tau_B-\tau_0-\sum_a\tau_{A_a}\alpha_{Fa}\bigr)$.
\end{proposition}

\begin{proof}
In the amorphic scheme $\{I,A,B_1\,,\dots,B_f\}$ that occurs in the proof of~\cite{vDKX}, Theorem~4.2 (with $A=A_1+\dots+A_e$), the matrix $B_j$ has
eigenvalue $-t_j$ on $E=\sum_iE_i$ and $n\delta_{jl}-t_j$ on $F_l$\,. Because $E_i\le E$ and $F_l\le F$, the eigenvalue of $\sum_xxx^\top$ on $E_i$ is
$\lambda_i$\,, and on $F_l$ it is $\tau_0+\sum_a\tau_{A_a}\alpha_{Fa}-t\bar\tau_B+n\tau_{B_l}=n(\tau_{B_l}-\tau_B^*)$. Apply
Proposition~\ref{prop:principle}.
\end{proof}

So the computation splits: the coarse scheme $\mathcal A$ is used with $\tau_B$ replaced by the weighted average $\bar\tau_B$\,, and
the fission adds a threshold sum over the pieces $B_j$\,.

\emph{An example.} Let $q$ be even, and let $V$ be the space of symmetric $3\times3$ matrices over $\GF(q)$, so $|V|=n^2$ with
$n\coloneqq q^3$. In the scheme of the symmetric bilinear forms graph (matrices adjacent when their difference has rank~$1$), the
distance-$3$ relation $B$ is strongly regular of Latin-square type with $t=q^3-q^2+1$; van Dam, Koolen, and
Xiong~\cite{vDKX} (their Section~5.1) split it into $B_1$ (difference alternating of rank $2$; a square spread, $t_1=1$)
and $B_2$ (difference of rank $3$; $t_2=q^3-q^2$). Over a perfect field of characteristic $2$, a nonzero symmetric form is
determined up to congruence by its rank and by whether it is alternating, so the orbits of the group generated by
translations and congruences $M\mapsto PMP^\top$ on pairs are exactly the relations of $\mathcal C$, and every point set invariant under
this group is concordant with $\mathcal C$.

\begin{corollary}\label{cor:symcodes}
Let $q$ be even. The convex hull of the indicators of the $q^3$-subsets of $V$ whose pairwise differences all have rank~$3$
has dimension ${(q^3-q^2)(q^3-1)}$. The convex hull of the indicators of the $q^3$-subsets whose pairwise differences all
lie in $B$ has dimension ${(q^3-q^2+1)(q^3-1)}$.
\end{corollary}

\begin{proof}
Both families are invariant, hence concordant. Their members are cliques of size $n$ in $B_2$\,, respectively in $B$, which are strongly
regular of Latin-square type. By Delsarte's clique bound, such a graph has no clique with more than $n$ vertices, and the centered
indicator of a clique with exactly $n$ vertices is an eigenvector for the eigenvalue $n-t_2$\,, respectively $n-t$; that is, it lies in
$F_2$\,, respectively in $F=F_1+F_2$\,. Hence $\lambda_i=0$ for every $i$, and $\dim(\mathbf 1^\top X)=0$. Both families are nonempty: for a basis $b_1\,,b_2\,,b_3$ of $\GF(q^3)$ over $\GF(q)$ and $x\in\GF(q^3)$, the Gram matrix $G_x\coloneqq(\operatorname{Tr}(xb_ib_j))_{i,j}$ of the symmetric form $(y,z)\mapsto\operatorname{Tr}(xyz)$ is symmetric, and $G_x-G_y=G_{x-y}$ is nonsingular for $x\neq y$, because the trace form is nondegenerate; so $\{G_x\}$ is a set of $q^3$ matrices whose pairwise differences have rank $3$, and it has a nonzero component in $F_2$\,. A coset of the subspace of
alternating matrices is a clique of size $n$ in $B_1$\,, with a nonzero component in $F_1$\,. So exactly the idempotents $F_2$\,,
respectively $F_1$ and $F_2$\,, are witnessed, of ranks $t_2(n-1)$, respectively $t_2(n-1)+t_1(n-1)$. Apply
Proposition~\ref{prop:fission}.
\end{proof}

\subsection{Symmetric matrices}\label{sec:symm}

\emph{A direct count.} The first dimension of Corollary~\ref{cor:symcodes} can also be obtained, for every $q$, without the fission, from a count of points and lines on conics. For odd $q$ this is the only route, because there are no nonzero alternating symmetric matrices, so the distance-$3$ relation consists of the nonsingular differences and the two families of Corollary~\ref{cor:symcodes} coincide. To cover every characteristic at once we write the characters of the space $V$ of symmetric $m\times m$ matrices as quadratic forms: for a quadratic form $Q(x)=x^\top Ux$ on $\GF(q)^m$ and $M\in V$ put $\langle Q,M\rangle\coloneqq\operatorname{tr}(U^\top M)$. This does not depend on the choice of $U$, which is determined up to an alternating matrix, pairing to $0$ with every symmetric $M$; it is nondegenerate, so $M\mapsto\psi(\langle Q,M\rangle)$ runs over all characters of $V$ as $Q$ runs over the quadratic forms; and $\langle Q\circ P,M\rangle=\langle Q,PMP^\top\rangle$ for invertible $P$, while $\langle Q,cM\rangle=\langle cQ,M\rangle$. For odd $q$ it is the pairing $\operatorname{tr}(WM)$ with $W\coloneqq(U+U^\top)/2$ the matrix of $Q$. A quadratic form $Q$ on $\GF(q)^3$ has rank $3$, $2$, or $1$ according as the conic $Q=0$ in $\PG(2,q)$ is nondegenerate, a pair of lines or a point, or a line; for odd $q$ this is the rank of $W$.

\begin{proposition}\label{prop:symodd}
Let $V$ be the space of symmetric $3\times3$ matrices over $\GF(q)$, let $n\coloneqq q^3$, and let $R$ be the set of nonsingular $M\in V$. Then, for every quadratic form $Q\neq0$ on $\GF(q)^3$,
\[
\sum_{M\in R}\psi(\langle Q,M\rangle)=\begin{cases}q^2,& \operatorname{rank}(Q)=3;\\-(q^3-q^2),& \operatorname{rank}(Q)\le2.\end{cases}
\]
So $R$ is a strongly regular relation of Latin-square type with $t=q^3-q^2$, and the polytope of the $q^3$-subsets of $V$ whose pairwise differences are nonsingular has dimension $(q^3-q^2)(q^3-1)$.
\end{proposition}

\begin{proof}
Let $\pi(Q)$ be the number of points and $\lambda(Q)$ the number of lines on the conic $Q=0$ in $\PG(2,q)$: $(\pi,\lambda)$ is $(q+1,0)$ for a nondegenerate conic, $(2q+1,2)$ for two lines, $(1,0)$ for a point, and $(q+1,1)$ for a line, in every characteristic. Let $F(Q)\coloneqq\sum_{M\notin R}\psi(\langle Q,M\rangle)$, so that the sum in the statement is $-F(Q)$. A singular $M$ has $1$, $q+1$, or $q^2+q+1$ kernel points according as its rank is $2$, $1$, or $0$, so $[M\notin R]=\sum_{\langle v\rangle}[Mv=0]-q[\operatorname{rank}(M)=1]-(q^2+q)[M=0]$, the sum over the points $\langle v\rangle$ of $\PG(2,q)$. The matrices with $Mv=0$ are the $\sum_{a,b}c_{ab}e_ae_b^\top$ with $c$ symmetric, for a basis $e_1\,,e_2$ of $v^\perp\coloneqq\{x:x^\top v=0\}$, a space of dimension $3$, and $\langle Q,\sum c_{ab}e_ae_b^\top\rangle=\sum_ac_{aa}Q(e_a)+c_{12}B_Q(e_1\,,e_2)$ with $B_Q(x,y)\coloneqq Q(x+y)-Q(x)-Q(y)$; so the sum of $\psi(\langle Q,M\rangle)$ over $Mv=0$ is $q^3$ if $Q$ vanishes on $v^\perp$, that is, if the line $v^\top x=0$ lies on the conic, and $0$ otherwise, which contributes $q^3\lambda(Q)$ in total. The matrices of rank $1$ are the $cvv^\top$, $c\neq0$, one for each point $\langle v\rangle$ with a fixed representative $v$ and each $c$; as $\langle Q,cvv^\top\rangle=cQ(v)$, the sum $\sum_{c\neq0}\psi(cQ(v))$ is $q-1$ or $-1$ according as $\langle v\rangle$ lies on the conic or not, so these matrices contribute $q\pi(Q)-(q^2+q+1)$. Hence
\[
F(Q)=q^3\lambda(Q)-q\bigl(q\pi(Q)-(q^2+q+1)\bigr)-(q^2+q)=q^3\lambda(Q)-q^2\pi(Q)+q^3,
\]
which is $-q^2$ for $Q$ of rank $3$ and $q^3-q^2$ for every other $Q\neq0$, as claimed. In particular the eigenvalues of $R$ other than $|R|$ are $n-t=q^2$ and $-t$, so $R$ is of Latin-square type.

Let $X$ be the set of indicators, which is nonempty by the Gram matrices $G_x$ in the proof of Corollary~\ref{cor:symcodes}. Its members are cliques of size $n$ for $R$, so by Delsarte's clique bound their centered indicators lie in the eigenspace for $n-t$, spanned by the characters $M\mapsto\psi(\langle Q,M\rangle)$ with $Q$ of rank $3$; so $\hat x(Q)=0$ for every $x\in X$ and every $Q\neq0$ of rank at most $2$. The set $X$ is invariant under translations, under $M\mapsto PMP^\top$ for invertible $P$, and under $M\mapsto cM$ for $c\neq0$. As in Steps~1--3 of the proof of Theorem~\ref{thm:affine}, with the pairing $\langle Q,M\rangle$ in place of $w\cdot z$, the dimension of $\conv(X)$ is the number of characters $Q$ with $\hat x(Q)\neq0$ for some $x\in X$, because $\mathbf 1^\top x\equiv q^3$, and this set of $Q$ is a union of orbits of the action $Q\mapsto Q\circ P^{-1}$, $Q\mapsto c^{-1}Q$. The $Q$ of rank $3$ form a single orbit: for even $q$ there is one class of nondegenerate ternary quadratic forms, and for odd $q$ the two classes, with discriminant a square or not, are exchanged by scaling with a nonsquare $c$, as $\det(cW)=c^3\det W$. A member of $X$ has a nonzero centered indicator, so some $Q$ of rank $3$, and hence every one, belongs to this set. So the dimension is $|R|=t(n-1)=(q^3-q^2)(q^3-1)$.
\end{proof}

\emph{Four-by-four matrices.} A relation of Latin-square type on the space of symmetric $m\times m$ matrices needs $q^{m(m+1)/2}$ to be a square, that is, $m\equiv0$ or $3$ modulo $4$. After $m=3$ the next case is $m=4$, with $n=q^5$. There the nonsingular differences do not have the right valency, but the differences of rank at least $3$ do: by Schmidt's bound for symmetric rank-metric codes, a set of symmetric $4\times4$ matrices whose pairwise differences have rank at least $3$ has at most $q^5$ elements (see~\cite{Sch}). The count of Proposition~\ref{prop:symodd} extends, with quadrics in $\PG(3,q)$. We say that a quadratic form $Q$ on $\GF(q)^4$ has rank $4$ if the quadric $Q=0$ in $\PG(3,q)$ is nondegenerate, hyperbolic or elliptic, rank $3$ if it is a cone over a conic, rank $2$ if it is a pair of planes or a line, and rank $1$ if it is a plane; for odd $q$ this is the rank of its matrix.

\begin{proposition}\label{prop:sym4}
Let $V$ be the space of symmetric $4\times4$ matrices over $\GF(q)$, let $n\coloneqq q^5$, and let $R$ be the set of $M\in V$ of rank at least $3$. Then, for every quadratic form $Q\neq0$,
\[
\sum_{M\in R}\psi(\langle Q,M\rangle)=\begin{cases}q^2,& \operatorname{rank}(Q)\ge3;\\-(q^5-q^2),& \operatorname{rank}(Q)\le2.\end{cases}
\]
So $R$ is a strongly regular relation of Latin-square type with $t=q^5-q^2$, and the polytope of the $q^5$-subsets of $V$ whose pairwise differences have rank at least $3$ has dimension $(q^5-q^2)(q^5-1)$.
\end{proposition}

\begin{proof}
Let $\pi(Q)$ be the number of points and $\lambda(Q)$ the number of lines on the quadric $Q=0$, and $N\coloneqq q^3+q^2+q+1$ the number of points of $\PG(3,q)$. The pair $(\pi,\lambda)$ is $((q+1)^2,2(q+1))$ for a hyperbolic quadric, $(q^2+1,0)$ for an elliptic one, $(q^2+q+1,q+1)$ for a cone, $(2q^2+q+1,2q^2+2q+1)$ for two planes, $(q+1,1)$ for a line, and $(q^2+q+1,q^2+q+1)$ for a plane, in every characteristic.

\emph{Eigenvalues.} Let $F(Q)\coloneqq\sum_{\operatorname{rank}(M)\le2}\psi(\langle Q,M\rangle)$, so that the sum in the statement is $-F(Q)$. A matrix of rank $2$, $1$, or $0$ has a kernel containing $1$, $q^2+q+1$, or $(q^2+1)(q^2+q+1)$ planes (subspaces of dimension $2$), so $[\operatorname{rank}(M)\le2]=\sum_L[ML=0]-(q^2+q)[\operatorname{rank}(M)=1]-\bigl((q^2+1)(q^2+q+1)-1\bigr)[M=0]$, the sum over the planes $L$ of $\GF(q)^4$. The matrices with $ML=0$ are the $M=\sum_{a,b}c_{ab}e_ae_b^\top$ with $c$ symmetric, for a basis $e_1\,,e_2$ of the $2$-dimensional space $L^\perp\coloneqq\{x:x^\top L=0\}$; they form a space of dimension $3$, and $\langle Q,M\rangle=\sum_ac_{aa}Q(e_a)+c_{12}B_Q(e_1\,,e_2)$ with $B_Q(x,y)\coloneqq Q(x+y)-Q(x)-Q(y)$; so $Q$ annihilates the space exactly when $Q$ vanishes on $L^\perp$, that is, when the line of $L^\perp$ lies on the quadric, and the first term contributes $q^3\lambda(Q)$. The matrices of rank $1$ are the $cvv^\top$, with $\langle Q,cvv^\top\rangle=cQ(v)$, and they contribute $q\pi(Q)-N$ as in Proposition~\ref{prop:symodd}. Hence
\[
F(Q)=q^3\lambda(Q)-q^2(q+1)\pi(Q)+q^5+q^4+q^3,
\]
which is $-q^2$ for the three kinds of $Q$ of rank at least $3$ and $q^5-q^2$ for the three kinds of rank $1$ or $2$.

\emph{A set attaining the bound.} Let $G_x$\,, $x\in\GF(q^5)$, be the Gram matrices of the forms $(y,z)\mapsto\operatorname{Tr}(xyz)$ with respect to a basis $b_1\,,\dots,b_5$ of $\GF(q^5)$ over $\GF(q)$, as in the proof of Corollary~\ref{cor:symcodes}, and let $C$ be the set of their upper left $4\times4$ submatrices, the restrictions to the span of $b_1\,,\dots,b_4$\,. Restricting a form to a hyperplane lowers its rank by at most $2$, and $G_x$ has rank $5$ for $x\neq0$; so $C$ is a linear space of $q^5$ matrices whose nonzero elements have rank at least $3$.

\emph{Dimension.} The set $X$ of indicators is invariant under translations, congruences $M\mapsto PMP^\top$, and scalings $M\mapsto cM$, and as in the proof of Proposition~\ref{prop:symodd}, by Delsarte's clique bound and Steps~1--3 of the proof of Theorem~\ref{thm:affine}, its dimension is the total number of characters $Q$, all of rank at least $3$, in the orbits at which some member of $X$ has a nonzero coefficient. The $Q$ of rank at least $3$ form three orbits, the cones, the hyperbolic, and the elliptic forms, because these are the congruence classes of quadratic forms of rank $3$ and $4$ and the quadric of $cQ$ is that of $Q$. The linear set $C$ has coefficient $|C|$ at the $Q$ in $C^\perp\coloneqq\{Q:\langle Q,M\rangle=0 \text{ for all } M\in C\}$ and $0$ elsewhere (Lemma~\ref{lem:linear}), so it suffices that $C^\perp$ meets all three orbits. Let $A_3$\,, $A_h$\,, and $A_e$ be the numbers of elements of $C^\perp$ in them; by the above, $C^\perp\setminus\{0\}$ has rank at least $3$, so $A_3+A_h+A_e=q^5-1$. For a function $f$ on $V$ that vanishes on $C$, Poisson summation $\sum_{M\in C}f(M)=q^{-5}\sum_{Q\in C^\perp}\hat f(Q)$, with $\hat f(Q)\coloneqq\sum_Mf(M)\psi(\langle Q,M\rangle)$, gives $\sum_{Q\in C^\perp}\hat f(Q)=0$. For $f$ the indicator of the matrices of rank $1$, $\hat f(Q)=q\pi(Q)-N$, and the equation reads $\sum A_k\pi_k=N(q^4-1)$. A third equation separates $A_h$ from $A_e$\,. For odd $q$, take $f$ the indicator of the matrices $c(ab^\top+ba^\top)$ with $a,b$ independent, the rank-$2$ forms that are products of two linear forms; choosing representatives of the points, each is $c(ab^\top+ba^\top)$ for exactly one unordered pair of distinct points and one $c\neq0$, and $\langle Q,c(ab^\top+ba^\top)\rangle=cB_Q(a,b)$, so $\hat f(Q)=qO(Q)-\binom N2$, where $O(Q)$ is the number of unordered pairs of distinct points with $B_Q(a,b)=0$. Counting, for each point $\langle a\rangle$, the points of the plane $B_Q(a,\cdot)=0$ (or of all of $\PG(3,q)$ if that form is zero), and using $B_Q(a,a)=2Q(a)$, we get $O(Q)=\frac12\bigl(N(q^2+q+1)-\pi(Q)\bigr)$ for $Q$ of rank $4$ and $O(Q)=\frac12\bigl(N(q^2+q+1)+q^3-\pi(Q)\bigr)$ for a cone, and the equation reads $\sum A_kO_k=\binom N2(q^4-1)$. For even $q$ the products are alternating, and we take instead $f$ the indicator of the matrices $c_1aa^\top+c_2bb^\top$ with $c_1\,,c_2\neq0$ and $a,b$ independent, the rank-$2$ forms that are not alternating; each arises from exactly $q/2$ unordered pairs of points, because its restriction to its support is a nondegenerate form $B$ with $B(x,x)$ the square of a linear form, so that of the $q+1$ points $q$ are anisotropic and they pair off into orthogonal pairs. As $\langle Q,c_1aa^\top+c_2bb^\top\rangle=c_1Q(a)+c_2Q(b)$, summing over $c_1$ and $c_2$ gives $\hat f(Q)=\frac2q\sum_{\{a,b\}}\bigl(q[Q(a)=0]-1\bigr)\bigl(q[Q(b)=0]-1\bigr)=\frac1q\bigl((q\pi-N)^2-(q-1)^2\pi-(N-\pi)\bigr)$, and $\hat f(0)=\frac2q\binom N2(q-1)^2$. In both cases the three equations give
\[
A_3=q^4-1,\qquad A_h=\tfrac12q^2(q-1)(q^2+1),\qquad A_e=\tfrac12q^2(q-1)^2(q+1),
\]
all positive. So all three orbits are witnessed, and $\dim(\conv(X))$ is $|R|=t(n-1)$, that is, ${(q^5-q^2)(q^5-1)}$.
\end{proof}

The pattern continues at least as far as Schmidt's bound goes: the largest additive set of symmetric $m\times m$ matrices with pairwise differences of rank at least $d$ has exactly $\sqrt{|V|}=q^{m(m+1)/4}$ elements for $d=(m+3)/2$ when $m\equiv3$ and for $d=(m+2)/2$ when $m\equiv0$ modulo $4$ (see~\cite{Sch}). Propositions~\ref{prop:symodd} and~\ref{prop:sym4} are the cases $m=3$ and $m=4$; whether the relation ``rank at least $d$'' is of Latin-square type for all such $m$ we leave open.
\section{Consequences for integer-programming formulations}\label{sec:ip}
We collect what the results of \S\S\ref{sec:affine}--\ref{sec:almost} say about integer-programming formulations.

\paragraph{Affinely tight relaxations.}
The natural integer-programming formulation of the maximum-weight transversal problem of a Latin square $L$ on $[n]$, equivalently of the maximum-weight rainbow perfect matching problem of the edge-colored $K_{n,n}\,$, with weights $c_{ab}$ on the cells, is
\begin{align*}
\max\ &\textstyle\sum_{a,b}c_{ab}x_{ab}\quad\text{s.t.}\quad x\in\{0,1\}^{n\times n},\\
&\textstyle\sum_bx_{ab}=1\ (a\in[n]),\ \ \sum_ax_{ab}=1\ (b\in[n]),\ \ \sum_{L(a,b)=g}x_{ab}=1\ (g\in[n]).
\end{align*}
Theorem~\ref{thm:main} says that for the Cayley table of an abelian group with a transversal the linear relaxation is \emph{affinely tight}: its affine hull is the affine hull of the integer points, and every valid equation is a combination of the $3n$ line sums. The same holds for the transversals of a Desarguesian net by Corollary~\ref{cor:nets}, for the multi-index assignment polytopes $P^{(k,s)}_q$ outside the exceptional families by Theorems~\ref{thm:assign} and~\ref{thm:qplus2}, and at the composite orders of Proposition~\ref{prop:products}. By Proposition~\ref{prop:stability} it also holds for every Latin square whose transversals have nearly uniform pair counts. For facet proofs, which need $\dim(P)$ affinely independent vertices on the face, this is the first piece of information required.

\paragraph{Explicit missing equations.}
In the exceptional families the formulation is not affinely tight, and the spectral method says exactly what is missing. For $(k,s)=(q+1,q-1)$, Lemma~\ref{lem:perfect} gives the equations~\eqref{eq:cuts} below, one for each of the $(q-1)^q$ directions $w$ of full weight and each $c\in\GF(q)$, a system of rank $(q-1)^{q+1}$, which is the whole deficit; for $(k,s)=(q+2,q-1)$ with $q\ge4$ even, Lemma~\ref{lem:puncture} gives the same equations for the directions of weight $q+1$, of rank $(q+2)(q-1)^{q+1}$. These equations hold at every integer point, are not implied by the relaxation, and cost nothing to write down. Explicitly, the $(k,s)$ assignment problem with weights $c_z\,$, $z\in[n]^k$, is
\begin{equation}\label{eq:ip}
\begin{aligned}
\max\ &\textstyle\sum_{z\in[n]^k}c_zx_z\quad\text{s.t.}\quad x\in\{0,1\}^{[n]^k},\\
&\textstyle\sum_{z:\,z_T=t}x_z=1\ \ \bigl(T\subseteq[k],\ |T|=s,\ t\in[n]^T\bigr),
\end{aligned}
\end{equation}
where $z_T$ is the restriction of $z$ to the coordinates in $T$; its linear relaxation replaces $x\in\{0,1\}^{[n]^k}$ by $x\ge0$ (the bounds $x\le1$ are implied), and the strengthened relaxation adds the equations
\begin{equation}\label{eq:cuts}
\sum_{z:\,w\cdot z=c}x_z=q^{q-2}\qquad\bigl(w\in\{1\}\times(\GF(q)^*)^{q},\ c\in\GF(q)\bigr),
\end{equation}
one representative $w$ per full-weight direction. Equivalently, since every array here is a perfect code, $x(B(z,1))=1$ for every $z$, with $B(z,1)$ the Hamming ball of radius $1$; as the Fourier transform of $1_{B(z,1)}$ at $w$ is $q(q-\mathrm{wt}(w))\psi(w\cdot z)$, these ball equations alone cut out the affine hull. For $k=q+2$, the cylinders $\{z:z_{-i}\in B(y,1)\}$ over such balls play the same role. For pairs of orthogonal Latin squares of order $3$, that is $(k,s,n)=(4,2,3)$ (Example~\ref{ex:tetra}), the $81$ variables $x_z\,$, $z\in\GF(3)^4$, are subject to the $54$ equations of~\eqref{eq:ip} ($6$ pairs $T$ of coordinates, $9$ values $t$) of rank $33$, and the $24$ equations~\eqref{eq:cuts} ($8$ directions $w\in\{1\}\times\{1,2\}^3$, $3$ values $c$, right-hand side $3$) have rank $16$ over them, lowering the dimension of the relaxation from $48$ to $32$, the dimension of the polytope. For $300$ random Gaussian objectives $c$, the LP bound of the plain relaxation exceeded the integer optimum (over the $72$ arrays) in $246$ cases, and in every one of these the added equations closed the gap entirely. The strengthened relaxation is not integral, but nearly: with random objectives it returned a fractional vertex in $1$ of $2000$ solves, and the plain relaxation in $390$ of $500$ solves.

\paragraph{Faces of stable-set polytopes.}
All of this can be read as a statement about the face $F_\Gamma\coloneqq\{x\in\mathrm{STAB}(\Gamma):x(E)=\alpha(\Gamma)\}$ of a stable-set polytope on which the rank inequality is tight. For a Latin square $L$, $\Gamma$ has the $n^2$ cells as vertices, two being adjacent when they share a row, a column, or a symbol; $\alpha(\Gamma)\le n$, with equality exactly when $L$ has a transversal, and then $F_\Gamma=\mathcal T(L)$. The $3n$ lines are cliques, and the $n$ row cliques partition $E$, so on $F_\Gamma$ the clique inequalities $x(C)\le1$ of all lines are forced to equality: $F_\Gamma$ has codimension at least $3n-2$, and Theorem~\ref{thm:main} says that for an abelian group table there is nothing more. For $\Gamma$ on $[n]^k$ with $z\sim z'$ when they agree in at least $s$ coordinates, $\alpha(\Gamma)\le n^s$ with equality exactly when an $\mathrm{OA}(n^s,k,n,s)$ exists, $F_\Gamma=P^{(k,s)}_n\,$, the sets $\{z:z_T=t\}$ are cliques partitioning $E$ for each $T$, and the same argument turns their clique inequalities into the equations of~\eqref{eq:ip}; the equations~\eqref{eq:cuts} of the exceptional cases are not combinations of these, but they too are clique equations: the Hamming balls of radius~$1$ are the type~II cliques of~\cite{AMMJ}, which define facets of the orthogonal Latin-square polytope for $n\ge5$, $n\neq6$, and are separated as cuts in~\cite{AMM04}. At the exceptional orders they hold with equality, and the spectral method identifies exactly which ones are needed. For the rank-distance graph on symmetric $m\times m$ matrices, with $M\sim M'$ when $\rk(M-M')\le2$, Schmidt's bound $\alpha(\Gamma)\le n=q^{m(m+1)/4}$ is tight for $m=3,4$, and $F_\Gamma$ has codimension $|E|-(n-q^2)(n-1)$ by Propositions~\ref{prop:symodd} and~\ref{prop:sym4}. In each case $\alpha(\Gamma)$ is an eigenvalue bound of Delsarte or Hoffman type, the rank inequality is far from facet-defining, and the dimension formulae identify the implicit equations of $F_\Gamma\,$, exactly those that a formulation of the maximum-weight largest-stable-set problem needs.

\subsection*{Declaration on the use of generative AI}
During the preparation of this work, the author used Claude to assist in developing and checking all of the work. The author takes full responsibility for the content of the article.

\bibliographystyle{amsplain-doi}
\bibliography{long}

@article{AMM,
  author  = {Appa, Gautam and Magos, Dimitris and Mourtos, Ioannis},
  title   = {On multi-index assignment polytopes},
  journal = {Linear Algebra and its Applications},
  volume  = {416},
  year    = {2006},
  pages   = {224--241},
  doi     = {10.1016/j.laa.2005.11.009}
}

@article{AMMJ,
  author  = {Appa, Gautam and Magos, Dimitris and Mourtos, Ioannis and Janssen, Jeannette C. M.},
  title   = {On the orthogonal {L}atin squares polytope},
  journal = {Discrete Mathematics},
  volume  = {306},
  year    = {2006},
  pages   = {171--187},
  doi     = {10.1016/j.disc.2005.10.020}
}

@article{BS,
  author  = {Balas, Egon and Saltzman, Matthew J.},
  title   = {Facets of the three-index assignment polytope},
  journal = {Discrete Applied Mathematics},
  volume  = {23},
  year    = {1989},
  pages   = {201--229},
  doi     = {10.1016/0166-218X(89)90014-0}
}

@article{BR,
  author  = {Balinski, Michel L. and Russakoff, Andrew},
  title   = {On the assignment polytope},
  journal = {SIAM Review},
  volume  = {16},
  year    = {1974},
  pages   = {516--525},
  doi     = {10.1137/1016083}
}

@article{Bul,
  author  = {Bulutoglu, Dursun A.},
  title   = {Finding the dimension of a non-empty orthogonal array polytope},
  journal = {Discrete Optimization},
  volume  = {45},
  year    = {2022},
  pages   = {100727},
  doi     = {10.1016/j.disopt.2022.100727}
}

@misc{vDKX,
  author = {van Dam, Edwin R. and Koolen, Jack H. and Xiong, Yanzhen},
  title  = {Almost amorphic association schemes},
  note   = {Preprint, arXiv:2604.06538},
  year   = {2026},
  doi    = {10.48550/arXiv.2604.06538}
}

@article{vDM,
  author  = {van Dam, Edwin R. and Muzychuk, Mikhail E.},
  title   = {Some implications on amorphic association schemes},
  journal = {Journal of Combinatorial Theory, Series A},
  volume  = {117},
  year    = {2010},
  pages   = {111--127},
  doi     = {10.1016/j.jcta.2009.03.018}
}

@article{EBG,
  author  = {Euler, Reinhardt and Burkard, Rainer E. and Grommes, Rainer},
  title   = {On {L}atin squares and the facial structure of related polytopes},
  journal = {Discrete Mathematics},
  volume  = {62},
  year    = {1986},
  pages   = {155--181},
  doi     = {10.1016/0012-365X(86)90116-0}
}

@misc{EMM,
  author = {Eberhard, Sean and Manners, Freddie and Mrazovi\'{c}, Rudi},
  title  = {An asymptotic for the {H}all--{P}aige conjecture},
  note   = {Preprint, arXiv:2003.01798},
  year   = {2020},
  doi    = {10.48550/arXiv.2003.01798}
}

@article{EGS,
  author  = {Ezerman, Martianus Frederic and Grassl, Markus and Sol\'{e}, Patrick},
  title   = {The weights in {MDS} codes},
  journal = {IEEE Transactions on Information Theory},
  volume  = {57},
  year    = {2011},
  pages   = {392--396},
  doi     = {10.1109/TIT.2010.2090246}
}

@article{HP,
  author  = {Hall, Marshall and Paige, Lowell J.},
  title   = {Complete mappings of finite groups},
  journal = {Pacific Journal of Mathematics},
  volume  = {5},
  year    = {1955},
  pages   = {541--549},
  doi     = {10.2140/pjm.1955.5.541}
}

@article{IRT,
  author  = {Itai, Alon and Rodeh, Michael and Tanimoto, Steven L.},
  title   = {Some matching problems for bipartite graphs},
  journal = {Journal of the ACM},
  volume  = {25},
  year    = {1978},
  pages   = {517--525},
  doi     = {10.1145/322092.322093}
}

@incollection{Karp,
  author    = {Karp, Richard M.},
  title     = {Reducibility among combinatorial problems},
  booktitle = {Complexity of Computer Computations},
  editor    = {Miller, Raymond E. and Thatcher, James W. and Bohlinger, Jean D.},
  publisher = {Plenum Press},
  address   = {New York},
  year      = {1972},
  pages     = {85--103},
  doi       = {10.1007/978-1-4684-2001-2_9}
}

@techreport{LL90,
  author      = {Leung, Janny and Lee, Jon},
  title       = {Reinforcing old fences gives new facets},
  institution = {Department of Operations Research, Yale University},
  number      = {90-22},
  year        = {1990},
  month       = nov,
  note        = {reprinted within~\cite{LL26}}
}

@misc{LL26,
  author = {Leung, Janny and Lee, Jon},
  title  = {Reinforcing old fences gives new facets and more},
  year   = {2026},
  month  = oct,
  note   = {\url{https://optimization-online.org/?p=37550}}
}

@article{L1,
  author  = {Lee, Jon},
  title   = {A spectral approach to polyhedral dimension},
  journal = {Mathematical Programming},
  volume  = {47},
  year    = {1990},
  pages   = {441--459},
  doi     = {10.1007/BF01580874}
}

@article{L2,
  author  = {Lee, Jon},
  title   = {Characterizations of the dimension for classes of concordant polytopes},
  journal = {Mathematics of Operations Research},
  volume  = {15},
  year    = {1990},
  pages   = {139--154},
  doi     = {10.1287/moor.15.1.139}
}

@article{L3,
  author  = {Lee, Jon},
  title   = {Canonical equation sets for classes of concordant polytopes},
  journal = {Discrete Applied Mathematics},
  volume  = {28},
  year    = {1990},
  pages   = {207--221},
  doi     = {10.1016/0166-218X(90)90003-U}
}

@article{LeeLee,
  author  = {Lee, Jon and Lee, Yew Sing},
  title   = {A spectral method for concordant polyhedral faces},
  journal = {Linear Algebra and its Applications},
  volume  = {159},
  year    = {1991},
  pages   = {11--41},
  doi     = {10.1016/0024-3795(91)90074-7}
}

@article{LP,
  author  = {Le, Van Bang and Pfender, Florian},
  title   = {Complexity results for rainbow matchings},
  journal = {Theoretical Computer Science},
  volume  = {524},
  year    = {2014},
  pages   = {27--33},
  doi     = {10.1016/j.tcs.2013.12.013}
}

@article{Mit,
  author  = {Mittenthal, Lothrop},
  title   = {Block substitutions using orthomorphic mappings},
  journal = {Advances in Applied Mathematics},
  volume  = {16},
  year    = {1995},
  pages   = {59--71},
  doi     = {10.1006/aama.1995.1003}
}

@article{Sch,
  author  = {Schmidt, Kai-Uwe},
  title   = {Symmetric bilinear forms over finite fields with applications to coding theory},
  journal = {Journal of Algebraic Combinatorics},
  volume  = {42},
  year    = {2015},
  pages   = {635--670},
  doi     = {10.1007/s10801-015-0595-0}
}

@incollection{Spi,
  author    = {Spieksma, Frits C. R.},
  title     = {Multi index assignment problems: complexity, approximation, applications},
  booktitle = {Nonlinear Assignment Problems: Algorithms and Applications},
  editor    = {Pardalos, Panos M. and Pitsoulis, Leonidas S.},
  publisher = {Kluwer},
  address   = {Dordrecht},
  year      = {2000},
  pages     = {1--12},
  doi       = {10.1007/978-1-4757-3155-2_1}
}

@book{Wan,
  author    = {Wan, Zhe-Xian},
  title     = {Lectures on Finite Fields and {G}alois Rings},
  publisher = {World Scientific},
  address   = {Singapore},
  year      = {2003},
  doi       = {10.1142/5350}
}

@article{BS91,
  author  = {Balas, Egon and Saltzman, Matthew J.},
  title   = {An algorithm for the three-index assignment problem},
  journal = {Operations Research},
  volume  = {39},
  year    = {1991},
  pages   = {150--161},
  doi     = {10.1287/opre.39.1.150}
}

@article{AMM04,
  author  = {Appa, Gautam and Magos, Dimitris and Mourtos, Ioannis},
  title   = {A branch \& cut algorithm for a four-index assignment problem},
  journal = {Journal of the Operational Research Society},
  volume  = {55},
  year    = {2004},
  pages   = {298--307},
  doi     = {10.1057/palgrave.jors.2601655}
}

@article{MM09,
  author  = {Magos, Dimitris and Mourtos, Ioannis},
  title   = {Clique facets of the axial and planar assignment polytopes},
  journal = {Discrete Optimization},
  volume  = {6},
  year    = {2009},
  pages   = {394--413},
  doi     = {10.1016/j.disopt.2009.05.001}
}

@article{BM08,
  author  = {Bulutoglu, Dursun A. and Margot, Fran\c{c}ois},
  title   = {Classification of orthogonal arrays by integer programming},
  journal = {Journal of Statistical Planning and Inference},
  volume  = {138},
  year    = {2008},
  pages   = {654--666},
  doi     = {10.1016/j.jspi.2006.12.003}
}

\end{document}